\documentclass[11pt,a4paper]{article}

\usepackage[english]{babel}
\usepackage{amsmath}
\usepackage{amssymb}
\usepackage{graphicx}        

\usepackage{cases}
\usepackage{amsfonts}
\usepackage{color}
\usepackage{verbatim}
\usepackage{caption}
\usepackage{subcaption}
\usepackage{fullpage}

\usepackage{array}

\usepackage{epsfig}
\usepackage{float}

\usepackage{algpseudocode}
\usepackage{url}

\usepackage{theorem}

\usepackage{fullpage}

\usepackage{hyperref}

\newcommand{\rosso}[1]{\color{black} #1} 

\usepackage{stackengine}
\newcommand{\tensorm}{
\setstackgap{S}{0.4ex}%
\mathrel{\Shortstack{{.} {.} {.}}}}

\newtheorem{rem}{Remark}[section]
\newtheorem{lem}{Lemma}[section]
\newtheorem{thm}{Theorem}[section]

\newenvironment{pf}{\noindent{\bf Proof. \/}\begin{small}\noindent%
}{\hfill\EndProofMarker\end{small}}
\newcommand{\EndProofMarker}{$\Box$}
\newcommand{\diver}{\ensuremath{\operatorname{div}}}

\def\XXint#1#2#3{{\setbox0=\hbox{$#1{#2#3}{\int}$ }
\vcenter{\hbox{$#2#3$ }}\kern-.6\wd0}}

\begin{document}
  
\title{A Cahn--Hilliard model coupled to viscoelasticity \\ with large deformations}
\author{A. Agosti$^{\sharp}$\footnote{Corresponding author. E-mail address: {\tt abramo.agosti@unipv.it} \newline \textit{Email addresses:} {\tt abramo.agosti@unipv.it} (A. Agosti), {\tt pierluigi.colli@unipv.it} (P. Colli), {\tt harald.garcke@ur.de} (H. Garcke), {\tt elisabetta.rocca@unipv.it} (E. Rocca)}, P. Colli$^{\sharp,\ddag}$, H. Garcke$^{\S}$, E. Rocca$^{\sharp,\ddag}$}

\maketitle 

\begin{center}
{\small $^\sharp$  Department of Mathematics, University of Pavia, 27100 Pavia, Italy.}
\end{center}

\begin{center}
{\small $^\ddag$  Research Associate at the IMATI--C.N.R. Pavia, 27100 Pavia, Italy.}
\end{center}

\begin{center}
{\small $^\S$  Fakult{\"a}t f{\"u}r Mathematik, Universit{\"a}t Regensburg, 93040 Regensburg, Germany.}
\end{center}

\maketitle

\date{}


\begin{abstract}
We propose a new class of phase field models coupled to viscoelasticity with large deformations, obtained from a diffuse interface mixture model composed by a phase with elastic properties and a liquid phase. The model is formulated in the Eulerian configuration and it is derived by imposing the mass balance for the mixture components and the momentum balance {\rosso that comes from} a generalized form of the principle of virtual powers. The latter considers the presence of a system of microforces and microstresses associated to the microscopic interactions
between the mixture's constituents together with a system of macroforces and macrostresses associated to their viscoelastic behavior, taking into account also the friction between the phases. The free energy density of the system is given as the sum of a Cahn--Hilliard term and an elastic polyconvex term, with a coupling between the phase field variable and the elastic deformation gradient in the elastic contribution. General constitutive assumptions complying with a mechanical version of the second law of thermodynamics in isothermal situations are taken. 
We study the global existence of a weak solution for a simplified and regularized version of the general model, which considers an incompressible elastic free energy
of Neo--Hookean type  with elastic coefficients depending on the phase field variable. The regularization is properly designed to deal with the coupling between the phase field variable and the elastic deformation gradient in the elastic energy density. The analysis is made both in two and three space dimensions.  
\end{abstract}



\section{Introduction}
In recent years, mathematical models based on the phase field approach which describe the phase separation phenomena for multicomponent materials in presence of elastic interactions between the materials constituents have been proposed and analyzed, starting from {\rosso the works \cite{collidreyer,miranville,garckeel,garckeel2}}. The dynamics for these materials is characterized, on short time scales, by spinodal decomposition and phase separation, where a microstructure of domains of separated phases is formed driven by the microscopic interactions between the individual components. On large time scales, the domains coarsen and grow driven by their interfacial energy and by the elastic properties of the constituents. These models, which take the form of modified Cahn--Hilliard equations coupled with viscoelastic systems, may find applications in describing the evolution of soft matter materials, e.g. in tumor growth dynamics \cite{garckeel3} and in polymeric soft gel dynamics \cite{lukacova1}.

The major modeling challenge in this framework is to combine the balance equations associated to the dynamics of the phase field variables, which are naturally formulated in the Eulerian configuration, with the constitutive assumptions for elasticity, which are typically formulated in the Lagrangian configuration. 
This issue is typically overcome by formulating the models in the small-strain setting, or by rephrasing the phase field dynamics in the Lagrangian configuration. 
In \cite{gurtin1}, a unified framework for the derivation of Cahn--Hilliard type equations coupled with elasticity is developed, based on standard force and momentum balances for elastic deformations and on balance laws for mass and microforces, given constitutive assumptions consistent with a mechanical version of the second law of thermodynamics.
The balance laws and the constitutive assumptions in \cite{gurtin1} are formulated in a fixed reference configuration, assuming an infinitesimal elastic deformation and neglecting inertia effects. A similar approach is {\rosso carried out in \cite{collidreyer,garckeel,garckeel2,garckeel3,heinemann,krejci}}, where a Cahn--Hilliard equation coupled with infinitesimal elasticity is introduced and analyzed. In these latter studies, the models are derived from a multiphasic mixture. In all the aforementioned works, the elasticity tensor {\rosso and the elastic coefficients} depend on the phase field variables. This coupling between the phase field variables and the elastic deformations is an important feature of the model, which eventually may describe the different elastic properties of each constituent of the mixtures.

Another possibility is to formulate both the phase field and the elasticity governing equations in the Eulerian reference configuration, employing the change of variables from the Lagrangian to the Eulerian configuration in the finite elasticity constitutive assumptions. In this framework, the main state variable for elasticity is the velocity field, and the deformation gradient, entering in the elasticity constitutive assumptions, is determined solving a transport equation in terms of the velocity and the velocity gradient. This modeling approach is employed e.g. in \cite{benesova}, which studies a problem in magnetoelasticity, and in \cite{lukacova1}, which studies a Cahn--Hilliard model coupled with viscoelasticity for a non--Newtonian mixture with a polymeric fluid component. In both the aforementioned works, the coupling between the phase field variable and the elastic deformation gradient in the elasticity constitutive assumptions is not taken into account. 

As for the analysis of the evolutionary problems defined by phase field models coupled with elasticity, we list the following results obtained in the aforementioned literature. In {\rosso \cite{collidreyer} existence and uniqueness results for solutions of a Cahn--Hilliard problem coupled with linear elasticity are obtained, while in \cite{garckeel3} existence, uniqueness and regularity results are proved for solutions of a tumor growth model consisting of a Cahn--Hilliard model coupled with linear elasticity and growth}. For what concerns the problem with nonlinear finite elasticity, in \cite{benesova} the existence of a global in time weak solution to a model of magnetoelasticity is obtained in two space dimensions, while strong well-posedness and stability results are obtained in two space dimensions for a similar problem in \cite{garcke4}. In \cite{lukacova1}, existence of a global in time weak solution to a Cahn--Hilliard model coupled with viscoelasticity is obtained in two space dimensions.  {\rosso Very recently, in \cite{garcke5} existence of a global in time weak solution to a Cahn--Hilliard model coupled with viscoelasticity and growth, extending the model considered in \cite{lukacova1}, is obtained in two space dimensions.}
We observe that {\rosso in all the mentioned works \cite{benesova}, \cite{lukacova1} and \cite{garcke5} an elliptic regularization in the transport equation for the deformation gradient is introduced, and the finite elastic energy density is considered to be of Neo--Hookean type, i.e. convex, and independent on the phase field variable. We further observe that analytical results for this kind of phase-field problems are obtained in the cited literature only in two space dimensions.}

{\rosso For what concerns the analysis of evolutionary viscoelastic problems alone, without the presence of a phase field dynamics, stronger results are available in literature. Regarding models with finite viscoelasticity formulated in the Eulerian configuration, we cite \cite{lukacova2}, where the existence and regularity of a global in time weak solution in three space dimension for an evolutionary viscoelastic problem with finite incompressible Neo--Hookean elasticity is proved. We also cite the recent work in \cite{rubicek}, where existence and regularity results of a distributional solution of an evolutionary non-isothermal viscoelastic problem with polyconvex and compressible elastic energy have been obtained, employing a regularization approach in which a dissipative contribution to the Cauchy stress tensor is introduced which involves high order nonlinear terms in the small strain rate (i.e. the symmetric part of the velocity gradient). We highlight the fact that the latter regularization approach, differently from the elliptic regularization in the transport equation for the deformation gradient introduced in \cite{benesova} and \cite{lukacova1}, preserve the kinematical relationship between the velocity and the deformation gradient variables. Its main drawback is the introduction of nonlinear terms involving the velocity gradient in the definition of the Cauchy stress tensor, making numerical approximations of the model not straightforward. Also, as we will see in the following, the regularization approach based on the elliptic regularization of the transport equation for the deformation gradient may be easily adapted to enhance the regularity in space of the combined product of a function of a phase field variable and a function of the deformation gradient involved in elasticity constitutive assumptions.
Regarding models with finite elasticity formulated in the Lagrangian configuration, we cite \cite{bonetti1}, where a local in time existence result for an evolutionary compressible large deformation problem with Neo--Hookean elasticity and third order gradient regularization in two space dimensions is obtained, and \cite{bonetti2}, where an analogous result is obtained for an evolutionary incompressible large deformation problem formulated in the Lagrangian configuration in three space dimensions. We also cite \cite{mielke}, where existence of a weak solution is proved for an evolutionary non-isothermal viscoelastic problem similar to the one studied in \cite{rubicek} with second order gradient regularization. For what concerns non-isothermal linear viscoelasticity, we cite e.g. \cite{grs}, where the asymptotic behavior of a system describing a thermoviscoelastic plate with hereditary effects is studied.}

In the present paper, we derive a Cahn--Hilliard model coupled to viscoelasticity with finite deformations, starting from a binary mixture dynamics, obtaining the mass and momentum balance of the mixture using a generalized form of the principles of virtual powers, giving constitutive assumptions satisfying the first and second law of thermodynamics in isothermal situations.  The model equations are formulated in the Eulerian configuration, taking both finite elasticity and inertia effects into accounts. Moreover, we further consider in the model derivation the friction between the mixture constituents. We thus propose a new general class of phase field models, expressed in terms of the dynamics of the phase field variables of the two components and in terms of their associated velocity fields, considering general compressible polyconvex elastic energy densities in which the phase field variables and the elastic deformation are inherently coupled. Then, we study the existence of a global in time weak solution for a simplified version of the general model, considering a Neo--Hookean incompressible law for elasticity with an elastic modulus depending on the phase field variable. We highlight that the main difficulties in the analysis of the model, especially in three space dimensions, arise when dealing with the nonlinear coupling betwen the phase field variable and the elastic deformation gradient in the elastic energy density. The {\rosso resulting PDE system for the latter simplified version of the general model} is the following: 
\begin{equation}
\label{eqn:systemintro}
\begin{cases}
-\nu\Delta \mathbf{v}+\nabla q=\mu \nabla \phi- \left(\frac{c}{2}f'(\phi)(\mathbf{F}\colon \mathbf{F}-d)\right)\nabla \phi+\diver\left(f(\phi)\mathbf{F}\mathbf{F}^T\right),\\
\diver\mathbf{v}=0,\\
\frac{\partial \mathbf{F}}{\partial t}+\left(\mathbf{v}\cdot \nabla\right)\mathbf{F}-(\nabla \mathbf{v})\mathbf{F}-\lambda\Delta \left(f(\phi)\mathbf{F}\right)=\boldsymbol{0},\\
\frac{\partial \phi}{\partial t}+\mathbf{v}\cdot \nabla \phi-\diver(b(\phi)\nabla \mu)=0,\\
\mu=\psi'(\phi)- \Delta \phi+\frac{c}{2}f'(\phi)(\mathbf{F}\colon \mathbf{F}-d),\\
\end{cases}
\end{equation}
valid in $\Omega \times [0,T]$, endowed with the boundary conditions
\begin{equation}
\label{eqn:systembc}
b(\phi)\nabla \mu\cdot \mathbf{n}=\nabla \phi \cdot \mathbf{n}=0, \quad \mathbf{v}=\mathbf{0}, \quad \left[\nabla (f(\phi)\mathbf{F})\right]\mathbf{n}=\boldsymbol{0}
\end{equation}
on $\partial \Omega \times [0,T]$, and with proper initial conditions. Here, $\mathbf{v}$ is the velocity field, $q$ is the pressure, $\mathbf{F}$ is the elastic deformation gradient, $\phi$ is the phase field variable and $\mu$ is the chemical potential. {\rosso Moreover, $\nu$ and $c$ are physical parameters, representing the viscosity and the Young elastic modulus of the viscoelastic phase in the mixture, while $\lambda$ is a positive regularization parameter. The function $b(\phi)$ represents a positive phase-dependent mobility, $\psi(\phi)$ represents the bulk interaction energy between the viscoelastic phase microscopic constituents, and $f(\phi)$ describes the dependence of the macroscopic elastic properties of the elastic phase on its concentration.}
In order to obtain analytic results, following ideas proposed in \cite{benesova}, we {\rosso design and} include a proper elliptic regularization in equation \eqref{eqn:systemintro}$_3$. We obtain existence results both in two and three space dimensions. In the latter case, we need to add a further viscous regularization {\rosso in the right hand side of \eqref{eqn:systemintro}$_5$} to obtain the analytic results. 

Hence, the novelties of the present work with respect to previous studies in literature are the following:
\begin{itemize}
\item Derivation of a {\rosso new class of} phase field models of Cahn--Hilliard type coupled with viscoelasticity, in Eulerian coordinates, which consider the microscopic coupling between the mixture constituents, their viscoelastic behaviour and the friction between the phases, with a general compressible polyconvex elastic energy which depends on both the elastic deformation and the phase field variable within a nonlinear coupling; 
\item {\rosso Proof} of existence results for a simplified and regularized incompressible version of the model, where the regularization substantially depends on the coupling between the phase field variable and the elastic deformation gradient in the elastic energy density, for both two and three space dimensions.
\end{itemize}
{\rosso Possible applications of the proposed models are the investigation of their predictive and diagnostic abilities for diseases in which the anatomical tissue viscosity and stiffness play a crucial role in determining the diseases' evolution. This is the case e.g. for neurological diseases \cite{neuro3,neuro2,neuro1}, where the brain tissues stiffness has been proved to be a robust biomarker for neurodegeneration; for tumor growth models \cite{neuro4}, where the mechanical feedback of the surrounding healthy tissues is relevant to predict the disease progression and the consequent tissue reorganization; for neuromuscular diseases \cite{neuro5}, where fat infiltrates the healthy tissues of muscles and the correct assessment of the amount of fat substitution is crucial to evaluate the disease progression and the outcome of possible therapeutic strategies. We highlight the fact that the proposed models may be quantitatively informed by medical imaging data, e.g. by Elastography Magnetic Resonance imaging \cite{neuro1}, which may give as a processing output the values of viscosity and shear elastic modulus of anatomical tissues.
}

The paper is organized as follows. In Section $2$ we introduce some notation regarding the employed tensor calculus and functional analysis. In Section $3$ we develop the model derivation for different class of general models. In Section $4$ we report the study of existence for a global weak solution for a simplified version of the general model. Finally, Section $5$ collects some considerations and conclusions. 

\section{Notation}
Let $\Omega \subset \mathbb{R}^d$ be an open bounded domain in $\mathbb{R}^d$, where $d$ is the space dimension. {\rosso Let $T>0$ denote some final time, and set $\Omega_T:=\Omega \times (0,T)$.}
We start by introducing the notation for vectorial and tensorial calculus. Given $\mathbf{a},\mathbf{b}\in \mathbb{R}^d$, we denote by $\mathbf{a}\cdot \mathbf{b} \in \mathbb{R}$ their canonical scalar product in $\mathbb{R}^d$, with associated norm $|\mathbf{a}|:=(\mathbf{a}\cdot \mathbf{a})^{\frac{1}{2}}$, and by $\mathbf{a}\otimes \mathbf{b}\in \mathbb{R}^{d\times d}$ their tensorial product. Given two second order tensors $\mathbf{A},\mathbf{B}\in \mathbb{R}^{d\times d}$, we denote by $\mathbf{A}\colon \mathbf{B}\in \mathbb{R}$ their Frobenius scalar product in $\mathbb{R}^{d\times d}$, i.e. by components $\mathbf{A}\colon \mathbf{B}:=\sum_{i,j=1}^dA_{ij}B_{ij}$, with associated norm $|\mathbf{A}|:=(\mathbf{A}\colon \mathbf{A})^{\frac{1}{2}}$. Given two third order tensors $\mathbf{C},\mathbf{D}\in \mathbb{R}^{d\times d \times d}$, we denote by $\mathbf{C}\mathbin{\tensorm}\mathbf{D}\in \mathbb{R}$ their scalar product in $\mathbb{R}^{d\times d \times d}$, i.e. by components 
\[\mathbf{C}\mathbin{\tensorm} \mathbf{D}:=\sum_{i,j,k=1}^dC_{ijk}D_{ijk}.\]
We also introduce the operation, defined by components, \[\left(\mathbf{C}\odot \mathbf{D}\right)_{kl}:=\sum_{ij}\mathbf{C}_{kij}\mathbf{D}_{lij},\] 
which contracts two third order tensors $\mathbf{C},\mathbf{D}\in \mathbb{R}^{d\times d \times d}$ to a second order tensor $\mathbf{C}\odot \mathbf{D}\in \mathbb{R}^{d\times d}$. 

We denote by $L^p(\Omega;K)$ and $W^{r,p}(\Omega;K)$ the standard Lebesgue and Sobolev spaces of functions defined on $\Omega$ with values in a set $K$, {\rosso where $K$ may be $\mathbb{R}$ or a multiple power of $\mathbb{R}$,} and by  $L^p(0,t;V)$ the Bochner space of functions defined on $(0,t)$ with values in the functional space $V$, with $1\leq p \leq \infty$, $r>0$. If $K\equiv \mathbb{R}$, we simply write  $L^p(\Omega)$ and $W^{r,p}(\Omega)$. 
Moreover, when stating general results which are valid for both functions with scalar or vectorial or tensorial values, we write $f\in L^p$, $f\in W^{r,p}$, without specifying if $f$ is a function with scalar, vectorial or tensorial values.
For a normed space $X$, the associated norm is denoted by $||\cdot||_X$. In the case $p=2$, we use the notations $H^1:=W^{1,2}$ and $H^2:=W^{2,2}$, and we denote by $(\cdot,\cdot)$ and $||\cdot||$ the $L^2$ scalar product and induced norm between functions with scalar, vectorial or tensorial values. The dual space of a Banach space $Y$ is denoted by $Y'$. The duality pairing between $H^1(\Omega;K)$ and $\left(H^1(\Omega;K)\right)'$ is denoted by $<\cdot,\cdot>$, while  the duality pairing between $H_N^2(\Omega;K)$ and $\left(H_N^2(\Omega;K)\right)'$ is denoted by $<\cdot,\cdot>_2$, where $H_N^2(\Omega;K):=\{f\in H_N^2(\Omega;K): \, \partial_{\mathbf{n}}f=0 \; \text{on} \; \partial \Omega\}$. Moreover, we denote by {\rosso $C^k([0,t];V)$, $k\geq 0$, and by $C_w^0([0,t];V)$} the space of {\rosso continuously differentiable functions up to} order $k$ and of weakly continuous functions, respectively, from $[0,t]$ to the space $V$. When $X\equiv \mathbb{R}$, we write $C^k([0,t])$. We finally introduce the spaces
\begin{align*}
& L_{\text{div}}^2(\Omega,\mathbb{R}^d):=\overline{\{\mathbf{u}\in C_c^{\infty}(\Omega,\mathbb{R}^d): \, \text{div}\mathbf{u}=0 \; \text{in} \; \Omega\}}^{||\cdot||},\\
& H_{0,\text{div}}^1(\Omega,\mathbb{R}^d):=\overline{\{\mathbf{u}\in C_c^{\infty}(\Omega,\mathbb{R}^d): \, \text{div}\mathbf{u}=0 \; \text{in} \; \Omega\}}^{||\cdot||_{H_0^1(\Omega;\mathbb{R}^d)}}.
\end{align*}
\noindent
In the following, $C$ denotes a generic positive constant independent of the unknown variables, the discretization and the regularization parameters, the value of which might change from line to line; $C_1, C_2, \dots$ indicate generic positive constants whose particular value must be tracked through the calculations; $C(a,b,\dots)$ denotes a constant depending on the nonnegative parameters $a,b,\dots$.

\section{Model derivation}
We consider a binary, saturated, closed and non reactive mixture in an open bounded domain $\Omega\subset \mathbb{R}^d$, $d=2,3$, fixed in time, composed by a  phase $c$ with volume fraction $\phi_c$ and a phase $l$ with volume fraction $\phi_l$. In applications, the phase $c$ may represent a solid or a fluid component (with elastic properties), while the phase $l$ may represent a liquid component of the mixture. We assume that both the phases have constant density $\gamma$. In life sciences applications, with the phase $c$ representing a soft biological tissue and the phase $l$ representing water and interstitial fluid, $\gamma$ should be equal to the water density, since biological cells are mostly composed by water.  We assume that the mixture dynamics is determined by the microscopic interactions between its constituents as well as by their macroscopic visco--elastic behaviour. In \cite{gurtin1}, a system of scalar microforces and vectorial microstresses transmitted across surfaces of uniform states, (i.e. level sets $\phi_c(\mathbf{x},t)=\text{constant}$), was introduced, which performs work against changes in the phase variables identified as movements of surfaces of uniform states along their normal. This means that the power expended by microforces and microstresses, computed using the velocity of surfaces of uniform state, is independent of the parameterizations of the state surface, and consequently is independent of the component of the velocity tangential to the state surfaces. Since the normal velocity of the isosurfaces $\phi_c(\mathbf{x},t)=\text{constant}$ is given by $\frac{\partial_t \phi_c}{|\nabla \phi_c|}$, the proper kinematic variables power--conjugated to the microforces and microstresses in this theoretical framework are $\partial_t \phi_c$ and $\nabla (\partial_t \phi_c)$ respectively.

The balance laws and the constitutive assumptions in \cite{gurtin1} were formulated in a fixed reference configuration, assuming an infinitesimal elastic deformation and neglecting inertia effects. In the following, we will consider inertia effects and assume that the elastic deformation is finite. We formulate the mixture dynamics in the Eulerian reference configuration, employing the change of variables from the Lagrangian to the Eulerian configuration in the finite elasticity constitutive assumptions.

Following \cite[Sections $5$ and $6$]{gurtin2} (see also \cite{dileo}), we derive the mass and momentum balance of the mixture using a generalized form of the principles of virtual powers, giving constitutive assumptions satisfying the first and second law of thermodynamics in isothermal situations. 
Hence, each phase satisfies a mass balance equation,
\begin{align}
\label{eqn:1a}
\gamma\frac{\partial \phi_c}{\partial t}+\gamma \diver(\phi_c \mathbf{v}_c)+\diver(\mathbf{h_c})=0,\\
\label{eqn:1b}
\gamma\frac{\partial \phi_l}{\partial t}+\gamma\diver(\phi_l \mathbf{v}_l)+\diver(\mathbf{h_l})=0,
\end{align}
where, as in \cite{gurtin1}, $\mathbf{h_c}$ and $\mathbf{h_l}$ are proper mass fluxes associated with the microsystem dynamics, whose contributions to the mass balance for each phase are added to the corresponding transport contributions determined by the macrosystem dynamics.
Since the mixture is closed, its mass $\rho=\gamma (\phi_c+\phi_l)$ is conserved, i.e. 
\[\frac{\partial \rho}{\partial t}+\diver(\rho \mathbf{v})=0,\]
where $\mathbf{v}=\frac{\gamma}{\rho}(\phi_c \mathbf{v}_c+\phi_l \mathbf{v}_l)$ is the volume-averaged mixture velocity. Hence, adding \eqref{eqn:1a} and \eqref{eqn:1b} gives that $\diver(\mathbf{h_c}+\mathbf{h_l})=0$.
Moreover, since the mixture is saturated, i.e. $\phi_c+\phi_l=1$, we have the incompressibility condition 
\begin{equation}
\label{eqn:2}
\diver\mathbf{v}=0.
\end{equation}
In the following, we will take $\gamma=1$ for simplicity.
\newline
We now introduce the power--conjugate variables:
\begin{align*}
 \mathbf{T}_c \longleftrightarrow \nabla \mathbf{v}_c, \quad &\text{Cauchy stress tensor of the phase $c$},\\
 \mathbf{T}_l \longleftrightarrow \nabla \mathbf{v}_l, \quad &\text{Cauchy stress tensor of the phase $l$},\\
 \mathbf{m} \longleftrightarrow (\mathbf{v}_c-\mathbf{v}_l), \quad &\text{internal force associated to the slip between the two phases},\\
 \pi \longleftrightarrow \frac{\partial \phi_c}{\partial t}, \quad &\text{internal microforce},\\
 \boldsymbol{\xi} \longleftrightarrow \nabla \frac{\partial \phi_c}{\partial t}, \quad &\text{microstress},
\end{align*}
where, as previously explained, the microsystem power is computed using the velocity of surfaces of uniform states along their normal.

Considering an arbitrary part of the mixture $R(t)\subset \Omega$, moving with the mixture, the principle of virtual powers states that the virtual internal power expended within $R(t)$ is equal to the virtual external powers expended in $R(t)$ by material external to $R(t)$ or by external forces. Moreover, as a consequence of frame--indifference, the internal power expended within $R(t)$ for rigid virtual velocities must be equal to zero.\\
We consider the fields $\mathbf{v}_c$, $\mathbf{v}_l$ and $\partial_t \phi_c$ as virtual velocities, to be specified independently. Denoting the virtual fields as $\hat{\mathbf{v}}_c$, $\hat{\mathbf{v}}_l$ and $\hat{\partial_t \phi}_c$, we define a generalized virtual velocity to be a list $\mathcal{V}=\{\hat{\mathbf{v}}_c,\hat{\mathbf{v}}_l,\hat{\partial_t \phi}_c\}$, with $\hat{\mathbf{v}}_c=\mathbf{0}$ on a portion of the domain boundary $\partial \Omega_{\text{D,c}}\subset \partial \Omega$, $\hat{\mathbf{v}}_l=\mathbf{0}$ on $\partial \Omega_{\text{D,l}}\subset \partial \Omega$. We also define a rigid macroscopic virtual velocity to be a list $\mathcal{V}_{\text{rigid}}=\{\hat{\mathbf{v}}_c=\mathbf{v}_{0c}+\mathbf{A}_{0c}\mathbf{x},\hat{\mathbf{v}}_l=\mathbf{v}_{0l}+\mathbf{A}_{0l}\mathbf{x}\}$, with $\mathbf{v}_{0c},\mathbf{v}_{0l}$ constant vectors and $\mathbf{A}_{0c},\mathbf{A}_{0l}$ constant antisymmetric tensors. Moreover, we define the internal and external virtual powers
\begin{align*}
&W_{\text{int}}(R(t),\mathcal{V})=\int_{R(t)}\left( \mathbf{T}_c\colon \nabla \hat{\mathbf{v}}_c+ \mathbf{T}_l\colon \nabla \hat{\mathbf{v}}_l+\mathbf{m}\cdot (\hat{\mathbf{v}}_c-\hat{\mathbf{v}}_l)+\pi \hat{\partial_t \phi}_c +\xi \cdot \nabla \hat{\partial_t \phi}_c\right),\\
& W_{\text{ext}}(R(t),\mathcal{V})=\int_{\partial R(t)}\left( \mathbf{t}_c\cdot \hat{\mathbf{v}}_c+ \mathbf{t}_l\cdot \hat{\mathbf{v}}_l+\eta \hat{\partial_t \phi}_c \right)-\frac{1}{2}\int_{R(t)}K(\phi_c,\phi_l)|\hat{\mathbf{v}}_c-\hat{\mathbf{v}}_l|^2,
\end{align*}
where $\mathbf{t}_c,\mathbf{t}_l$ are tractions on the boundary of $R(t)$, $\eta$ is a scalar microscopic traction and the last term represents friction between the phases exerted on $R(t)$ by material external to it. 
The function $K(\phi_c,\phi_l)$ is a positive coefficient representing the effects of the mixture viscosity and permeability, whose form must be constitutively assigned. Then, the principle of virtual powers states that
\begin{align}
\label{eqn:virtpow}
&W_{\text{int}}(R(t),\mathcal{V})=W_{\text{ext}}(R(t),\mathcal{V}), \quad \forall \, \mathcal{V}, \; \forall R(t) \subset \Omega,\\
\label{eqn:virtpow2}
& W_{\text{int}}(R(t),\mathcal{V})=0 \quad \forall \, \mathcal{V} \equiv \mathcal{V}_{\text{rigid}}, \; \forall R(t) \subset \Omega.
\end{align}
Condition \eqref{eqn:virtpow} gives the momentum balances:
\begin{equation}
\label{eqn:tc}
\begin{cases}
\diver\mathbf{T}_c-K(\phi_c)(\mathbf{v}_c-\mathbf{v}_l)-\mathbf{m}=\boldsymbol{0} \; \; \text{in} \; \Omega,\\
\mathbf{T}_c\mathbf{n}=\mathbf{t}_c \; \; \text{on} \; \partial\Omega_{N,c},\\
\mathbf{v}_c=\mathbf{0} \; \; \text{on} \; \partial\Omega_{D,c},
\end{cases}
\end{equation}
with $\partial\Omega_{D,c} \cup \partial\Omega_{N,c}=\partial\Omega$,
\begin{equation}
\label{eqn:tl}
\begin{cases}
\diver\mathbf{T}_l+K(\phi_c)(\mathbf{v}_c-\mathbf{v}_l)+\mathbf{m}=\boldsymbol{0} \; \; \text{in} \; \Omega,\\
\mathbf{T}_l\mathbf{n}=\mathbf{t}_l \; \; \text{on} \; \partial\Omega_{N,l},\\
\mathbf{v}_l=\mathbf{0} \; \; \text{on} \; \partial\Omega_{D,l},
\end{cases}
\end{equation}
with $\partial\Omega_{D,l} \cup \partial\Omega_{N,l}=\partial\Omega$,
\begin{equation}
\label{eqn:xi}
\begin{cases}
\diver \boldsymbol{\xi}-\pi=0 \; \; \text{in} \; \Omega,\\
\boldsymbol{\xi} \cdot \mathbf{n}=\eta \; \; \text{on} \; \partial\Omega,
\end{cases}
\end{equation}
where $\mathbf{n}$ is the outward normal to $\Omega$.
The condition in \eqref{eqn:virtpow2} gives the relations $\mathbf{T}_c\colon \mathbf{A}_c=0$ and $\mathbf{T}_l\colon \mathbf{A}_l=0$ for any antisymmetric tensors $\mathbf{A}_c,\mathbf{A}_l$, i.e. $\mathbf{T}_c$ and $\mathbf{T}_l$ are symmetric tensors.
 
In order to proceed, we need to introduce the free energy $E$ of the system in $R(t)$. We assume that elasticity acts only on the phase $c$, with a negligible contribution of the phase $l$ to the elastic energy. We thus postulate the following form for $E$:
\begin{align}
\label{eqn:free}
&E(\phi_c,\nabla \phi_c,\mathbf{F}_c,\nabla \mathbf{F}_c)=\int_{R(t)}e(\phi_c,\nabla \phi_c,\mathbf{F}_c,\nabla \mathbf{F}_c)\\
& \notag =\int_{R(t)}\left(\frac{1}{\det\mathbf{F}_c}w(\phi_c,\nabla \phi_c, \mathbf{F}_c, \nabla \mathbf{F}_c)+\Pi \Phi(\phi_c,\nabla \phi_c)\right),
\end{align}
where
\begin{equation}
\label{eqn:Phidef}
\Phi(\phi_c,\nabla \phi_c)=\frac{\epsilon}{2}|\nabla \phi_c|^2+\frac{1}{\epsilon}\psi(\phi_c).
\end{equation}
Here, $\mathbf{F}_c$ is the deformation gradient associated to the motion of the phase $c$, $\frac{1}{\det\mathbf{F}_c}w$ is the hyperelastic free energy density, with units of $[Pa]$, in Eulerian coordinates, and $\Pi$ is proportional to the surface tension, in units of $[N/m]$. The surface energy of the interface between the phases is expressed through a diffuse interface approach, with $\psi(\phi_c)$ representing a bulk energy due to the mechanical interactions of the micro--components. The parameter $\epsilon$ has units of $[m]$ and is related to the interface thickness. The dependence of $w$ on $\nabla \mathbf{F}_c$ is associated to elastic energy contributions from interfaces in the elastic material, while the dependence on $\nabla \phi_c$ is associated to the energy contributions from elastic deformations of material surfaces on the isosurfaces $\phi_c=\text{constant}$. In the following, we will take $\Pi=1$ for simplicity. 

The first and second law of thermodynamics in isothermal situations take the form of the following dissipation inequality \cite{gurtin1}:
\begin{align}
\label{eqn:dissipation}
&\frac{d}{dt}\int_{R(t)}e(\phi_c,\nabla \phi_c,\mathbf{F}_c,\nabla \mathbf{F}_c)+\int_{\partial R(t)}\mathbf{J}_e\cdot \mathbf{n}\\
& \notag \leq W_{\text{ext}}(R(t))-\int_{\partial R(t)}\mu_c(\phi_c\mathbf{v}_c+\mathbf{h}_c)\cdot \mathbf{n}+\int_{R(t)}p\diver(\phi_c\mathbf{v}_c+\phi_l\mathbf{v}_l),
\end{align}
which means that the free energy of $R(t)$ increases at a rate not greater than the working of all forces external to $R(t)$ plus the free energy carried into $R(t)$ by
mass transport. We note that we have introduced in \eqref{eqn:dissipation} the constraint \eqref{eqn:2} as a mass contribution for the mixture through a Lagrange multiplier $p$. Here, $\mu_c$ is the chemical potential of the phase $c$, and $p$ is a Lagrange multiplier for the incompressibility condition for the mixture \eqref{eqn:2}. Using \eqref{eqn:virtpow} and \eqref{eqn:1a} in \eqref{eqn:dissipation} and the following relations for the material derivative of the deformation gradient, for the gradient of the deformation gradient (which are valid if the deformation gradient is sufficiently regular), and for $\partial_{\mathbf{F}_c}\det\mathbf{F}_c$: 
\begin{equation}
\label{eqn:fdot}
\dot{\mathbf{F}}_c=(\nabla \mathbf{v}_c)\mathbf{F}_c,
\end{equation}
\begin{equation}
\label{eqn:fradfdot}
\dot{(\nabla \mathbf{F}_c)}=\nabla \dot{\mathbf{F}}_c - \bigl[\nabla \mathbf{F}_c\bigr]\nabla \mathbf{v}_c,
\end{equation}
\begin{equation}
\label{eqn:detdot}
\partial_{\mathbf{F}_c}\det\mathbf{F}_c=\det\mathbf{F}_c\mathbf{F}_c^{-T}=:\text{cof}\mathbf{F}_c,
\end{equation}
where $\dot{\mathbf{F}}_c$ and $\dot{(\nabla \mathbf{F}_c)}$ indicate  the material time derivatives of $\mathbf{F}$ and $\nabla \mathbf{F}$ respectively, we get, with the help of the Reynolds transport theorem,
\begin{align*}
&\int_{R(t)}\left((\partial_{\phi_c}e)\frac{\partial \phi_c}{\partial t}+(\partial_{\nabla \phi_c}e)\cdot \frac{\partial  \nabla \phi_c}{\partial t}+(\partial_{\phi_c}e)\mathbf{v}_c\cdot \nabla \phi_c+(\partial_{\nabla \phi_c}e)\cdot \nabla (\nabla \phi_c)\mathbf{v}_c +e\mathbf{I}\colon \nabla \mathbf{v}_c\right)\\
&+\int_{R(t)}\left(\frac{1}{\det\mathbf{F}_c}(\partial_{\mathbf{F}_c}w)\mathbf{F}_c^T -\frac{w\mathbf{I}}{\det\mathbf{F}_c}\right)\colon \nabla \mathbf{v}_c+ \int_{R(t)}\biggl\{\biggr(-\diver\left(\frac{1}{\det\mathbf{F}_c}\partial_{\nabla \mathbf{F}_c}w\right) \mathbf{F}_c^T\\
&-\frac{1}{\det\mathbf{F}_c}\nabla \mathbf{F}_c \odot \partial_{\nabla \mathbf{F}_c}w\biggr)\colon \nabla \mathbf{v}_c \biggr\}+\int_{\partial R(t)}\left(\mathbf{J}_e+\frac{1}{\det\mathbf{F}_c}[\partial_{\nabla \mathbf{F}_c}w]^T\dot{\mathbf{F}}_c\right)\cdot \mathbf{n} \\
& \leq \int_{R(t)}\biggl( \mathbf{T}_c\colon \nabla \mathbf{v}_c+ \mathbf{T}_l\colon \nabla \mathbf{v}_l+\mathbf{m}\cdot (\mathbf{v}_c-\mathbf{v}_l)+(\pi+\mu_c) \frac{\partial \phi_c}{\partial t}+\boldsymbol{\xi} \cdot \nabla \frac{\partial \phi_c}{\partial t}-\phi_c\mathbf{v}_c\cdot \nabla \mu_c \\
& -\mathbf{h}_c \cdot \nabla \mu_c +p\diver(\phi_c\mathbf{v}_c+\phi_l\mathbf{v}_l)\biggr)=\int_{R(t)}\biggl( \mathbf{T}_c\colon \nabla \mathbf{v}_c+ \mathbf{T}_l\colon \nabla \mathbf{v}_l+\mathbf{m}\cdot (\mathbf{v}_c-\mathbf{v}_l)\\
&+(\pi+\mu_c) \left(\frac{\partial \phi_c}{\partial t}+\mathbf{v}_c\cdot \nabla \phi_c\right)-\pi\mathbf{v}_c\cdot \nabla \phi_c+\left(\mu_c\phi_c+p\phi_c\right)\mathbf{I}\colon \nabla \mathbf{v}_c+p\phi_l\mathbf{I}\colon \nabla \mathbf{v}_l-\mathbf{h}_c \cdot \nabla \mu_c\\
&+p\nabla \phi_c \cdot \mathbf{v}_c+p\nabla \phi_l \cdot \mathbf{v}_l+\boldsymbol{\xi} \cdot \nabla \frac{\partial \phi_c}{\partial t}\biggr)-\int_{\partial R(t)}\mu_c\phi_c\mathbf{v}_c\cdot \mathbf{n}.
\end{align*}
\begin{rem}
\label{rem:mixture}
We note that, in accordance with the axioms of the theory of mixtures \cite{allen} and due to the fact that the energy density $e$ is a quantity associated with the phase $c$, the advective velocity used in the Reynolds transport theorem is $\mathbf{v}_c$. This means that the quantity $e$ in the first term in \eqref{eqn:dissipation} is integrated over a material volume for the phase $c$, namely $R_c(t)$, which coincides with the material volume $R(t)$ for the whole mixture by assumption.
\end{rem}
Integrating by parts, and using \eqref{eqn:xi}, in the term 
\[
\int_{R(t)}\pi\mathbf{v}_c\cdot \nabla \phi_c=\int_{\partial R(t)}\left(\mathbf{v}_c\cdot \nabla \phi_c\right)\boldsymbol{\xi}\cdot \mathbf{n}-\int_{R(t)}\left(\boldsymbol{\xi}\cdot \nabla (\nabla \phi_c)\mathbf{v}_c+\nabla \phi_c \otimes \boldsymbol{\xi} \colon \nabla \mathbf{v}_c\right)
\]
we get
\begin{align}
\label{eqn:dissineq}
&\int_{R(t)}\left(\left(\partial_{\phi_c}e-\pi-\mu_c\right)\left(\frac{\partial \phi_c}{\partial t}+\mathbf{v}_c\cdot \nabla \phi_c\right)+\left(\partial_{\nabla \phi_c}e-\boldsymbol{\xi}\right)\cdot\left(\frac{\partial  \nabla \phi_c}{\partial t}+\nabla (\nabla \phi_c)\mathbf{v}_c\right)\right)\\
& \notag+\int_{R(t)}\biggl(\frac{1}{\det\mathbf{F}_c}(\partial_{\mathbf{F}_c}w)\mathbf{F}_c^T-\diver\left(\frac{1}{\det\mathbf{F}_c}\partial_{\nabla \mathbf{F}_c}w\right) \mathbf{F}_c^T-\frac{1}{\det\mathbf{F}_c}\nabla \mathbf{F}_c \odot \partial_{\nabla \mathbf{F}_c}w\\
& \notag -\nabla \phi_c \otimes \boldsymbol{\xi} +\biggl[\Phi(\phi_c, \nabla \phi_c)-\phi_c\mu_c-p\phi_c\biggr]\mathbf{I}-\mathbf{T}_c\biggr)\colon \nabla \mathbf{v}_c-\int_{R(t)}(\mathbf{T}_l+p\phi_l\mathbf{I})\colon \nabla \mathbf{v}_l\\
& \notag -\int_{R(t)}(\mathbf{m}+p\nabla \phi_c)\cdot \mathbf{v}_c+\int_{R(t)}(\mathbf{m}-p\nabla (1-\phi_c))\cdot \mathbf{v}_l+\int_{R(t)}\mathbf{h}_c \cdot \nabla \mu_c\\
& \notag +\int_{\partial R(t)}\left(\mathbf{J}_e+\frac{1}{\det\mathbf{F}_c}[\partial_{\nabla \mathbf{F}_c}w]^T\dot{\mathbf{F}}_c+(\mathbf{v}_c\cdot \nabla \phi_c)\boldsymbol{\xi}+\phi_c\mu_c\mathbf{v}_c\right)\cdot \mathbf{n}\leq 0.
\end{align}
Inequality \eqref{eqn:dissineq} is satisfied if we make the following general constitutive assumptions:
\begin{equation}
\label{eqn:constitutive}
\begin{cases}
&\mu_c=\partial_{\phi_c}e-\pi+\delta\left(\frac{\partial \phi_c}{\partial t}+\mathbf{v}_c\cdot \nabla \phi_c\right)=\partial_{\phi_c}e-\diver\boldsymbol{\xi}+\delta\left(\frac{\partial \phi_c}{\partial t}+\mathbf{v}_c\cdot \nabla \phi_c\right),\\
&\boldsymbol{\xi}=\partial_{\nabla \phi_c}e+\theta\left(\frac{\partial  \nabla \phi_c}{\partial t}+\nabla (\nabla \phi_c)\mathbf{v}_c\right),\\
&\mathbf{T}_c=\frac{1}{\det\mathbf{F}_c}(\partial_{\mathbf{F}_c}w)\mathbf{F}_c^T-\diver\left(\frac{1}{\det\mathbf{F}_c}\partial_{\nabla \mathbf{F}_c}w\right) \mathbf{F}_c^T-\frac{1}{\det\mathbf{F}_c}\nabla \mathbf{F}_c \odot \partial_{\nabla \mathbf{F}_c}w\\
&- \nabla \phi_c \otimes \boldsymbol{\xi}+\biggl[\Phi(\phi_c, \nabla \phi_c)-\phi_c\mu_c-p\phi_c\biggr]\mathbf{I}+\nu_c\left(\nabla \mathbf{v}_c+\left(\nabla \mathbf{v}_c\right)^T\right),\\
&\mathbf{T}_l=-p(1-\phi_c)\mathbf{I}+\nu_l\left(\nabla \mathbf{v}_l+\left(\nabla \mathbf{v}_l\right)^T\right),\\
&\mathbf{m}=-p\nabla \phi_c,\\
&\mathbf{h}_c =-b(\phi_c)\nabla \mu_c,\\
&\mathbf{J}_e=-\frac{1}{\det\mathbf{F}_c}[\partial_{\nabla \mathbf{F}_c}w]^T\dot{\mathbf{F}}_c-(\mathbf{v}_c\cdot \nabla \phi_c)\boldsymbol{\xi}-\phi_c\mu_c\mathbf{v}_c,
\end{cases}
\end{equation}
where \eqref{eqn:xi} has been used in the first relation of \eqref{eqn:constitutive}, $\delta, \theta, \nu_c, \nu_l$ are positive coefficients (which we assume to be constant) and $b(\phi_c)$ is a positive mobility function, considered as dependent only on $\phi_c$.

Inserting \eqref{eqn:constitutive} in \eqref{eqn:1a}, \eqref{eqn:tc} and \eqref{eqn:tl}, considering $\partial \Omega_{D,c}=\partial \Omega_{D,l}=\partial \Omega$, $\eta=0$, and imposing Neumann boundary conditions on $\partial_{\nabla \mathbf{F}_c}w$ and on the flux $\mathbf{h}_c$, we finally obtain the system 
\begin{equation}
\label{eqn:final}
\begin{cases}
\mathbf{T}_c=\frac{1}{\det\mathbf{F}_c}(\partial_{\mathbf{F}_c}w)\mathbf{F}_c^T-\diver\left(\frac{1}{\det\mathbf{F}_c}\partial_{\nabla \mathbf{F}_c}w\right) \mathbf{F}_c^T-\frac{1}{\det\mathbf{F}_c}\nabla \mathbf{F}_c \odot \partial_{\nabla \mathbf{F}_c}w- \nabla \phi_c \otimes \partial_{\nabla \phi_c}e\\
- \theta \nabla \phi_c \otimes \left(\frac{\partial  \nabla \phi_c}{\partial t}+\nabla (\nabla \phi_c)\mathbf{v}_c\right)+\biggl[\Phi(\phi_c, \nabla \phi_c)-\phi_c\mu_c-p\phi_c\biggr]\mathbf{I}+\nu_c\left(\nabla \mathbf{v}_c+\left(\nabla \mathbf{v}_c\right)^T\right),\\
\diver\mathbf{T}_c-K(\phi_c)(\mathbf{v}_c-\mathbf{v}_l)+p\nabla \phi_c=\boldsymbol{0},\\
\mathbf{T}_l=-p(1-\phi_c)\mathbf{I}+\nu_l\left(\nabla \mathbf{v}_l+\left(\nabla \mathbf{v}_l\right)^T\right),\\
\diver\mathbf{T}_l+K(\phi_c)(\mathbf{v}_c-\mathbf{v}_l)-p\nabla \phi_c=\boldsymbol{0},\\
\diver(\phi_c\mathbf{v}_c+\phi_l\mathbf{v}_l)=0,\\
\dot{\mathbf{F}}_c=(\nabla \mathbf{v}_c)\mathbf{F}_c,\\
\frac{\partial \phi_c}{\partial t}+\diver(\phi_c \mathbf{v}_c)-\diver(b(\phi_c)\nabla \mu_c)=0,\\
\mu_c=\partial_{\phi_c}e-\diver\partial_{\nabla \phi_c}e-\theta \diver\left(\frac{\partial  \nabla \phi_c}{\partial t}+\nabla (\nabla \phi_c)\mathbf{v}_c\right)+\delta \left(\frac{\partial \phi_c}{\partial t}+\mathbf{v}_c\cdot \nabla \phi_c\right),
\end{cases}
\end{equation}
valid in $\Omega \times [0,T]$, where $T>0$ is the final time of evolution of the system, endowed with the boundary conditions
\begin{equation}
\label{eqn:finalbc}
b(\phi_c)\nabla \mu_c\cdot \mathbf{n}=\left(\partial_{\nabla \phi_c}e+\theta\left(\frac{\partial  \nabla \phi_c}{\partial t}+\nabla (\nabla \phi_c)\mathbf{v}_c\right)\right) \cdot \mathbf{n}=0, \; \left[\partial_{\nabla \mathbf{F}_c}w\right]\mathbf{n}=\mathbf{0},\; \mathbf{v}_c=\mathbf{v}_l=\mathbf{0},
\end{equation}
on $\partial \Omega \times [0,T]$, and with initial conditions for $\mathbf{F}_c(\mathbf{x},0)$ and $\phi_c(\mathbf{x},0)$, $\mathbf{x}\in \Omega$. Here,
\begin{equation}
\label{eqn:finalbis}
\partial_{\phi_c}e=\frac{\psi'(\phi_c)}{\epsilon}+\frac{1}{\det\mathbf{F}_c}\partial_{\phi_c}w, \quad \partial_{\nabla \phi_c}e=\epsilon \nabla \phi_c +\frac{1}{\det\mathbf{F}_c} \partial_{\nabla \phi_c} w.
\end{equation}
With the choice of boundary conditions in \eqref{eqn:final}, we have that $\mathbf{J}_e\cdot \mathbf{n}=0$ on $\partial \Omega$. We note that the boundary condition $\mathbf{v}_c=\mathbf{v}_l=\mathbf{0}$ on $\partial \Omega\times [0,T]$ implies that the boundary of the domain $\Omega$ is not deformed over time by the elastic deformation, hence the domain $\Omega$ is fixed in time.
We moreover note that the elastic contribution $\frac{1}{\det\mathbf{F}_c}(\partial_{\mathbf{F}_c}w)\mathbf{F}_c^T$ in the Cauchy stress tensor $\mathbf{T}_c$ is the push--forward of the Piola stress tensor $\partial_{\mathbf{F}_c}w$ from the Lagrangian to the Eulerian configuration.
\begin{rem}
\label{rem:viscous}
The terms in \eqref{eqn:final} which are multiplied by the coefficients $\delta$ and $\theta$, which are dissipative contributions in \eqref{eqn:dissineq} and \eqref{eqn:constitutive}, can be identified as viscous regularizations of the material time derivatives of $\phi_c$ and $\nabla \phi_c$ in the system dynamics.
\end{rem}
\begin{rem}
\label{rem:secondorder}
The dependence of $w$ on $\nabla \mathbf{F}_c$ in \eqref{eqn:free} could be associated, in a framework of a second--gradient theory, to a third order stress tensor power--conjugated to $\nabla \nabla \mathbf{v}_c$, instead of being incorporated in the second order stress tensor $\mathbf{T}_c$ as done in \eqref{eqn:constitutive}.
\end{rem}
\subsection{Energy estimate}
In this section we prove that a solution of system \eqref{eqn:final} formally satisfies an energy equality, as a direct consequence of \eqref{eqn:dissineq} and of the constitutive assumptions \eqref{eqn:constitutive}.
Let us multiply \eqref{eqn:final}$_1$ and \eqref{eqn:final}$_2$ by $\nabla \mathbf{v}_c$ and $\mathbf{v}_c$ respectively, \eqref{eqn:final}$_3$ and \eqref{eqn:final}$_4$ by $\nabla \mathbf{v}_l$ and $\mathbf{v}_l$ respectively, \eqref{eqn:final}$_5$ by $p$, \eqref{eqn:final}$_6$ by $\partial_{\mathbf{F}_c}\left(\frac{w}{\det\mathbf{F}_c}\right)-\diver\left(\frac{1}{\det\mathbf{F}_c}\partial_{\nabla \mathbf{F}_c}w\right)$, \eqref{eqn:final}$_7$ by $\mu_c$ and \eqref{eqn:final}$_8$ by $\partial_t \phi_c+\mathbf{v}_c\cdot \nabla \phi_c$, and integrate over $\Omega$. We obtain:
\begin{equation}
\label{eqn:energy}
\begin{cases}
\biggl(\partial_{\mathbf{F}_c}\left(\frac{w}{\det\mathbf{F}_c}\right)\mathbf{F}_c^T- \left(\diver\left(\frac{1}{\det\mathbf{F}_c}\partial_{\nabla \mathbf{F}_c}w\right) \mathbf{F}_c^T+\frac{1}{\det\mathbf{F}_c}\nabla \mathbf{F}_c \odot \partial_{\nabla \mathbf{F}_c}w\right) 
-\nabla \phi_c\otimes \partial_{\nabla \phi_c}e\\  
- \theta \nabla \phi_c \otimes \left(\frac{\partial  \nabla \phi_c}{\partial t}+\nabla (\nabla \phi_c)\mathbf{v}_c\right)+\biggl[e(\phi_c, \nabla \phi_c,\mathbf{F}_c)-\phi_c\mu_c-p\phi_c\biggr]\mathbf{I}\\
+\nu_c\left(\nabla \mathbf{v}_c+\left(\nabla \mathbf{v}_c\right)^T\right),\nabla \mathbf{v}_c\biggr)+(K(\phi_c)(\mathbf{v}_c-\mathbf{v}_l),\mathbf{v}_l)-(p\nabla \phi_c,\mathbf{v}_c)=0,\\ \\
\left(-p(1-\phi_c)\mathbf{I}+\nu_l\left(\nabla \mathbf{v}_l+\left(\nabla \mathbf{v}_l\right)^T\right),\nabla \mathbf{v}_l\right)-(K(\phi_c)(\mathbf{v}_c-\mathbf{v}_l),\mathbf{v}_l)+(p\nabla \phi_c,\mathbf{v}_l)=0,\\ \\
\bigl(\phi_c\diver\mathbf{v}_c+\mathbf{v}_c\cdot \nabla \phi_c+(1-\phi_c)\diver\mathbf{v}_l+\mathbf{v}_l\cdot \nabla (1-\phi_c),p\bigr)=0,\\ \\
\left(\nabla \mathbf{v}_c,\partial_{\mathbf{F}_c}\left(\frac{w}{\det\mathbf{F}_c}\right)\mathbf{F}_c^{T}-\diver\left(\frac{1}{\det\mathbf{F}_c}\partial_{\nabla \mathbf{F}_c}w\right) \mathbf{F}_c^T\right)=\left(\dot{\mathbf{F}}_c,\partial_{\mathbf{F}_c}\left(\frac{w}{\det\mathbf{F}_c}\right)\right)\\
+\left(\nabla \dot{\mathbf{F}}_c,\frac{1}{\det\mathbf{F}_c}\partial_{\nabla \mathbf{F}_c}w\right)=\left(\dot{\mathbf{F}}_c,\partial_{\mathbf{F}_c}\left(\frac{w}{\det\mathbf{F}_c}\right)\right)
+\left(\dot{(\nabla \mathbf{F}_c)},\frac{1}{\det\mathbf{F}_c}\partial_{\nabla \mathbf{F}_c}w\right)\\
+\left(\frac{1}{\det\mathbf{F}_c}\nabla \mathbf{F}_c \odot \nabla \partial_{\nabla \mathbf{F}_c}w,\nabla \mathbf{v}_c\right),\\ \\
\left(\frac{\partial \phi_c}{\partial t},\mu_c\right)+(\phi_c \diver\mathbf{v}_c,\mu_c)+(\mathbf{v}_c\cdot \nabla \phi_c,\mu_c)+(b(\phi_c)\nabla \mu_c,\nabla \mu_c)=0,\\ \\
\left(\mu_c,\frac{\partial \phi_c}{\partial t}+\mathbf{v}_c\cdot \nabla \phi_c\right)= \left(\partial_{\phi_c}e,\frac{\partial \phi_c}{\partial t}+\mathbf{v}_c\cdot \nabla \phi_c\right)+ \left(\partial_{\nabla \phi_c}e,\nabla\left(\frac{\partial \phi_c}{\partial t}+\mathbf{v}_c\cdot \nabla \phi_c\right)\right)\\
+\theta\left(\frac{\partial  \nabla \phi_c}{\partial t}+\nabla (\nabla \phi_c)\mathbf{v}_c,\nabla\left(\frac{\partial \phi_c}{\partial t}+\mathbf{v}_c\cdot \nabla \phi_c\right)\right)+\delta\left|\left|\frac{\partial \phi_c}{\partial t}+\mathbf{v}_c\cdot \nabla \phi_c\right|\right|^2\\
=\left(\partial_{\phi_c}e,\frac{\partial \phi_c}{\partial t}+\mathbf{v}_c\cdot \nabla \phi_c\right)+ \left(\partial_{\nabla \phi_c}e,\frac{\partial  \nabla \phi_c}{\partial t}+\nabla (\nabla \phi_c)\mathbf{v}_c\right)+\theta\left|\left|\frac{\partial  \nabla \phi_c}{\partial t}+\nabla (\nabla \phi_c)\mathbf{v}_c\right|\right|^2 \\
+\delta\left|\left|\frac{\partial \phi_c}{\partial t}+\mathbf{v}_c\cdot \nabla \phi_c\right|\right|^2
+\left(\nabla \phi_c\otimes \left(\partial_{\nabla \phi_c}e+\theta\left(\frac{\partial  \nabla \phi_c}{\partial t}+\nabla (\nabla \phi_c)\mathbf{v}_c\right)\right),\nabla \mathbf{v}_c\right),
\end{cases}
\end{equation}
where in the {\rosso seventh line} we used \eqref{eqn:fradfdot}, and in {\rosso the last equality} we used the relation
\[
 \left(\partial_{\nabla \phi_c}e,\nabla\left(\mathbf{v}_c\cdot \nabla \phi_c\right)\right)=\left(\partial_{\nabla \phi_c}e,[\nabla \nabla \phi_c]\mathbf{v}_c\right)+\left(\nabla \phi_c\otimes \partial_{\nabla \phi_c}e,\nabla \mathbf{v}_c\right),
\]
and a similar relation for the term multiplied by the coefficient $\theta$. Summing all the equations in \eqref{eqn:energy} and using the boundary conditions \eqref{eqn:finalbc}, we obtain the Energy equality
\begin{multline}
\label{eqn:energy3}
\frac{d}{dt}\int_{\Omega}e(\phi_c,\nabla \phi_c,\mathbf{F}_c,\nabla \mathbf{F}_c)+\int_{\Omega}b(\phi_c)|\nabla \mu_c|^2+\frac{\nu_c}{2}\int_{\Omega}\left|\nabla \mathbf{v}_c+\left(\nabla \mathbf{v}_c\right)^T\right|^2+\frac{\nu_l}{2}\int_{\Omega}\left|\nabla \mathbf{v}_l+ \left(\nabla \mathbf{v}_l\right)^T\right|^2\\
+\frac{1}{2}\int_{\Omega}K(\phi_c)|\mathbf{v}_c-\mathbf{v}_l|^2+\theta\left|\left|\frac{\partial  \nabla \phi_c}{\partial t}+\nabla (\nabla \phi_c)\mathbf{v}_c\right|\right|^2 +\delta\left|\left|\frac{\partial \phi_c}{\partial t}+\mathbf{v}_c\cdot \nabla \phi_c\right|\right|^2=0.
\end{multline}
\subsection{Constitutive assumptions for $w$ and $K(\phi_c)$}
\label{sec:constitutive}
In the following, we specify the model derivation by considering an hyperelastic free energy density which depends only on $\phi_c$ and $\mathbf{F}_c$, i.e. $w=w(\phi_c,\mathbf{F}_c)$, and which satisfies objectivity, i.e. 
\[
w(\phi_c,\mathbf{F}_c)=\hat{w}(\phi_c,\mathbf{B}_c),
\]
where $\mathbf{B}_c=\mathbf{F}_c\mathbf{F}_c^T$. Of particular interest for both analytical aspects of the theory of nonlinear elasticity \cite{ball} and for applications \cite{ogden} are the policonvexity and the coercivity conditions for $w$. The function $w$ is said to be polyconvex if (in $3$ space dimensions) there exists a function $g:\Omega \times \mathbb{R}^{19}\to \mathbb{R}$ with $g(\mathbf{x},\cdot,\cdot,\cdot)$ convex for a.e.  $\mathbf{x}\in \Omega$ and 
\begin{equation}
\label{eqn:poly}
w(\phi_c,\mathbf{F}_c)=g(\phi_c,\mathbf{F}_c,\text{cof}\mathbf{F}_c,\det\mathbf{F}_c).
\end{equation}
Moreover, if $w$ is polyconvex, the following coercivity condition on $g(\phi_c,\cdot,\cdot,\cdot)$ is typically required: there exist $\alpha >0$, $\beta \in \mathbb{R}$, $p\geq 2$, $q\geq \frac{p}{p-1}$, $r>0$ such that, for each $\mathbf{F},\mathbf{C}\in \mathcal{L}(\mathbb{R}^3)$, $\eta >0$,
\begin{equation}
\label{eqn:coer}
g(\phi_c,\mathbf{F},\mathbf{C},\eta)\geq \alpha \left(|\mathbf{F}|^p+|\mathbf{C}|^q+\eta^r\right)+\beta.
\end{equation}
Many realistic models of nonlinear elastic materials satisfy the polyconvexity and the coercivity properties, e.g. the generalized Ogden model and the Mooney-Rivlin model \cite{ogden}. Specifically, the isotropic Mooney--Rivlin model, which is a particular case of the Ogden model,  has been used in literature to describe the mechanical behavior of soft brain tissues, and it has been experimentally calibrated and validated on in--vivo data \cite{balbi}.

In light of these considerations, in the following we specify the model derivation by considering the following Mooney--Rivlin form (in $3$ space dimensions) for $w$,
\begin{equation}
\label{eqn:freemr}
w(\phi_c,\mathbf{F}_c)= \frac{c_1}{2}f(\phi_c)(\mathbf{F}_c\colon \mathbf{F}_c-3)+\frac{c_2}{2}g(\phi_c)(\text{cof}(\mathbf{F}_c)\colon \text{cof}(\mathbf{F}_c)-3)+c_3h(\det\mathbf{F}_c),
\end{equation}
with
\begin{align}
\label{eqn:dfreemr}
\partial_{\mathbf{F}_c}w(\phi_c,\mathbf{F}_c)=&c_1f(\phi_c)\mathbf{F}_c+c_2g(\phi_c)\left(\left(\text{cof}(\mathbf{F}_c)\colon \text{cof}(\mathbf{F}_c)\right)\mathbf{F}_c^{-T}-\text{cof}(\mathbf{F}_c)\left(\text{cof}(\mathbf{F}_c)\right)^T\mathbf{F}_c^{-T}\right)\\
& \notag +c_3h'(\det\mathbf{F}_c)\det\mathbf{F}_c\mathbf{F}_c^{-T},
\end{align}
where $c_1,c_2,c_3$ are material parameters, in units of $[Pa]$, the functions $f(\phi_c), g(\phi_c)$ represent the change in stiffness of the material depending on the concentration $\phi_c$, and $h:\mathbb{R}^+\to \mathbb{R}$ is a convex function which penalizes the infinite compression situation $\det\mathbf{F_c}=0$ in the energy density and which eventually contains a weak compressibility contribution. 
\begin{rem}
\label{rem:fg}
Realistic constitutive assumptions for the functions $f(\phi_c), g(\phi_c)$, if the mixture is composed by an elastic and an anelastic components, should satisfy the property $f(0)=g(0)=0$, with the particular form of $f,g:[0,1]\to [0,1]$ determined by empirical observations on the considered materials. 
\end{rem}
We make the simplifying assumption that the momentum transfer in the mixture due to shear stresses in the liquid is negligible with respect to the momentum transfer between the solid and the liquid components, i.e. we take
\begin{equation}
\label{eqn:tldarcy}
\mathbf{T}_l=-p(1-\phi_c)\mathbf{I}.
\end{equation}
We moreover assume that the momentum exchange between the phases $c$ and $l$, expressed as a Stokes drag in \eqref{eqn:tc} and \eqref{eqn:tl}, is induced by a Darcy--like flow of the viscous fluid associated to the phase $l$ through the porous-permeable solid matrix associated to the phase $c$, hence (see e.g. \cite{allen})
\begin{equation}
\label{eqn:kphi}
K(\phi_c)=\frac{\mu_l (1-\phi_c)}{k(\phi_c)},
\end{equation}
where $\mu_l$ is the viscosity of the phase $l$ and $k(\phi_c)$ is the intrinsic permeability of the phase $c$. Assuming $\mu_l=\mu \phi_l$ and $k(\phi_c)=k$ constant, we obtain
\begin{equation}
\label{eqn:kphi2}
K(\phi_c)=\frac{\mu (1-\phi_c)^2}{k}.
\end{equation}
Inserting \eqref{eqn:tldarcy} and \eqref{eqn:kphi2} in the third and fourth equations of system \eqref{eqn:final}, we obtain the Darcy--like relation 
\begin{equation}
\label{eqn:darcy}
\diver\mathbf{T}_l-\frac{\mu}{k}(1-\phi_c)^2(\mathbf{v}_l-\mathbf{v}_c)-p\nabla \phi_c=0 \to (\mathbf{v}_l-\mathbf{v}_c)=-\frac{k}{\mu(1-\phi_c)}\nabla p.
\end{equation}
We finally substitute \eqref{eqn:darcy} in the fifth equation of \eqref{eqn:final}, thus eliminating the variables $\mathbf{T}_l$ and $\mathbf{v}_l$ from the system, and consider the Mooney--Rivlin model \eqref{eqn:freemr}, with \eqref{eqn:dfreemr}, taking also $\delta=\theta=0$. 
Hence, system \eqref{eqn:final} reduces to:
\begin{equation}
\label{eqn:final2}
\begin{cases}
\mathbf{T}_c=\frac{1}{\det\mathbf{F}_c}\left(c_1f(\phi_c)\mathbf{F}_c\mathbf{F}_c^T-c_2g(\phi_c)\text{cof}(\mathbf{F}_c)(\text{cof}(\mathbf{F}_c))^T\right)
- \epsilon \nabla \phi_c \otimes \nabla \phi_c\\
+\left[\frac{c_2}{\det\mathbf{F}_c}g(\phi_c)(\text{cof}(\mathbf{F_c})\colon \text{cof}(\mathbf{F_c}))+\Phi(\phi_c, \nabla \phi_c)-\phi_c\mu_c+c_3h'(\det\mathbf{F}_c)\right]\mathbf{I}+\nu_c\left(\nabla \mathbf{v}_c+\left(\nabla \mathbf{v}_c\right)^T\right),\\
\diver\mathbf{T}_c-\nabla p=0,\\
\diver\left(\frac{k}{\mu}\nabla p\right)=\diver\mathbf{v}_c,\\
\dot{\mathbf{F}}_c=(\nabla \mathbf{v}_c)\mathbf{F}_c,\\
\frac{\partial \phi_c}{\partial t}+\diver(\phi_c \mathbf{v}_c)-\diver(b(\phi_c)\nabla \mu_c)=0,\\
\mu_c=\frac{\psi'(\phi_c)}{\epsilon}-\epsilon \Delta \phi_c+\frac{c_1}{2\det\mathbf{F}_c}f'(\phi_c)(\mathbf{F}_c\colon \mathbf{F}_c-3)+\frac{c_2}{2\det\mathbf{F}_c}g'(\phi_c)(\text{cof}(\mathbf{F_c})\colon \text{cof}(\mathbf{F_c})-3),
\end{cases}
\end{equation}
valid in $\Omega \times [0,T]$, endowed with the boundary conditions
\begin{equation}
\label{eqn:final2bc}
b(\phi_c)\nabla \mu_c\cdot \mathbf{n}=\nabla \phi_c \cdot \mathbf{n}=\nabla p \cdot \mathbf{n}=0, \; \mathbf{v}_c=\mathbf{0},
\end{equation}
on $\partial \Omega \times [0,T]$, and with initial conditions for $\mathbf{F}_c(\mathbf{x},0)$ and $\phi_c(\mathbf{x},0)$, $\mathbf{x}\in \Omega$.
\begin{rem}
\label{rem:darcyvariants}
Different possible constitutive assumptions on the form of the intrinsic permeability $k(\phi_c)$ in \eqref{eqn:kphi} could be given, depending on the particular microstructure of the considered constituent $c$. For instance, $k(\phi_c)\propto (1-\phi_c)$ if the phase $c$ consists of homogeneous and isotropic parallel cylindrical pores, and using the Poiseuille formula for a capillary tube \cite{holmes}, while $k(\phi_c)\propto \frac{(1-\phi_c)^2}{\phi_c}$ when a Kozeny--Carman permeability-porosity relation is assumed for the porous medium \cite{kozeny,carman}. In these situations, \eqref{eqn:final2}$_3$ would take the form of a degenerate elliptic equation for $p$, with a mobility depending on $\phi_c$. 
\end{rem}
A solution of \eqref{eqn:final2} formally satisfies an energy equality analogous to \eqref{eqn:energy3}. Indeed, multiplying \eqref{eqn:final2}$_1$ and \eqref{eqn:final2}$_2$ by $\nabla \mathbf{v}_c$ and $\mathbf{v}_c$ respectively, \eqref{eqn:final2}$_3$ by $p$, \eqref{eqn:final2}$_4$ by $\partial_{\mathbf{F}_c}\left(\frac{w}{\det\mathbf{F}_c}\right)$, \eqref{eqn:final2}$_5$ by $\mu_c$ and \eqref{eqn:final2}$_6$ by $\partial_t \phi_c+\mathbf{v}_c\cdot \nabla \phi_c$, integrating over $\Omega$ and employing similar calculations as in \eqref{eqn:energy} and the boundary conditions, we obtain
\begin{equation}
\label{eqn:energy3darcy}
\frac{d}{dt}\int_{\Omega}e(\phi_c,\nabla \phi_c,\mathbf{F}_c)+\int_{\Omega}b(\phi_c)|\nabla \mu_c|^2+\frac{\nu_c}{2}\int_{\Omega}\left|\nabla \mathbf{v}_c+\left(\nabla \mathbf{v}_c\right)^T\right|^2+\int_{\Omega}\frac{k}{\mu}|\nabla p|^2=0.
\end{equation}
\subsection{Elastic incompressibility constraint for the cell phase}
Let us introduce a solenoidal  constraint for the velocity of the phase $c$, i.e. $\diver\mathbf{v}_c=0$, with a corresponding Lagrange multiplier in the dissipation inequality \eqref{eqn:dissipation},
\begin{align}
\label{eqn:dissipation2}
&\frac{d}{dt}\int_{R(t)}e(\phi_c,\nabla \phi_c,\mathbf{F}_c,\nabla \mathbf{F}_c)+\int_{\partial R(t)}\mathbf{J}_e\cdot \mathbf{n}\leq W_{\text{ext}}(R(t))-\int_{\partial R(t)}\mu_c(\phi_c\mathbf{v}_c+\mathbf{h}_c)\cdot \mathbf{n}\\
& \notag +\int_{R(t)}p\diver(\phi_c\mathbf{v}_c+\phi_l\mathbf{v}_l)+\int_{R(t)}q\diver\mathbf{v}_c.
\end{align}
This constraint implies that the partial volume occupied by the phase $c$ within a material volume changes only due to the presence of the micro--flux. Analogously to the mixture model introduced in \cite{frieboes}, this modeling choice introduces two kind of pressures in the system dynamics, a liquid pressure $p$ and a solid pressure $q$.
\begin{rem}
\label{rem:incompressibility}
The constraint $\diver\mathbf{v}_c=0$ implies $\det\mathbf{F}_c=\text{constant}$. Indeed, taking the Frobenius scalar product of equation \eqref{eqn:fdot} with $\partial_{\mathbf{F}_c}\det\mathbf{F}_c$, we have that
\begin{equation}
\label{eqn:incsol}
\frac{d (\det\mathbf{F}_c)}{dt}=\dot{\mathbf{F}}_c\colon \det\mathbf{F}_c\mathbf{F}_c^{-T}=(\nabla \mathbf{v}_c)\mathbf{F}_c\colon \det\mathbf{F}_c\mathbf{F}_c^{-T}=(\nabla \mathbf{v}_c)\colon \mathbf{I} \det\mathbf{F}_c=\diver\mathbf{v}_c\det\mathbf{F}_c=0.
\end{equation}
\end{rem}
Starting from \eqref{eqn:dissipation2}, using similar calculations as in \eqref{eqn:dissineq} and \eqref{eqn:constitutive}, imposing $\det\mathbf{F}_c(\mathbf{x},t)=\det\mathbf{F}_c(\mathbf{x},0)=1$, and redefining $q=q-\Phi(\phi_c, \nabla \phi_c)+\phi_c\mu_c$, we obtain the system
\begin{equation}
\label{eqn:finalinc}
\begin{cases}
\mathbf{T}_c=(\partial_{\mathbf{F}_c}w)\mathbf{F}_c^T-\diver\left(\partial_{\nabla \mathbf{F}_c}w\right) \mathbf{F}_c^T-\nabla \mathbf{F}_c \odot \partial_{\nabla \mathbf{F}_c}w- \nabla \phi_c \otimes \partial_{\nabla \phi_c}e\\
- \theta \nabla \phi_c \otimes \left(\frac{\partial  \nabla \phi_c}{\partial t}+\nabla (\nabla \phi_c)\mathbf{v}_c\right)+\left[-p\phi_c-q\right]\mathbf{I}+\nu_c\left(\nabla \mathbf{v}_c+\left(\nabla \mathbf{v}_c\right)^T\right),\\
\diver\mathbf{T}_c-K(\phi_c)(\mathbf{v}_c-\mathbf{v}_l)+p\nabla \phi_c=0,\\
\mathbf{T}_l=-p(1-\phi_c)\mathbf{I}+\nu_l\left(\nabla \mathbf{v}_l+\left(\nabla \mathbf{v}_l\right)^T\right),\\
\diver\mathbf{T}_l+K(\phi_c)(\mathbf{v}_c-\mathbf{v}_l)-p\nabla \phi_c=0,\\
\diver(\phi_c\mathbf{v}_c+\phi_l\mathbf{v}_l)=0,\\
\diver\mathbf{v}_c=0,\\
\dot{\mathbf{F}}_c=(\nabla \mathbf{v}_c)\mathbf{F}_c,\\
\frac{\partial \phi_c}{\partial t}+\diver(\phi_c \mathbf{v}_c)-\diver(b(\phi_c)\nabla \mu_c)=0,\\
\mu_c=\partial_{\phi_c}e-\diver\partial_{\nabla \phi_c}e-\theta \diver\left(\frac{\partial  \nabla \phi_c}{\partial t}+\nabla (\nabla \phi_c)\mathbf{v}_c\right)+\delta \left(\frac{\partial \phi_c}{\partial t}+\mathbf{v}_c\cdot \nabla \phi_c\right),\\
\end{cases}
\end{equation}
valid in $\Omega \times [0,T]$, endowed with the boundary conditions
\begin{equation}
\label{eqn:finalincbc}
b(\phi_c)\nabla \mu_c\cdot \mathbf{n}=\left(\partial_{\nabla \phi_c}e+\theta\left(\frac{\partial  \nabla \phi_c}{\partial t}+\nabla (\nabla \phi_c)\mathbf{v}_c\right)\right) \cdot \mathbf{n}=0, \; \left[\partial_{\nabla \mathbf{F}_c}w\right]\mathbf{n}=\mathbf{0},\; \mathbf{v}_c=\mathbf{v}_l=\mathbf{0}
\end{equation}
on $\partial \Omega \times [0,T]$, and with initial conditions for $\mathbf{F}_c(\mathbf{x},0)$ and $\phi_c(\mathbf{x},0)$, $\mathbf{x}\in \Omega$.
Note that the energy estimate \eqref{eqn:energy3} is still formally satisfied by a solution of \eqref{eqn:finalinc}.
\begin{rem}
\label{rem:independence}
We observe that $\diver(\phi_c\mathbf{v}_c+\phi_l\mathbf{v}_l)=0$ and $\diver\mathbf{v}_c=0$ are independent constraints, as it is easily verifiable by direct calculations.
\end{rem}
We finally rewrite the system \eqref{eqn:finalinc} by introducing the constitutive assumption \eqref{eqn:kphi2} and considering the Mooney--Rivlin model \eqref{eqn:freemr} without the dependence on $h(\det\mathbf{F}_c)$, due to \eqref{eqn:incsol}, taking also $\delta=\theta=0$. We observe that \eqref{eqn:final2}$_3$, which now takes the form
\[
\diver\left(\frac{k}{\mu}\nabla p\right)=\diver\mathbf{v}_c=0,
\] 
gives, together with the boundary condition $\nabla p \cdot \mathbf{n}=0$ on $\partial \Omega$, that $p=\text{constant}$. This would be true also if we would consider an homogeneous degenerate elliptic equation for $p$ starting from different constitutive assumptions for the friction function $K(\phi_c)$ than \eqref{eqn:kphi2}, as described in Remark \ref{rem:darcyvariants}. Hence, the contribution of the variable $p$ disappears from \eqref{eqn:final2}$_2$. Redefining the variable $q=q-c_2g(\phi_c)(\text{cof}(\mathbf{F_c})\colon \text{cof}(\mathbf{F_c}))-\Phi(\phi_c, \nabla \phi_c)+\phi_c\mu_c$, and considering that 
\[\diver\left(\nabla \mathbf{v}_c+\left(\nabla \mathbf{v}_c\right)^T\right)=\Delta \mathbf{v}_c + \nabla \left(\diver\mathbf{v}_c\right),\]
we obtain
\begin{equation}
\label{eqn:finalinc2}
\begin{cases}
-\nu_c\Delta \mathbf{v}_c+\nabla q=-\epsilon \diver\left(\nabla \phi_c \otimes \nabla \phi_c\right)+\diver\left(c_1f(\phi_c)\mathbf{F}_c\mathbf{F}_c^T-c_2g(\phi_c)\text{cof}(\mathbf{F}_c)(\text{cof}(\mathbf{F}_c))^T\right),\\
\diver\mathbf{v}_c=0,\\
\dot{\mathbf{F}}_c=(\nabla \mathbf{v}_c)\mathbf{F}_c,\\
\frac{\partial \phi_c}{\partial t}+\mathbf{v}_c\cdot \nabla \phi_c-\diver(b(\phi_c)\nabla \mu_c)=0,\\
\mu_c=\frac{\psi'(\phi_c)}{\epsilon}-\epsilon \Delta \phi_c+\frac{c_1}{2}f'(\phi_c)(\mathbf{F}_c\colon \mathbf{F}_c-3)+\frac{c_2}{2}g'(\phi_c)(\text{cof}(\mathbf{F_c})\colon \text{cof}(\mathbf{F_c})-3),
\end{cases}
\end{equation}
valid in $\Omega \times [0,T]$, endowed with the boundary conditions
\begin{equation}
\label{eqn:finalinc2bc}
b(\phi_c)\nabla \mu_c\cdot \mathbf{n}=\nabla \phi_c \cdot \mathbf{n}=0, \; \mathbf{v}_c=\mathbf{0}
\end{equation}
on $\partial \Omega \times [0,T]$, and with initial conditions for $\mathbf{F}_c(\mathbf{x},0)$ and $\phi_c(\mathbf{x},0)$, $\mathbf{x}\in \Omega$.
To system \eqref{eqn:finalinc2}, with boundary conditions \eqref{eqn:finalinc2bc}, is associated an analogous energy estimate to \eqref{eqn:energy3darcy}, without the dissipative term depending on $\nabla p$, 
\begin{equation}
\label{eqn:energy3darcyinc}
\frac{d}{dt}\int_{\Omega}e(\phi_c,\nabla \phi_c,\mathbf{F}_c)+\int_{\Omega}b(\phi_c)|\nabla \mu_c|^2+\nu_c\int_{\Omega}\left|\nabla \mathbf{v}_c\right|^2=0.
\end{equation}

\section{Global existence of a weak solution}
\label{sec:2}
In this section we develop the analysis of existence of a global weak solution to a simplified version of system \eqref{eqn:finalinc2}, both in $2$ and $3$ space dimensions, obtained by regularizing the system and considering a Neo--Hookean model for the finite elasticity,
\begin{equation}
\label{eqn:neohook}
w(\phi_c,\mathbf{F}_c)= \frac{c_1}{2}f(\phi_c)(\mathbf{F}_c\colon \mathbf{F}_c-d),
\end{equation} 
where $d=2,3$, which is both polyconvex ($w(\cdot,\mathbf{F}_c)$ is strictly convex) and coercive, with $p=2$, in the sense of \eqref{eqn:poly} and \eqref{eqn:coer}. 
In particular, we add to system \eqref{eqn:finalinc2} {\rosso a regularizing term $\delta \frac{\partial \phi_c}{\partial t}$ in the expression for the chemical potential $\mu_c$, in order to enhance the regularity in time of the variable $\phi_c$, while we add a regularizing term $\lambda\Delta \left(f(\phi)\mathbf{F}\right)$ in the transport equation for $\mathbf{F}_c$, in order to enhance the regularity in space of the variable $\mathbf{F}_c$, hence improving the compactness properties of sequences of approximating solutions to the system. We will obtain existence of a weak solution in two space dimensions in the case $\lambda>0,\,\delta=0$ (i.e. no time derivative regularization is needed in $2$-d), while we will show the existence of a weak solution in three space dimensions in the case $\lambda>0,\,\delta>0$ (i.e. the viscous regularization term in the expression of the chemical potential $\mu_c$ is needed to recover existence in $3$-d)}.

{\rosso As \[\mu_c=-\epsilon \Delta \phi_c+\frac{1}{\epsilon}\psi'(\phi_c)+\partial_{\phi_c}w,\] we can obtain the relation}
\begin{align}
\label{eqn:gradgradparts}
&-\epsilon \diver\left(\nabla \phi_c \otimes \nabla \phi_c\right)=-\nabla\left(\frac{\epsilon}{2}|\nabla \phi_c|^2+\frac{1}{\epsilon}\psi(\phi_c)\right)+\mu_c\nabla \phi_c-\partial_{\phi_c}w\nabla \phi_c\\
& \notag =-\nabla\left(\frac{\epsilon}{2}|\nabla \phi_c|^2+\frac{1}{\epsilon}\psi(\phi_c)+w(\phi_c,\mathbf{F}_c)\right)+\mu_c\nabla \phi_c+\left(\nabla \mathbf{F}_c\right)^T\partial_{\mathbf{F}_c}w,
\end{align}
to rewrite the first term on the right hand side of \eqref{eqn:finalinc2}$_1$. Omitting the subscript $c$, setting $\epsilon=c_1=1$ for ease of notation {\rosso and employing \eqref{eqn:neohook}, we thus arrive at} the following system
{\rosso
\begin{equation}
\label{eqn:systemnoreg}
\begin{cases}
-\nu\Delta \mathbf{v}+\nabla q=\mu \nabla \phi- \left(\frac{1}{2}f'(\phi)(\mathbf{F}\colon \mathbf{F}-d)\right)\nabla \phi+\diver\left(f(\phi)\mathbf{F}\mathbf{F}^T\right),\\
\diver\mathbf{v}=0,\\
\frac{\partial \mathbf{F}}{\partial t}+\left(\mathbf{v}\cdot \nabla\right)\mathbf{F}-(\nabla \mathbf{v})\mathbf{F}-\lambda\Delta \left(f(\phi)\mathbf{F}\right)=0,\\
\frac{\partial \phi}{\partial t}+\mathbf{v}\cdot \nabla \phi-\diver(b(\phi)\nabla \mu)=0,\\
\mu=\psi'(\phi)- \Delta \phi+\frac{1}{2}f'(\phi)(\mathbf{F}\colon \mathbf{F}-d)+\delta \frac{\partial \phi}{\partial t},\\
\end{cases}
\end{equation}
valid in $\Omega_T$, endowed with the boundary conditions
\begin{equation}
\label{eqn:systembc}
b(\phi)\nabla \mu\cdot \mathbf{n}=\nabla \phi \cdot \mathbf{n}=0, \quad \mathbf{v}=\mathbf{0}, \quad \left[\nabla (f(\phi)\mathbf{F})\right]\mathbf{n}=\boldsymbol{0}
\end{equation}
on $\partial \Omega_T$, and with initial conditions for $\mathbf{F}(\mathbf{x},0)=\mathbf{F}_0$ and $\phi(\mathbf{x},0)=\phi_0$, $\mathbf{x}\in \Omega$. Here, {\rosso in view of \eqref{eqn:Phidef}, we have redefined} $q=q+\Phi(\phi_c,\nabla \phi_c)$.}
\begin{rem}
\label{rem:reg2}
As we will see in the next paragraphs, the viscous regularization inserted in the chemical potential in \eqref{eqn:systemnoreg}$_5$ aims at bounding the $L^2(\Omega_T)$ norm of  $\frac{\partial \phi}{\partial t}$ using the a--priori energy estimates associated to the system. A different viscous regularization of the time derivative of $\phi$ was introduced in \eqref{eqn:constitutive}, which was compatible with the dissipation inequality \eqref{eqn:dissineq} and which aimed at bounding the $L^2(\Omega_T)$ {\rosso norm} of the material time derivative of $\phi$ through \eqref{eqn:energy3}, obtained when the first term on the right hand side of \eqref{eqn:finalinc2}$_1$ is written as $-\epsilon \diver\left(\nabla \phi_c \otimes \nabla \phi_c\right)$. Since we will need to enhance the regularity of the partial time derivative of $\phi$ in order to obtain proper compactness results for the convergence of approximate solutions in $3$-d, we have rewritten the first term in the right hand side of \eqref{eqn:finalinc2}$_1$ using \eqref{eqn:gradgradparts} and inserted by hand a viscous regularization involving only the partial time derivative of $\phi$ in the chemical potential.
\end{rem}
\begin{rem}
\label{rem:reg3}
The regularization term $\lambda\Delta \left(f(\phi)\mathbf{F}\right)$ in \eqref{eqn:systemnoreg}$_3$ is introduced, following \cite{benesova}, to enhance the regularity in space of the variable $\mathbf{F}$, without {\rosso considering} a dependence on $\nabla \mathbf{F}$ in the elastic energy density $w$. Differently from \cite{benesova,garcke2} and due to the combined dependence of $w$ on $\phi$ and $\mathbf{F}$ in \eqref{eqn:neohook}, this elliptic regularization term must depend also on $\phi$ in order to be able to derive an a--priori energy estimates in a closed form, and its particular form prevents to obtain elliptic regularity for the variable $\mathbf{F}$. We note that, due to the presence of this regularization term, the constraint $\diver\mathbf{v}_c=0$ no longer implies, from \eqref{eqn:incsol}, that $\det\mathbf{F}_c=\text{constant}$. The kinematic relationship between the variables $\mathbf{v}_c$ and $\mathbf{F}_c$ is therefore lost due to the presence of the regularization term. 
\end{rem}
In order to proceed, we make the following assumptions:
\begin{itemize}
\item[\bf{A0}] {\rosso $\Omega$ is smooth and the boundary $\partial \Omega$ of the domain is of class $C^{\infty}$};
\item[\bf{A1}] $b\in C^{0}(\mathbb{R})$ and there exist $b_0,b_1>0$ such that $b_0\leq b(s) \leq b_1,$ $\forall s \in \mathbb{R}$;
\item[\bf{A2}] $\psi \in C^{1}(\mathbb{R})$ and there exist $c_1,c_2>0$ such that $|\psi'(s)|\leq c_1\left(|s|^p+1\right)$, $\psi(s)\geq -c_2$, $\forall s \in \mathbb{R}$, with $p\in [0,6)$ if $d=3$ and $p\in \mathbb{R}^+$ arbitrary if $d=2$. Moreover, there exists a convex decomposition of $\psi=\psi_++\psi_-$, where $\psi_+$ is convex and $\psi_-$ is concave, such that $|\psi_-''(s)|\leq c_1\left(|s|^q+1\right)$, with $q\in [0,4)$ if $d=3$ and $q\in \mathbb{R}^+$ arbitrary if $d=2$;
\item[\bf{A3}] $f\in C^{1,1}(\mathbb{R})$, and there exists a constant $f_{\text{min}} \in \mathbb{R}^+$ such that $0<f_{\text{min}}\leq f(s)\leq 1,$ $\forall s \in \mathbb{R}$. We also assume that there exists a constant $f_{1} \in \mathbb{R}^+$ such that $ |f'(s)| \leq f_1,$ $\forall s \in \mathbb{R}$;
\item[\bf{A4}] {\rosso The initial data have the regularity $\mathbf{F}_0\in L^2(\Omega;\mathbb{R}^{d\times d})$, $\phi_0\in H^1(\Omega)$.}
\end{itemize}
Note that through \textbf{A3} we assume the function $f$ to be uniformly positive, which, as will be evident in the following, is a fundamental property to obtain analytic results for system \eqref{eqn:systemnoreg}, even if realistic constitutive assumptions for $f$ should entail the degeneracy $f(0)=0$ (see Remark \ref{rem:fg}). 

We state here the main theorems of the present work concerning the existence of a global weak solution to \eqref{eqn:systemnoreg} in $2$ and $3$ space dimensions, which will be proved in the forthcoming sections.
\begin{thm}[$2$-d]
\label{thm:2d}
Let $d=2$ and let the assumptions \textbf{A0}--\textbf{A4} be satisfied. Then, there exists a weak solution $(\mathbf{v}, \mathbf{F}, \phi, \mu)$ of \eqref{eqn:systemnoreg}-\eqref{eqn:systembc} {\rosso for $\delta=0$, $\lambda>0$}, with
\begin{align*}
&\mathbf{v}\in L^{2}\left(0,T;H_{0,\diver}^1\left(\Omega;\mathbb{R}^{2}\right)\right),\\
&\mathbf{F} \in L^{\infty}(0,T;L^2(\Omega;\mathbb{R}^{2\times 2}))\cap L^{2}(0,T;H^1(\Omega;\mathbb{R}^{2\times 2})),\\
&\partial_t \mathbf{F} \in L^{2-h}\left(0,T;\left(H^1(\Omega;\mathbb{R}^{2\times 2})\right)'\right), \; \; h\in (0,1],\\
&\phi \in L^{\infty}(0,T;H^1(\Omega))\cap L^{2}(0,T;H^2(\Omega)),\\
&\partial_t \phi \in L^{2}\left(0,T;\left(H^1(\Omega)\right)'\right),\\
&\mu \in L^{2}\left(0,T;H^1\left(\Omega\right)\right),
\end{align*}
such that
\begin{equation}
\label{eqn:continuous2d}
\begin{cases}
\displaystyle \nu \int_{\Omega}\nabla \mathbf{v}\colon \nabla \mathbf{u}=\int_{\Omega}\mu \nabla \phi \cdot \mathbf{u}-\int_{\Omega}\frac{f'\left(\phi \right)}{2}\left(\mathbf{F}\colon \mathbf{F}-2\right)\nabla \phi \cdot \mathbf{u}-\int_{\Omega}f\left(\phi \right)\mathbf{F}\mathbf{F}^T\colon \nabla \mathbf{u},\\
\displaystyle  \left<\partial_t \mathbf{F}, \boldsymbol{\Sigma} \right>-\int_{\Omega}\mathbf{F}\otimes \mathbf{v}\mathbin{\tensorm} \nabla \boldsymbol{\Sigma}-\int_{\Omega}\left(\nabla {\mathbf{v}}\right)\mathbf{F}\colon \boldsymbol{\Sigma}+\lambda \int_{\Omega}\nabla \left(f({\phi})\mathbf{F}\right)\mathbin{\tensorm} \nabla \boldsymbol{\Sigma}=0,\\
 \displaystyle <\partial_t\phi, q>+\int_{\Omega}\left(\mathbf{v} \cdot \nabla \phi\right) q+\int_{\Omega}b\left(\phi \right)\nabla \mu \cdot \nabla q=0,\\
\displaystyle \int_{\Omega}\mu r=\int_{\Omega}\nabla \phi\cdot \nabla r+\int_{\Omega}\psi'\left(\phi\right)r+\int_{\Omega}\frac{f'\left(\phi\right)}{2}\left(\mathbf{F}\colon \mathbf{F}-2\right)r,
\end{cases}
\end{equation} 
a.e. in $t\in (0,T)$ and for all $\mathbf{u}\in H_{0,\diver}^1\left(\Omega;\mathbb{R}^{2}\right)$, $\boldsymbol{\Sigma}\in H^1\left(\Omega;\mathbb{R}^{2\times 2}\right)$, $q,r\in H^1(\Omega)$, satisfying the initial conditions $\mathbf{F}(\cdot,0)=\mathbf{F}_0$ a.e. in $\Omega$ and $\phi(\cdot,0)=\phi_0$ a.e. in $\Omega$.
\end{thm}

\begin{thm}[$3$-d]
\label{thm:3d}
Let $d=3$ and let the assumptions \textbf{A0}--\textbf{A4} be satisfied. Then, there exists a weak solution $(\mathbf{v}, \mathbf{F}, \phi, \mu)$ of \eqref{eqn:systemnoreg}-\eqref{eqn:systembc}, with {\rosso $\lambda,\, \delta>0$}, with
\begin{align*}
&\mathbf{v}\in L^{2}\left(0,T;H_{0,\diver}^1\left(\Omega;\mathbb{R}^{3}\right)\right),\\
&f(\phi)\mathbf{F} \in L^{\infty}(0,T;L^2(\Omega;\mathbb{R}^{3\times 3}))\cap L^{2}(0,T;H^1(\Omega;\mathbb{R}^{3\times 3})),\\
&\mathbf{F} \in L^{\infty}(0,T;L^2(\Omega;\mathbb{R}^{3\times 3})),\\
&\partial_t \mathbf{F} \in L^{\frac{4}{3}}\left(0,T;\left(H^1(\Omega;\mathbb{R}^{3\times 3})\right)'\right),\\
&\phi \in L^{\infty}(0,T;H^1(\Omega))\cap L^{\frac{4}{3}}(0,T;H^2(\Omega)),\\
&\partial_t \phi \in L^{2}\left(\Omega_T\right),\\
&\mu \in L^{2}\left(0,T;H^1\left(\Omega\right)\right),
\end{align*}
such that
\begin{equation}
\label{eqn:continuous3d}
\begin{cases}
\displaystyle \nu \int_{\Omega}\nabla \mathbf{v}\colon \nabla \mathbf{u}=\int_{\Omega}\mu \nabla \phi \cdot \mathbf{u}-\int_{\Omega}\frac{f'\left(\phi \right)}{2}\left(\mathbf{F}\colon \mathbf{F}-3\right)\nabla \phi \cdot \mathbf{u}-\int_{\Omega}f\left(\phi \right)\mathbf{F}\mathbf{F}^T\colon \nabla \mathbf{u},\\
\displaystyle  \left<\partial_t \mathbf{F}, \boldsymbol{\Sigma} \right>-\int_{\Omega}\mathbf{F}\otimes \mathbf{v}\mathbin{\tensorm} \nabla \boldsymbol{\Sigma}-\int_{\Omega}\left(\nabla {\mathbf{v}}\right)\mathbf{F}\colon \boldsymbol{\Sigma}+\lambda \int_{\Omega}\nabla \left(f({\phi})\mathbf{F}\right)\mathbin{\tensorm} \nabla \boldsymbol{\Sigma}=0,\\
 \displaystyle \int_{\Omega}\partial_t\phi \, q+\int_{\Omega}\left(\mathbf{v} \cdot \nabla \phi\right) q+\int_{\Omega}b\left(\phi \right)\nabla \mu \cdot \nabla q=0,\\
\displaystyle \int_{\Omega}\mu r=\int_{\Omega}\nabla \phi\cdot \nabla r+\int_{\Omega}\psi'\left(\phi\right)r+\int_{\Omega}\frac{f'\left(\phi\right)}{2}\left(\mathbf{F}\colon \mathbf{F}-3\right)r+\delta \int_{\Omega}\partial_t\phi \,r,
\end{cases}
\end{equation} 
a.e. in $t\in (0,T)$ and for all $\mathbf{u}\in H_{0,\diver}^1\left(\Omega;\mathbb{R}^{3}\right)$, $\boldsymbol{\Sigma}\in H^1\left(\Omega;\mathbb{R}^{3\times 3}\right)$, $q,r\in H^1(\Omega)$, satisfying the initial conditions $\mathbf{F}(\cdot,0)=\mathbf{F}_0$ a.e. in $\Omega$ and $\phi(\cdot,0)=\phi_0$ a.e. in $\Omega$..
\end{thm}
In the following, we will define a proper Faedo--Galerkin approximation of \eqref{eqn:systemnoreg}, proving the existence of a discrete solution and studying its convergence to a continuous weak solution in both $2$ and $3$ space dimensions. {\rosso In order to prove the well posedness of the Faedo--Galerkin approximation, we will need to add some further regularizing terms $\gamma \frac{\partial \mathbf{v}}{\partial t}$ and $\theta \frac{\partial \mu}{\partial t}$ to \eqref{eqn:systemnoreg}$_1$ and \eqref{eqn:systemnoreg}$_4$, respectively, with regularization parameters $\gamma, \, \theta$. The limit problem of the approximation will be then investigated in a first step as $\gamma, \, \theta\to 0$, and in a second step as the discretization parameter tends to infinity.}
Before proceeding, we state some preliminary results which will be used in the analysis.
\subsection{Preliminary lemmas}
We state here some Sobolev embedding and interpolation results which will be used in the following calculations. We start by recalling the Gagliardo-Nirenberg inequality {\rosso (see e.g. \cite{brezis,gagliardo,leoni,nirenberg})}.
\begin{lem}
\label{lem:gagliardoniremberg}
Let $\Omega \subset \mathbb{R}^d$, $d=2,3$, be a bounded domain with Lipschitz boundary and $f\in W^{m,r}\cap L^q$, $q\geq 1$, $r\leq \infty$, where $f$ can be a function with scalar, vectorial or tensorial values. For any integer $j$ with $0 \leq j < m$, suppose there is $\alpha \in \mathbb{R}$ such that
\[
j-\frac{d}{p}=\left(m-\frac{d}{r}\right)\alpha+(1-\alpha)\left(-\frac{d}{q}\right), \quad \frac{j}{m}\leq \alpha \leq 1.
\]
Then, there exists a positive constant $C$ depending on $\Omega$, d, m, j, q, r, and $\alpha$ such that
\begin{equation}
\label{eqn:gagliardoniremberg}
||D^jf||_{L^p}\leq C||f||_{W^{m,r}}^{\alpha}||f||_{L^q}^{1-\alpha}.
\end{equation}
\end{lem}
We also state the following interpolation results.
\begin{lem}
\label{lem:interpolation2d3d}
Let $\Omega \subset \mathbb{R}^d$, $d=2,3$, be a bounded domain with Lipschitz boundary and $f \in L^{\infty}(0,T;L^2)\cap L^2(0,T;H^1)$, where $f(\mathbf{x},t)$, with $t\in (0,T)$, $\mathbf{x}\in \Omega$, may be a scalar, a vector or a tensor. 
Then, there exists a positive constant $C$ depending on $\Omega$ and d such that
{\rosso
\begin{equation}
\label{eqn:gninterpolation}
\int_0^T||f||_{L^{2+h}}^{\frac{2(2+h)}{h}}\leq C\int_0^T||f||_{L^{2}}^{\frac{4}{h}}||f||_{H^1}^2, \quad d=2, \; h> 0,
\end{equation}
\begin{equation}
\label{eqn:gninterpolation3d}
\int_0^T||f||_{L^{2+h}}^{\frac{4(2+h)}{3h}}\leq C\int_0^T||f||_{L^{2}}^{\frac{2(4-h)}{3h}}||f||_{H^1}^2, \quad d=3, \; h\in (0,4].
\end{equation}
}%
\end{lem}
We observe that \eqref{eqn:gninterpolation} and \eqref{eqn:gninterpolation3d} are consequences of the Gagliardo--Nirenberg inequality \eqref{eqn:gagliardoniremberg} with $j=0$, $m=1$, $p=2+h$, $r=2$, $q=2$.
\noindent
We moreover recall the following interpolation inequality.
\begin{lem}
\label{lem:interpolationlp}
Let $\Omega \subset \mathbb{R}^d$ be a bounded domain and $f \in L^q$, $q\geq 1$, where $f$ can be a function with scalar, vectorial or tensorial values. Let also $s\leq r \leq q$. Then, there exists a positive constant $C$ depending on $\Omega$ and d such that
\begin{equation}
\label{eqn:lpinterpolation}
\int_{\Omega}|f|^r\leq \left(\int_{\Omega}|f|^s\right)^{\frac{q-r}{q-s}}\left(\int_{\Omega}|f|^q\right)^{\frac{r-s}{q-s}}.
\end{equation}
\end{lem}
We finally state two Agmon type inequalities in two and three space dimensions {\rosso (see e.g. \cite{agmon})} which will be used in the following calculations.
\begin{lem}
\label{lem:agmon2d3d}
Let $\Omega \subset \mathbb{R}^d$, $d=2,3$, be a bounded domain with Lipschitz boundary and $f \in H^2(\Omega)$. 
Then, there exists a positive constant $C$ depending on $\Omega$ and d such that
\begin{equation}
\label{eqn:agmon2d}
||f||_{L^{\infty}(\Omega)}\leq C||f||_{H^1(\Omega)}^{1-\theta}||f|||_{H^2(\Omega)}^{\theta}, \quad \forall \theta \in (0,1], d=2,
\end{equation}
\begin{equation}
\label{eqn:agmon3d}
||f||_{L^{\infty}(\Omega)}\leq C||f||_{H^1(\Omega)}^{\frac{1}{2}}||f|||_{H^2(\Omega)}^{\frac{1}{2}}, \quad d=3.
\end{equation}
\end{lem}
\subsection{Faedo--Galerkin approximation scheme}
{
\rosso Adding the regularizing terms $\gamma \frac{\partial \mathbf{v}}{\partial t}, \theta \frac{\partial \mu}{\partial t}$ to \eqref{eqn:systemnoreg}$_1$ and \eqref{eqn:systemnoreg}$_4$ respectively, we formulate a Faedo--Galerkin approximation scheme for the following system:
\begin{equation}
\label{eqn:system}
\begin{cases}
\gamma \frac{\partial \mathbf{v}}{\partial t}-\nu\Delta \mathbf{v}+\nabla q=\mu \nabla \phi- \left(\frac{1}{2}f'(\phi)(\mathbf{F}\colon \mathbf{F}-d)\right)\nabla \phi+\diver\left(f(\phi)\mathbf{F}\mathbf{F}^T\right),\\
\diver\mathbf{v}=0,\\
\frac{\partial \mathbf{F}}{\partial t}+\left(\mathbf{v}\cdot \nabla\right)\mathbf{F}-(\nabla \mathbf{v})\mathbf{F}-\lambda\Delta \left(f(\phi)\mathbf{F}\right)=0,\\
\frac{\partial \phi}{\partial t}+\mathbf{v}\cdot \nabla \phi-\diver(b(\phi)\nabla \mu)+\theta \frac{\partial \mu}{\partial t}=0,\\
\mu=\psi'(\phi)- \Delta \phi+\frac{1}{2}f'(\phi)(\mathbf{F}\colon \mathbf{F}-d)+\delta \frac{\partial \phi}{\partial t},\\
\end{cases}
\end{equation}
endowed with the same boundary conditions as \eqref{eqn:systemnoreg}, and with initial conditions $\mathbf{v}(\mathbf{x},0)=\mathbf{v}_0$, $\mathbf{F}(\mathbf{x},0)=\mathbf{F}_0$, $\phi(\mathbf{x},0)=\phi_0$, $\mu(\mathbf{x},0)=\mu_0$ for $\mathbf{x}\in \Omega$. As previously explained, the regularizing terms $\gamma \frac{\partial \mathbf{v}}{\partial t}$ and $\theta \frac{\partial \mu}{\partial t}$ are introduced to formulate a well-posed approximation scheme, and imply the specification of artificial initial conditions for $\mathbf{v}$ and $\mu$. In this context, we modify Assumption \textbf{A4} by the following
\begin{itemize}
\item[\bf{A4 Bis}] The initial data have the regularity $\mathbf{v}_0\in L_{\text{div}}^2(\Omega;\mathbb{R}^{d})$, $\mathbf{F}_0\in L^2(\Omega;\mathbb{R}^{d\times d})$, $\phi_0\in H^1(\Omega)$, $\mu_0\in L^2(\Omega)$.
\end{itemize}
We will then obtain existence of a discrete solution for $\gamma, \, \theta \to 0$, while we will study the limit problem for the discretization parameter tending to infinity in the case $\delta \to 0, \, \lambda >0$ in two space dimensions, and in the case $\delta, \, \lambda >0$ in three space dimensions.
}%

In view of obtaining a--priori energy estimates associated to system \eqref{eqn:system}, observe that the equation for $\mathbf{F}$ should be tested by $\partial_{\mathbf{F}}w=f(\phi)\mathbf{F}$ (see \eqref{eqn:energy}), which, when $\phi$ and $\mathbf{F}$ are projected onto proper discrete spaces, does not belong to the discrete space associated to the variable $\mathbf{F}$. Then, we formulate a Galerkin approximation where the third equation of \eqref{eqn:system} is solved in a continuous space. 
\newline
We start by defining the finite dimensional spaces which will be used to formulate the Galerkin ansatz to approximate the solutions of system \eqref{eqn:system}.
Let $\{\boldsymbol{\eta_i}\}_{i\in \mathbb{N}}$ be the eigenfunctions of the Stokes operator with homogeneous Dirichlet boundary conditions, i.e.
\[
P_L(-\Delta)\boldsymbol{\eta_i}=\beta_i \boldsymbol{\eta_i} \quad \text{in} \; \Omega, \quad \boldsymbol{\eta_i}=\mathbf{0} \quad \text{on} \; \partial \Omega,
\]
where $P_L:L^2(\Omega;\mathbb{R}^d)\to L_{\text{div}}^2(\Omega;\mathbb{R}^d)$ is the Leray projection operator, with $0<\beta_0\leq \beta_1 \leq \dots \leq \beta_m\to \infty$. The sequence $\{\boldsymbol{\eta_i}\}_{i\in \mathbb{N}}$ can be chosen as an orthonormal basis in $L_{\diver}^2(\Omega;\mathbb{R}^d)$ and an orthogonal basis in $H_{0,\diver}^1(\Omega;\mathbb{R}^d)$, and, thanks to Assumption \textbf{A0}, we have that $\{\boldsymbol{\eta_i}\}_{i\in \mathbb{N}}\subset C^{\infty}(\bar{\Omega};\mathbb{R}^d)$. We introduce the projection operator
\[
P_m^S:H_{0,\diver}^1(\Omega;\mathbb{R}^d)\to \text{span}\{\boldsymbol{\eta_0},\boldsymbol{\eta_1},\dots,\boldsymbol{\eta_m}\}.
\]
Let moreover $\{\xi_i\}_{i\in \mathbb{N}}$ be the eigenfunctions of the Laplace operator with homogeneous Neumann boundary conditions, i.e.
\[
-\Delta \xi_i=\alpha_i \xi_i \quad \text{in} \; \Omega, \quad \nabla \xi_i \cdot \mathbf{n}=\mathbf{0} \quad \text{on} \; \partial \Omega,
\]
with $0=\alpha_0< \alpha_1 \leq \dots \leq \alpha_m\to \infty$. The sequence $\{\xi_i\}_{i\in \mathbb{N}}$ can be chosen as an orthonormal basis in $L^2(\Omega)$ and an orthogonal basis in $H^1(\Omega)$, and, thanks to Assumption \textbf{A0}, $\{\xi_i\}_{i\in \mathbb{N}}\subset C^{\infty}(\bar{\Omega})$. Without loss of generality, we assume $\alpha_0=0$. We introduce the projection operator
\[
P_m^L:H^1(\Omega)\to \text{span}\{\xi_0,\xi_1,\dots,\xi_m\}.
\]
We make the Galerkin ansatz $\mathbf{v}_m=\sum_i d_i^m(t)\boldsymbol{\eta_i}(\mathbf{x})$, $\phi_m=\sum_i a_i^m(t)\xi_i(\mathbf{x})$, $\mu_m=\sum_i b_i^m(t)\xi_i(\mathbf{x})$ to approximate the solutions $\mathbf{v}, \phi, \mu$ of system \eqref{eqn:system}. Hence, we project the equation for $\mathbf{v}_m$ onto \\
$\text{span}\left\{\boldsymbol{\eta}_0,\boldsymbol{\eta}_1,\dots,\boldsymbol{\eta}_m\right\}$ and the equations for $\phi_m$ and $\mu_m$ onto $\text{span}\left\{\xi_0,\xi_1, \dots, \xi_m\right\}$, and formulate the following Galerkin approximation of system \eqref{eqn:system}: 
\begin{equation}
\label{eqn:galerkin}
\begin{cases}
 \displaystyle \gamma\int_{\Omega}\partial_t\mathbf{v}_m\cdot \boldsymbol{\eta}_i+\nu \int_{\Omega}\nabla \mathbf{v}_m\colon \nabla \boldsymbol{\eta}_i =\int_{\Omega}\mu_m \nabla \phi_m \cdot \boldsymbol{\eta}_i-\int_{\Omega}\frac{f'\left(\phi_m\right)}{2}\left(\mathbf{F}_m\colon \mathbf{F}_m-d\right)\nabla \phi_m\cdot \boldsymbol{\eta}_i\\
 \displaystyle -\int_{\Omega}f\left(\phi_m\right)\mathbf{F}_m\mathbf{F}_m^T\colon \nabla \boldsymbol{\eta}_i,\\
\displaystyle  \left<\partial_t \mathbf{F}_m, \boldsymbol{\Sigma} \right>+\int_{\Omega}\left({\mathbf{v}}_m \cdot \nabla \right)\mathbf{F}_m\colon \boldsymbol{\Sigma}-\int_{\Omega}\left(\nabla {\mathbf{v}}_m\right)\mathbf{F}_m\colon \boldsymbol{\Sigma}+\lambda \int_{\Omega}\nabla \left(f({\phi}_m)\mathbf{F}_m\right)\mathbin{\tensorm} \nabla \boldsymbol{\Sigma}=0,\\
 \displaystyle \int_{\Omega}\partial_t\phi_m \xi_i+\int_{\Omega}\left(\mathbf{v}_m \cdot \nabla \phi_m\right) \xi_i+\int_{\Omega}b\left(\phi_m\right)\nabla \mu_m \cdot \nabla \xi_i +\theta \int_{\Omega}\partial_t\mu_m \xi_i=0,\\
\displaystyle \int_{\Omega}\mu_m \xi_i=\int_{\Omega}\nabla \phi_m\cdot \nabla \xi_i+\int_{\Omega}\psi'\left(\phi_m\right)\xi_i+\int_{\Omega}\frac{f'\left(\phi_m\right)}{2}\left(\mathbf{F}_m\colon \mathbf{F}_m-d\right)\xi_i+\delta \int_{\Omega}\partial_t\phi_m \xi_i,
\end{cases}
\end{equation} 
{\rosso a.e. in $[0,t)$}, with $0<t\leq T$, for $i=0, \dots,m$ and $\forall \, \boldsymbol{\Sigma} \in H^1(\Omega;\mathbb{R}^{d\times d})$, with $\mathbf{F}_m$ belonging to a proper continuous space, which will be defined in the following. Here, we take the initial conditions $\phi_m(\mathbf{x},0)=P_m^L(\phi_0), \mu_m(\mathbf{x},0)=P_m^L(\mu_0), \mathbf{v}_m(\mathbf{x},0)=P_m^S(\mathbf{v}_0), \mathbf{F}_m(\mathbf{x},0)=\mathbf{F}_0$.

In order to prove the existence of a solution of \eqref{eqn:galerkin} and its convergence to a continuous weak solution of system \eqref{eqn:system}, following \cite{benesova}, we employ a fixed point strategy, which in our case can be developed in the following steps:
\begin{itemize}
\item[\bf{Step 1}] In the second equation of \eqref{eqn:galerkin}, we take fixed values $\bar{\phi}_m$ and $\bar{\mathbf{v}}_m$ of the variables $\phi_m$ and $\mathbf{v}_m$ in proper convex subsets of the corresponding Galerkin spaces, and solve for $\mathbf{F}_m$ in a proper continuous space;
\item[\bf{Step 2}] Given the solution $\mathbf{F}_m$ obtained in the previous step, we solve for the variables $\phi_m$, $\mu_m$ and $\mathbf{v}_m$, hence a fixed point map is defined from $(\bar{\phi}_m,\bar{\mathbf{v}}_m)$ to $(\phi_m,\mathbf{v}_m)$;
\item[\bf{Step 3}] We prove the existence of a fixed point of the map defined in the previous step, {\rosso which yields a solution of the original approximating problem} \eqref{eqn:galerkin};
\item[\bf{Step 4}] We obtain a--priori estimates, uniformly independent of the discretization parameter, for the discrete solutions, and study the limit problem as the discretization parameter tends to zero.
\end{itemize}
\subsection{Step 1: Solution for $\mathbf{F}$ at the Galerkin level}
We define the following spaces, for $t>0$,
\begin{align}
\label{eqn:vmpm}
& V_{m}(t):=\left\{\mathbf{v}=\sum_{i=0}^m d_i^m(s)\boldsymbol{\eta_i}(\mathbf{x}): \sup_{s\in[0,t]}\left(\sum_{i=0}^m\left|d_i^m(s)\right|^2\right)^{\frac{1}{2}}\leq L_v:=\left|\left|P_m^S(\mathbf{v}_0)\right|\right|_{L^2(\Omega;\mathbb{R}^d)}+1\right\},\\
&\notag \Phi_{m}(t):=\left\{\phi=\sum_{i=0}^m a_i^m(s)\xi(\mathbf{x}): \sup_{s\in[0,t]}\left(\sum_{i=0}^m\left|a_i^m(s)\right|^2\right)^{\frac{1}{2}}\leq L_{\phi}:=\left|\left|P_m^L(\phi_0)\right|\right|_{L^2(\Omega)}+1\right\},
\end{align}
which are closed and convex subsets of the spaces $C^0([0,t];\text{span}\{\boldsymbol{\eta_0},\boldsymbol{\eta_1},\dots,\boldsymbol{\eta_m}\})$ and \\
$C^0([0,t];\text{span}\{\xi_0,\xi_1,\dots,\xi_m\})$ respectively. 
\newline
Given fixed $\bar{\mathbf{v}}_m\in V_m(t_0)$ and $\bar{\phi}_m\in \Phi_m(t_0)$, with $0<t_0\leq T$, we can define a proper weak solution of the third equation of system \eqref{eqn:system} through the following theorem.
\begin{thm}
\label{thm:fgalerkin}
Given $\bar{\mathbf{v}}_m\in V_m(t_0)$ and $\bar{\phi}_m\in \Phi_m(t_0)$, there exist a $t_1$, with $0<t_1\leq t_0$, and a unique $\mathbf{F}_m$, with
\begin{equation}
\label{eqn:fspace}
\mathbf{F}_m\in L^{\infty}\left(0,t_1;L^2(\Omega;\mathbb{R}^{d\times d})\right) \cap L^{2}\left(0,t_1;H^1(\Omega;\mathbb{R}^{d\times d})\right) \cap H^1\left(0,t_1;\left(H^1(\Omega;\mathbb{R}^{d\times d})\right)'\right),
\end{equation}
and
\begin{equation}
\label{eqn:fnorm}
||\mathbf{F}_m||_{L^{\infty}\left(0,t_1;L^2(\Omega;\mathbb{R}^{d\times d})\right)}+ ||\mathbf{F}_m||_{L^{2}\left(0,t_1;H^1(\Omega;\mathbb{R}^{d\times d})\right)}+||\mathbf{F}_m||_{H^1\left(0,t_1;\left(H^1(\Omega;\mathbb{R}^{d\times d})\right)'\right)}\leq C,
\end{equation}
where $C$ is independent of the regularization parameters $\gamma, \theta, \delta$, satisfying the weak formulation
\begin{equation}
\label{eqn:fweak}
\left<\partial_t \mathbf{F}_m, \boldsymbol{\Sigma} \right>+\int_{\Omega}\left(\bar{\mathbf{v}}_m \cdot \nabla \right)\mathbf{F}_m\colon \boldsymbol{\Sigma}-\int_{\Omega}\left(\nabla \bar{\mathbf{v}}_m\right)\mathbf{F}_m\colon \boldsymbol{\Sigma}+\lambda \int_{\Omega}\nabla \left(f(\bar{\phi}_m)\mathbf{F}_m\right)\mathbin{\tensorm} \nabla \boldsymbol{\Sigma}=0,
\end{equation}
for all $\boldsymbol{\Sigma} \in H^1(\Omega;\mathbb{R}^{d\times d})$ and a.e. $t \in [0,t_1)$, with initial condition $\mathbf{F}_m(\mathbf{x},0)=\mathbf{F}_0(\mathbf{x})$. Here, $t_1$ depends on the discretization parameter $m$, on $L_v$, $L_{\phi}$, and on the initial conditions $\mathbf{v}_0,\phi_0$.
\end{thm}
\begin{pf}
Let $\{\boldsymbol{\Sigma}_i\}_{i\in \mathbb{N}}$ be the eigenfunctions of the Laplace operator with homogeneous Neumann boundary conditions, i.e.
\[
-\Delta \boldsymbol{\Sigma}_i=\gamma_i \boldsymbol{\Sigma}_i \quad \text{in} \; \Omega, \quad \nabla \boldsymbol{\Sigma}_i  \mathbf{n}=\mathbf{0} \quad \text{on} \; \partial \Omega,
\]
with $0=\gamma_0< \gamma_1 \leq \dots \leq \gamma_m\to \infty$. The sequence $\{\boldsymbol{\Sigma}_i\}_{i\in \mathbb{N}}$ can be chosen as an orthonormal basis in $L^2(\Omega;\mathbb{R}^{d\times d}))$ and an orthogonal basis in $H^1(\Omega;\mathbb{R}^{d\times d}))$, and, thanks to Assumptions \textbf{A0} and \textbf{A1}, $\{\boldsymbol{\Sigma}_i\}_{i\in \mathbb{N}}\subset C^{\infty}(\bar{\Omega};\mathbb{R}^{d\times d})$. Without loss of generality, we assume $\gamma_0=0$. We introduce the projection operator
\[
P_n^{\Sigma}:H^1(\Omega;\mathbb{R}^{d\times d})\to \text{span}\{\boldsymbol{\Sigma}_0,\boldsymbol{\Sigma}_1,\dots,\boldsymbol{\Sigma}_n\}.
\]
We make the Galerkin ansatz
\[
\mathbf{F}_n=\sum_i f_i^n(t)\boldsymbol{\Sigma_i}(\mathbf{x})
\]
for the solution of the third equation of \eqref{eqn:system}, and project the equation onto $\text{span}\{\boldsymbol{\Sigma_0},\boldsymbol{\Sigma_1},\dots,\boldsymbol{\Sigma_n}\}$, i.e. 
\begin{equation}
\label{eqn:fgalerkin}
\int_{\Omega}\partial_t \mathbf{F}_n \colon \boldsymbol{\Sigma}_i+\int_{\Omega}\left(\bar{\mathbf{v}}_m \cdot \nabla \right)\mathbf{F}_n\colon \boldsymbol{\Sigma}_i-\int_{\Omega}\left(\nabla \bar{\mathbf{v}}_m\right)\mathbf{F}_n\colon \boldsymbol{\Sigma}_i+\lambda \int_{\Omega}\nabla \left(f(\bar{\phi}_m)\mathbf{F}_n\right)\mathbin{\tensorm} \nabla \boldsymbol{\Sigma}_i=0, \quad i=0, \dots, n,
\end{equation}
with initial condition $\mathbf{F}_n(\mathbf{x},0)=\sum_i \left(\mathbf{F}_0,\boldsymbol{\Sigma_i}\right)\boldsymbol{\Sigma_i}(\mathbf{x})$, where $\mathbf{F}_0(\mathbf{x})=\mathbf{F}(\mathbf{x},0)$. System \eqref{eqn:fgalerkin} defines a collection of initial value problems for the ODEs
\begin{equation}
\label{eqn:fode}
\begin{cases}
\displaystyle \frac{d}{dt}f_i^n=\sum_l \biggl(-\int_{\Omega}\left(\bar{\mathbf{v}}_m \cdot \nabla \right)\boldsymbol{\Sigma}_l\colon \boldsymbol{\Sigma}_i+\int_{\Omega}\left(\nabla \bar{\mathbf{v}}_m\right)\boldsymbol{\Sigma}_l\colon \boldsymbol{\Sigma}_i-\lambda \int_{\Omega}f(\bar{\phi}_m)\nabla \boldsymbol{\Sigma}_l \mathbin{\tensorm} \nabla \boldsymbol{\Sigma}_i\\
\displaystyle \qquad \qquad \qquad -\lambda \int_{\Omega}f'(\bar{\phi}_m)\boldsymbol{\Sigma}_l \otimes \nabla \bar{\phi}_m \mathbin{\tensorm} \nabla \boldsymbol{\Sigma}_i\biggr)f_l^n,\\
\displaystyle f_i^n(0)=\left(\mathbf{F}_0,\boldsymbol{\Sigma_i}\right), \quad i=0, \dots, n.
\end{cases}
\end{equation} 
Thanks to the regularity of the fields $\bar{\mathbf{v}}_m\in V_m(t_0)$, $\bar{\phi}_m\in \Phi_m(t_0)$ and $\boldsymbol{\Sigma}_i\in C^{\infty}(\bar{\Omega};\mathbb{R}^{d\times d})$, we can apply the Peano existence theorem to infer that there exists a sufficiently small $t_{*}$ with $0<t_{*}\leq t_1$ and a unique local solution $f_i^n\in C^1([0,t_{*}))$ of \eqref{eqn:fode}, for $i=0, \dots, n$. 
\begin{rem}
\label{rem:fic}
We observe that $t_{*}$ depends on $n$, on the initial condition $\mathbf{F}_0$, on the $L^{\infty}(0,t_0;H_0^1(\Omega;\mathbb{R}^d))$ norm of $\bar{\mathbf{v}}_m$ and  on the $L^{\infty}(0,t_0;H^1(\Omega))$ norm of $\bar{\phi}_m$, which in turn, given the definitions \eqref{eqn:vmpm}, depend on $L_v$ and $L_{\phi}$, hence on the initial conditions $\mathbf{v}_0$, $\phi_0$ and on $m$.
\end{rem}
We now derive a--priori estimates uniform in $n$ for the Galerkin solution $\mathbf{F}_n$, in order to pass to the limit as $n\to \infty$ and to prove \eqref{eqn:fweak}.
Combining the $n+1$ equations \eqref{eqn:fgalerkin}, for $i=0, \dots, n$, we get
\begin{align}
\label{eqn:fweak2}
&\int_{\Omega}\partial_t \mathbf{F}_n \colon \boldsymbol{\Sigma}+\int_{\Omega}\left(\bar{\mathbf{v}}_m \cdot \nabla \right)\mathbf{F}_n\colon \boldsymbol{\Sigma}-\int_{\Omega}\left(\nabla \bar{\mathbf{v}}_m\right)\mathbf{F}_n\colon \boldsymbol{\Sigma}+\lambda \int_{\Omega}f(\bar{\phi}_m)\nabla \mathbf{F}_n\mathbin{\tensorm} \nabla \boldsymbol{\Sigma}\\
&\notag +\lambda \int_{\Omega}f'(\bar{\phi}_m)\mathbf{F}_n \otimes \nabla \bar{\phi}_m \mathbin{\tensorm} \nabla \boldsymbol{\Sigma}=0, 
\end{align}
with $\boldsymbol{\Sigma}\in \text{span}\{\boldsymbol{\Sigma}_0,\boldsymbol{\Sigma}_1,\dots,\boldsymbol{\Sigma}_n\}$. Taking $\boldsymbol{\Sigma}=\mathbf{F}_n$, considering that 
\[\int_{\Omega}\bar{\mathbf{v}}_m \cdot \nabla \left(\frac{\left|\mathbf{F}_n\right|^2}{2}\right)=-\int_{\Omega}\diver\bar{\mathbf{v}}_m\frac{\left|\mathbf{F}_n\right|^2}{2}+\int_{\partial \Omega}\frac{\left|\mathbf{F}_n\right|^2}{2}\bar{\mathbf{v}}_m \cdot \mathbf{n}=0,
\]
using the Cauchy--Schwarz and the Young inequalities and Assumption \textbf{A3}, we obtain
\begin{align*}
&\frac{1}{2}\frac{d}{dt}\int_{\Omega}\left|\mathbf{F}_n\right|^2+\lambda \int_{\Omega}f(\bar{\phi}_m)\left|\nabla \mathbf{F}_n\right|^2=\int_{\Omega}\left(\nabla \bar{\mathbf{v}}_m\right)\colon \mathbf{F}_n\mathbf{F}_n^T-\lambda \int_{\Omega}f'(\bar{\phi}_m)\mathbf{F}_n \otimes \nabla \bar{\phi}_m \mathbin{\tensorm} \nabla \mathbf{F}_n \\
&\leq \left|\left|\nabla \bar{\mathbf{v}}_m\right|\right|_{L^{\infty}(\Omega;\mathbb{R}^{d\times d})}\int_{\Omega}\left|\mathbf{F}_n\right|^2+\lambda f_1\left|\left|\nabla \bar{\phi}_m\right|\right|_{L^{\infty}(\Omega;\mathbb{R}^{d})}\left|\left|\mathbf{F}_n\right|\right|_{L^2(\Omega;\mathbb{R}^{d\times d})} \left|\left|\nabla \mathbf{F}_n\right|\right|_{L^2(\Omega;\mathbb{R}^{d\times d\times d})} \\
&\leq \frac{\lambda}{2}f_{\text{min}}\int_{\Omega}\left|\nabla \mathbf{F}_n\right|^2+\left(\left|\left|\nabla \bar{\mathbf{v}}_m\right|\right|_{L^{\infty}(\Omega;\mathbb{R}^{d\times d})}+\frac{\lambda f_1^2}{2f_{\text{min}}}\left|\left|\nabla \bar{\phi}_m\right|\right|_{L^{\infty}(\Omega;\mathbb{R}^{d})}^2\right)\int_{\Omega}\left|\mathbf{F}_n\right|^2.
\end{align*}
Integrating in time over the interval $(0,t_{*})$, using again Assumption \textbf{A3} and the fact that $\bar{\mathbf{v}}_m\in V_m(t_0)$, $\bar{\phi}_m\in \Phi_m(t_0)$, we get
\begin{equation}
\label{eqn:fapriori}
\frac{1}{2}\int_{\Omega}\left|\mathbf{F}_n\right|^2+\frac{\lambda}{2}f_{\text{min}}\int_0^{t_{*}}\int_{\Omega}\left|\nabla \mathbf{F}_n\right|^2\leq \frac{1}{2}\int_{\Omega}\left|P_n^{\Sigma}\left(\mathbf{F}_0\right)\right|^2+C(L_v,L_{\phi},m)\int_0^{t_{*}}\int_{\Omega}\left|\mathbf{F}_n\right|^2.
\end{equation}
Using a Gronwall argument, we get from \eqref{eqn:fapriori} the uniform in $n$ estimate
{\rosso
\begin{equation}
\label{eqn:fapriori2}
\left|\left|\mathbf{F}_n\right|\right|_{C^{0}\left([0,t_{*}];L^2\left(\Omega;\mathbb{R}^{d\times d}\right)\right)}^2+\left|\left|\nabla \mathbf{F}_n\right|\right|_{L^{2}\left(0,t_{*};L^2\left(\Omega;\mathbb{R}^{d\times d\times d}\right)\right)}^2\leq C\left(L_v,L_{\phi},m,\mathbf{F}_0,t_{*}\right).
\end{equation}
}%
Thanks to the a--priori estimate \eqref{eqn:fapriori2} and using {\rosso an extension argument for the local solution, it turns out that the maximal solution of \eqref{eqn:fode} lives in some interval $[0,t_1)$, with $0<t_1\leq t_0$ depending only on $L_v$,  $L_{\phi}$ and on $m$}. In order to derive a uniform in $n$ bound for $\partial_t \mathbf{F}_n$, we multiply \eqref{eqn:fweak2} by a time function $\zeta \in L^2(0,t_1)$, integrate in time from $0$ to $t_1$ and take $\boldsymbol{\Sigma}=P_n^{\Sigma}\left(\boldsymbol{\Pi}\right)$, with $\boldsymbol{\Pi}\in H^1(\Omega;\mathbb{R}^{d\times d})$, obtaining, using the Cauchy--Schwarz inequality, 
\begin{align}
\label{eqn:faprioridual}
&{\rosso \int_0^{t_1}\int_{\Omega}\partial_t \mathbf{F}_n \colon \boldsymbol{\Pi}\zeta=\int_0^{t_1}\int_{\Omega}\partial_t \mathbf{F}_n \colon P_n^{\Sigma}\left(\boldsymbol{\Pi}\right)\zeta}\\
&\notag =-\int_0^{t_1}\int_{\Omega}\left(\bar{\mathbf{v}}_m \cdot \nabla \right)\mathbf{F}_n\colon P_n^{\Sigma}\left(\boldsymbol{\Pi}\right)\zeta+\int_0^{t_1}\int_{\Omega}\left(\nabla \bar{\mathbf{v}}_m\right)\mathbf{F}_n\colon P_n^{\Sigma}\left(\boldsymbol{\Pi}\right)\zeta\\
& \notag -\lambda \int_0^{t_1}\int_{\Omega}f(\bar{\phi}_m)\nabla \mathbf{F}_n\mathbin{\tensorm} \nabla \left(P_n^{\Sigma}\left(\boldsymbol{\Pi}\right)\right)\zeta-\lambda \int_0^{t_1}\int_{\Omega}f'(\bar{\phi}_m)\mathbf{F}_n \otimes \nabla \bar{\phi}_m \mathbin{\tensorm} \nabla \left(P_n^{\Sigma}\left(\boldsymbol{\Pi}\right)\right)\zeta \\
& \notag \leq \int_0^{t_1}\biggl(\left|\left|\bar{\mathbf{v}}_m\right|\right|_{L^{\infty}(\Omega;\mathbb{R}^{d})}\left|\left|\nabla \mathbf{F}_n\right|\right|_{L^2(\Omega;\mathbb{R}^{d\times d\times d})}+\left|\left|\nabla \bar{\mathbf{v}}_m\right|\right|_{L^{\infty}(\Omega;\mathbb{R}^{d\times d})}\left|\left| \mathbf{F}_n\right|\right|_{L^2(\Omega;\mathbb{R}^{d\times d})}\biggr)|\zeta|\left|\left| P_n^{\Sigma}\left(\boldsymbol{\Pi}\right)\right|\right|_{L^2(\Omega;\mathbb{R}^{d\times d})}\\
& \notag +\int_0^{t_1}\biggl(\lambda \left|\left|\nabla \mathbf{F}_n\right|\right|_{L^2(\Omega;\mathbb{R}^{d\times d\times d})}+\lambda f_1\left|\left|\nabla \bar{\phi}_m\right|\right|_{L^{\infty}(\Omega;\mathbb{R}^{d})}\left|\left| \mathbf{F}_n\right|\right|_{L^2(\Omega;\mathbb{R}^{d\times d})}\biggr)|\zeta|\left|\left| \nabla P_n^{\Sigma}\left(\boldsymbol{\Pi}\right)\right|\right|_{L^2(\Omega;\mathbb{R}^{d\times d \times d})} \\
& \notag {\rosso \leq C(L_{v},L_{\phi},m)\left|\left| \mathbf{F}_n\right|\right|_{L^2(0,t_1;H^1\left(\Omega;\mathbb{R}^{d\times d})\right)}||\zeta||_{L^2(0,t_1)}\left|\left| P_n^{\Sigma}\left(\boldsymbol{\Pi}\right)\right|\right|_{H^1(\Omega;\mathbb{R}^{d\times d})}}\\
& \notag {\rosso \leq C(L_{v},L_{\phi},m,\mathbf{F}_0)||\zeta||_{L^2(0,t_1)}\left|\left|\boldsymbol{\Pi}\right|\right|_{H^1(\Omega;\mathbb{R}^{d\times d})}}.
\end{align}
Using \eqref{eqn:fapriori2}, extended to the time interval $[0,t_1)$, in \eqref{eqn:faprioridual}, we obtain the bound 
\begin{equation}
\label{eqn:faprioridual2}
\left|\left|\partial_t \mathbf{F}_n\right|\right|_{L^{2}\left(0,t_1;\left(H^1\left(\Omega;\mathbb{R}^{d\times d}\right)\right)'\right)}\leq C\left(L_v,L_{\phi},m,\mathbf{F}_0\right).
\end{equation}
The bounds \eqref{eqn:fapriori2}, extended to the time interval $[0,t_1)$, and \eqref{eqn:faprioridual2} prove \eqref{eqn:fspace}. Hence, applying the Banach--Alaoglu theorem, there exists a subsequence of $\left(\mathbf{F}_n\right)_{n\in \mathbb{N}}$, which we still call $\left(\mathbf{F}_n\right)_{n\in \mathbb{N}}$, such that
\begin{align}
\label{eqn:fnconv}
&\partial_t \mathbf{F}_n \rightharpoonup \partial_t \mathbf{F}_m \quad \text{in} \quad L^{2}\left(0,t_1;\left(H^1\left(\Omega;\mathbb{R}^{d\times d}\right)\right)'\right);\\
& \notag \mathbf{F}_n \overset{\ast}{\rightharpoonup} \mathbf{F}_m \quad \text{in} \quad L^{\infty}\left(0,t_1;L^2\left(\Omega;\mathbb{R}^{d\times d}\right)\right);\\
& \notag \nabla \mathbf{F}_n \rightharpoonup \nabla \mathbf{F}_m \quad \text{in} \quad L^{2}\left(0,t_1;L^2\left(\Omega;\mathbb{R}^{d\times d\times d}\right)\right).
\end{align}
Passing to the limit as $n\to \infty$ in \eqref{eqn:fweak2} and using the aforementioned convergence properties, on noting that \eqref{eqn:fweak2} is linear in $\mathbf{F}_n$, {\rosso and using a density argument on the test elements $\boldsymbol{\Sigma}$}, we get that the limit point $\mathbf{F}_m$ is the unique solution of \eqref{eqn:fweak}. Finally, from the facts that $\mathbf{F}_n(\mathbf{x},0)\to \mathbf{F}_0$ in $L^2(\Omega;\mathbb{R}^{d \times d})$ and that {\rosso $\mathbf{F}_n \to \mathbf{F}_m$ in $C_w^0([0,t_1];L^2(\Omega;\mathbb{R}^{d \times d}))$}, which follows from \eqref{eqn:fnconv} and e.g. \cite[Ch. 3, Lemma 1.4]{temam}, we get that $\mathbf{F}_m(\mathbf{x},0)=\mathbf{F}_0$.
\end{pf}
\subsection{Step $2$: Definition of the fixed point map}
Given $(\bar{\mathbf{v}}_m,\bar{\phi}_m)\in V_m(t_0)\times \Phi_m(t_0)$ and $\mathbf{F}_m$ as in \eqref{eqn:fspace} a solution of \eqref{eqn:fweak} in the time interval $[0,t_1)$, with $0<t_1\leq t_0$, the equations for $\mathbf{v}_m$, $\phi_m$ and $\mu_m$ in \eqref{eqn:galerkin} define an initial value problem for a system of coupled ODEs,
\begin{equation}
\label{eqn:odegalerkin}
\begin{cases}
 \displaystyle \gamma\frac{d}{dt}d_i^m=-\nu \beta_id_i^m+\sum_{l,s}\left[\int_{\Omega}\xi_l \nabla \xi_s \cdot \boldsymbol{\eta}_i\right]b_l^ma_s^m-\\
\displaystyle \sum_{l}\left[\int_{\Omega}\frac{f'\left(\sum_{s} a_s^m\xi_s\right)}{2}\left(\mathbf{F}_m\colon \mathbf{F}_m-d\right)\nabla \xi_l\cdot \boldsymbol{\eta}_i\right]a_l^m -\int_{\Omega}f\left(\sum_{s} a_s^m\xi_s\right)\mathbf{F}_m\mathbf{F}_m^T\colon \nabla \boldsymbol{\eta}_i,\\
 \displaystyle \frac{d}{dt}a_i^m=-\sum_{l,s}\left[\int_{\Omega}\left(\boldsymbol{\eta}_l \cdot \nabla \xi_s\right) \cdot \xi_i\right]d_l^ma_s^m-
\sum_{l}\left[\int_{\Omega}b\left(\sum_{s} a_s^m\xi_s\right)\nabla \xi_l \cdot \nabla \xi_i \right]b_l^m-\theta \frac{d}{dt}b_i^m\\
\displaystyle b_i^m=\alpha_ia_i^m+\int_{\Omega}\psi'\left(\sum_{s} a_s^m\xi_s\right)\xi_i+\int_{\Omega}\frac{f'\left(\sum_{s} a_s^m\xi_s\right)}{2}\left(\mathbf{F}_m\colon \mathbf{F}_m-d\right)\xi_i+\delta \frac{d}{dt}a_i^m\\
\displaystyle a_i^m(0)=\left(\phi_0,\xi_i\right), \; b_i^m(0)=\left(\mu_0,\xi_i\right), \; d_i^m(0)=\left(\mathbf{v}_0,\boldsymbol{\eta}_i\right), \quad i=0, \dots, m.
\end{cases}
\end{equation} 
We can substitute $\frac{d}{dt}a_i^m$ from the third equation in the second equation of system \eqref{eqn:odegalerkin}, and obtain a system of coupled first order ODEs, which, upon introducing the vectors of time coefficients $\mathbf{a}=\left(a_0^m,\dots,a_m^m\right)\in \mathbb{R}^{m+1}$, $\mathbf{b}=\left(b_0^m,\dots,b_m^m\right)\in \mathbb{R}^{m+1}$, $\mathbf{d}=\left(d_0^m,\dots,d_m^m\right)\in \mathbb{R}^{m+1}$, takes the form  
\begin{equation}
\label{eqn:odegalerkin2}
\begin{cases}
 \displaystyle \gamma\frac{d}{dt}\mathbf{d}=-\nu \text{diag}(\boldsymbol{\beta})\mathbf{d}+\mathbf{D}(\mathbf{b}\otimes \mathbf{a})-\mathbf{E}\left(\mathbf{a},\mathbf{F}_m\right)\mathbf{a}-\mathbf{p}\left(\mathbf{a},\mathbf{F}_m\right),\\
 \displaystyle \delta \frac{d}{dt}\mathbf{a}=\mathbf{b}-\text{diag}(\boldsymbol{\alpha})\mathbf{a}-\mathbf{r}(\mathbf{a})-\mathbf{q}\left(\mathbf{a},\mathbf{F}_m\right),\\
\displaystyle \theta \frac{d}{dt}\mathbf{b}=-\frac{1}{\delta}\left(\mathbf{b}-\text{diag}(\boldsymbol{\alpha})\mathbf{a}-\mathbf{r}(\mathbf{a})-\mathbf{q}\left(\mathbf{a},\mathbf{F}_m\right)\right)-\mathbf{B}(\mathbf{d}\otimes \mathbf{a})-\mathbf{C}(\mathbf{a})\mathbf{b},\\
\displaystyle \mathbf{a}(0)=\mathbf{a}_0, \mathbf{b}(0)=\mathbf{b}_0, \mathbf{d}(0)=\mathbf{d}_0,
\end{cases}
\end{equation} 
where we defined the third order and second order tensors and vectors as (in components, for $i,j,k=0, \dots, m$)
\begin{align}
\label{eqn:tensors}
& \mathbf{D}_{ijk}:=\int_{\Omega}\xi_k \nabla \xi_j \cdot \boldsymbol{\eta}_i; \quad \mathbf{B}_{ijk}:=\int_{\Omega}\left(\boldsymbol{\eta}_k \cdot \nabla \xi_j\right) \cdot \xi_i;\\
& \notag \mathbf{E}\left(\mathbf{a},\mathbf{F}_m\right)_{ij}:=\int_{\Omega}\frac{f'\left(\sum_{s} a_s^m\xi_s\right)}{2}\left(\mathbf{F}_m\colon \mathbf{F}_m-d\right)\nabla \xi_j\cdot \boldsymbol{\eta}_i; \; \; \mathbf{C}(\mathbf{a})_{ij}:=\int_{\Omega}b\left(\sum_{s} a_s^m\xi_s\right)\nabla \xi_j \cdot \nabla \xi_i;\\
& \notag \mathbf{p}\left(\mathbf{a},\mathbf{F}_m\right)_i:=\int_{\Omega}f\left(\sum_{s} a_s^m\xi_s\right)\mathbf{F}_m\mathbf{F}_m^T\colon \nabla \boldsymbol{\eta}_i; \; \; \mathbf{r}(\mathbf{a})_i:=\int_{\Omega}\psi'\left(\sum_{s} a_s^m\xi_s\right)\xi_i;\\
& \notag \mathbf{q}\left(\mathbf{a},\mathbf{F}_m\right)_i:=\int_{\Omega}\frac{f'\left(\sum_{s} a_s^m\xi_s\right)}{2}\left(\mathbf{F}_m\colon \mathbf{F}_m-d\right)\xi_i; \; \, \mathbf{a}_{0i}:=\left(\phi_0,\xi_i\right); \; \mathbf{b}_{0i}:=\left(\mu_0,\xi_i\right); \; \mathbf{d}_{0i}:=\left(\mathbf{v}_0,\boldsymbol{\eta}_i\right).
\end{align}
Due to the Assumptions \textbf{A1, A2, A3} on the regularity of the functions $b,\psi,f$, the regularity in space of the functions $\xi_i, \boldsymbol{\eta}_i$ and the regularity result $\mathbf{F}_m\in L^{\infty}(0,t_1;L^2(\Omega;\mathbb{R}^{d\times d}))$, the right hand side of \eqref{eqn:odegalerkin2} depends continuously on the variables $\mathbf{a}, \mathbf{b}, \mathbf{d}$, with bounded coefficients in $L^{\infty}(0,t_1)$. Hence, we can apply the Carath\'eodory's existence theorem to obtain the existence of a sufficiently small value $t_2$, with $0<t_2\leq t_1$, and of a unique Lipschitz continuous solution $\mathbf{a}, \mathbf{b}, \mathbf{d}\in W^{1,\infty}(0,t_2)$.
\begin{rem}
\label{rem:gic}
The Lipschitz constant of the solution of the system \eqref{eqn:odegalerkin2}, as well as the parameter $t_2$, depend on $m$, on the initial conditions $\mathbf{a}_0, \mathbf{b}_0, \mathbf{d}_0$ and on the $L^{\infty}(0,t_2)$ norm of the coefficients, hence on the  $L^{\infty}(0,t_2;L^2(\Omega;\mathbb{R}^{d\times d}))$ norm of $\mathbf{F}_m$. The latter norm depends, thanks to \eqref{eqn:fapriori2} and \eqref{eqn:fspace}, on $L_v$, $L_{\phi}$, $m$ and $\mathbf{F}_0$.
\end{rem}
Given the solution $\mathbf{F}_m$ of \eqref{eqn:fweak}, with fixed $\left(\bar{\mathbf{v}}_m,\bar{\phi}_m\right)\in V_m(t_0)\times \Phi(t_0)$, the solution of system \eqref{eqn:odegalerkin2} defines a map
\begin{align}
\label{eqn:map}
&\mathcal{L}: V_m(t_0)\times \Phi(t_0) \to W^{1,\infty}(0,t_2;\text{span}\left\{\boldsymbol{\eta}_0, \boldsymbol{\eta}_1, \dots, \boldsymbol{\eta}_m\right\})\times W^{1,\infty}(0,t_2;\text{span}\left\{\xi_0, \xi_1, \dots, \xi_m\right\}), \\
& \notag \mathcal{L}\left(\bar{\mathbf{v}}_m,\bar{\phi}_m\right)=\left(\mathbf{v}_m,\phi_m\right).
\end{align} 
In the following, we will prove the existence of a fixed point of $\mathcal{L}$, which is a solution of system \eqref{eqn:galerkin}, by applying the Schauder fixed point theorem. 
\subsection{Step $3$: Existence of a solution of the Faedo--Galerkin approximation}
We prove the following lemma.
\begin{lem}
\label{lem:fixedpoint}
There exists a small enough value $t_3$, with $0<t_3\leq t_2$, such that $(\mathbf{v}_m,\phi_m)$, solution of system \eqref{eqn:odegalerkin2}, belongs to the space $V_m(t_3)\times \Phi_m(t_3)$. Hence, 
\[\mathcal{L}: V_m(t_3)\times \Phi_m(t_3) \to V_m(t_3)\times \Phi_m(t_3).\]
Moreover, the image of the map $\mathcal{L}$ is precompact in the space
\[X:=C^0([0,t_3];\text{span}\{\boldsymbol{\eta_0},\boldsymbol{\eta_1},\dots,\boldsymbol{\eta_m}\})\times C^0([0,t_3];\text{span}\{\xi_0,\xi_1,\dots,\xi_m\}),\]
and $\mathcal{L}$ is continuous on $V_m(t_3)\times \Phi_m(t_3)$ in the topology of $X$.
\end{lem}
\begin{pf}
We separate the proof in three parts.
\newline
\textit{Part I}: $\mathcal{L}$ maps $V_m(t_3)\times \Phi_m(t_3)$ into itself for a small enough value $0<t_3\leq t_2$.
\newline
Combining the first and second equations of \eqref{eqn:odegalerkin} for $i=0,\dots, m$, we get
\begin{align*}
 \displaystyle &\gamma\int_{\Omega}\partial_t\mathbf{v}_m\cdot \boldsymbol{\eta}+\nu \int_{\Omega}\nabla \mathbf{v}_m\colon \nabla \boldsymbol{\eta} =\int_{\Omega}\mu_m \nabla \phi_m \cdot \boldsymbol{\eta}-\int_{\Omega}\frac{f'\left(\phi_m\right)}{2}\left(\mathbf{F}_m\colon \mathbf{F}_m-d\right)\nabla \phi_m\cdot \boldsymbol{\eta}\\
 \displaystyle &-\int_{\Omega}f\left(\phi_m\right)\left(\mathbf{F}_m\mathbf{F}_m^T\right)\colon \nabla \boldsymbol{\eta},\\
 \displaystyle &\int_{\Omega}\partial_t\phi_m \xi+\int_{\Omega}\left(\mathbf{v}_m \cdot \nabla \phi_m\right) \xi+\int_{\Omega}b\left(\phi_m\right)\nabla \mu_m \cdot \nabla \xi +\theta \int_{\Omega}\partial_t\mu_m \xi=0,
\end{align*}
with $\boldsymbol{\eta} \in \text{span}\left\{\boldsymbol{\eta}_0, \boldsymbol{\eta}_1, \dots, \boldsymbol{\eta}_m\right\}$, $\xi \in \text{span}\left\{\xi_0, \xi_1, \dots, \xi_m\right\}$.
Taking $\boldsymbol{\eta}=\mathbf{v}_m$, considering the regularity in space of $\boldsymbol{\eta}_i$ and $\xi_i$, the fact that $a_i^m, b_i^m \in C^0([0,t_2])$, Remark \ref{rem:gic}, Assumption \textbf{A3} and the inverse inequalities $\left|\left|\boldsymbol{\eta}_i\right|\right|_{L^{\infty}(\Omega;\mathbb{R}^d)}\leq C(m)\left|\left|\boldsymbol{\eta}_i\right|\right|_{L^{2}(\Omega;\mathbb{R}^d)}$, $\left|\left|\nabla \boldsymbol{\eta}_i\right|\right|_{L^{\infty}(\Omega;\mathbb{R}^{d\times d})}\leq C(m)\left|\left|\boldsymbol{\eta}_i\right|\right|_{L^{2}(\Omega;\mathbb{R}^d)}$, we get
\begin{align*}
&\frac{\gamma}{2}\frac{d}{dt}\left|\left|\mathbf{v}_m\right|\right|_{L^{2}(\Omega;\mathbb{R}^d)}^2+\nu \left|\left|\nabla \mathbf{v}_m\right|\right|_{L^{2}(\Omega;\mathbb{R}^{d\times d})}^2 \leq \left|\left|\mu_m \nabla \phi_m\right|\right|_{L^{2}(\Omega;\mathbb{R}^d)}\left|\left|\mathbf{v}_m\right|\right|_{L^{2}(\Omega;\mathbb{R}^d)}+\\
&C\left|\left|\nabla \phi_m\right|\right|_{L^{\infty}(\Omega;\mathbb{R}^d)}\left|\left|\mathbf{F}_m\right|\right|_{L^{\infty}\left(0,t_2;L^{2}(\Omega;\mathbb{R}^{d\times d}\right)}^2\left|\left|\mathbf{v}_m\right|\right|_{L^{\infty}(\Omega;\mathbb{R}^d)}+C\left|\left|\mathbf{F}_m\right|\right|_{L^{\infty}\left(0,t_2;L^{2}(\Omega;\mathbb{R}^{d\times d}\right)}^2\left|\left|\nabla \mathbf{v}_m\right|\right|_{L^{\infty}(\Omega;\mathbb{R}^{d\times d})} \\
&\leq C\left(m,L_v,L_{\phi},\mathbf{F}_0\right)\left|\left|\mathbf{v}_m\right|\right|_{L^{2}(\Omega;\mathbb{R}^d)}.
\end{align*}
Integrating in time and applying a Gronwall--Bellman type inequality, we obtain that
\[
{\rosso
\sup_{[0,t]}\left|\left|\mathbf{v}_m\right|\right|_{L^{2}(\Omega;\mathbb{R}^d)}\leq \left|\left|P_m^S\left(\mathbf{v}_0\right)\right|\right|_{L^{2}(\Omega;\mathbb{R}^d)}+C_1\left(m,L_v,L_{\phi},\mathbf{F}_0\right)t,}
\]
for $t\in [0,t_2]$.
Taking moreover $\xi=\phi_m$, considering the regularity in space of $\xi_i$, the fact that $ b_i^m \in W^{1,\infty}([0,t_2])$, Remark \ref{rem:gic}, Assumption \textbf{A1}, the fact that $\int_{\Omega}\left(\mathbf{v}_m \cdot \nabla \frac{|\phi_m|^2}{2}\right)=0$ and the inverse inequality $\left|\left|\nabla \xi_i\right|\right|_{L^{2}(\Omega;\mathbb{R}^d)}\leq C(m)\left|\left|\xi_i\right|\right|_{L^{2}(\Omega)}$, we get that 

\begin{align*}
&\frac{1}{2}\frac{d}{dt}\left|\left|\phi_m\right|\right|_{L^{2}(\Omega)}^2\leq b_1 \left|\left|\nabla \mu_m\right|\right|_{L^{2}(\Omega;\mathbb{R}^d)}\left|\left|\nabla \phi_m\right|\right|_{L^{2}(\Omega;\mathbb{R}^d)}+\theta \left|\left|\partial_t \mu_m\right|\right|_{L^{\infty}\left(0,t_2;L^{2}(\Omega)\right)}\left|\left|\phi_m\right|\right|_{L^{2}(\Omega)}\\
&\leq C\left(m,L_v,L_{\phi},\mathbf{F}_0\right)\left|\left|\phi_m\right|\right|_{L^{2}(\Omega)}.
\end{align*}
Integrating in time and applying a Gronwall--Bellman type inequality, we obtain that
\[
{\rosso
\sup_{t\in[0,t]}\left|\left|\phi_m\right|\right|_{L^{2}(\Omega)}\leq \left|\left|P_m^L\left(\phi_0\right)\right|\right|_{L^{2}(\Omega)}+C_2\left(m,L_v,L_{\phi},\mathbf{F}_0\right)t,}
\]
for $t\in [0,t_2]$.
Taking finally $t_3$ to be the maximum value such that $C_1t_3\leq1$ and $C_2t_3\leq 1$, from the definitions \eqref{eqn:vmpm} we obtain that $(\mathbf{v}_m,\phi_m)\in V_m(t_3)\times \Phi_m(t_3)$.\\
\newline
\textit{Part II}: the image of $\mathcal{L}$ is precompact in $X$.
\newline
Thanks to Remark \ref{rem:gic} and to the Ascoli--Arzel\'a theorem, the image of $\mathcal{L}$ is precompact in $X$, since the solution of \eqref{eqn:odegalerkin2} is uniformly Lipschitz continuous with respect to $(\bar{\mathbf{v}}_m,\bar{\phi}_m)$.
\\
\newline
\textit{Part III}: the map $\mathcal{L}$ is continuous on $V_m(t_3)\times \Phi_m(t_3)$.
\newline
In the rest of proof, we will omit to write the subscript $m$ for the discrete variables $\mathbf{v}_m, \mathbf{F}_m, \phi_m, \mu_m$ for ease of notation. Let us take a sequence $\{\bar{\mathbf{v}}_l=\sum_i \left(\bar{d}_i^m\right)_l\boldsymbol{\eta}_i\}_{l\in \mathbb{N}}$ in $V_m(t_3)$ converging to $\bar{\mathbf{v}}=\sum_i \bar{d}_i^m\boldsymbol{\eta}_i$ in $V_m(t_3)$, in the sense that $\left(\bar{d}_i^m\right)_l\to \bar{d}_i^m$ in $C^0([0,t_3])$ as $l\to \infty$. Let us moreover take a sequence $\{\bar{\phi}_l=\sum_i \left(\bar{a}_i^m\right)_l\xi_i\}_{l\in \mathbb{N}}$ in $\Phi_m(t_3)$ converging to $\bar{\phi}=\sum_i \bar{a}_i^m\xi_i$ in $\Phi_m(t_3)$, in the sense that $\left(\bar{a}_i^m\right)_l\to \bar{a}_i^m$ in $C^0([0,t_3])$ as $l\to \infty$. Then, if $\mathbf{F}_l$ is the solution of \eqref{eqn:fweak} with fixed values $\bar{\mathbf{v}}_l$ and $\bar{\phi}_l$, and if $\mathbf{F}$ is the solution of \eqref{eqn:fweak} with fixed values $\bar{\mathbf{v}}$ and $\bar{\phi}$, we have that
\begin{equation}
\label{eqn:fltof}
\mathbf{F}_l \to \mathbf{F} \quad \text{in} \quad L^{\infty}(0,t_3;L^2(\Omega;\mathbb{R}^{d\times d})).
\end{equation}
Indeed, taking the difference between equations \eqref{eqn:fweak} for $\mathbf{F}_l$ and for $\mathbf{F}$, and choosing $\boldsymbol{\Sigma}=\mathbf{F}_l-\mathbf{F}$, we get
\begin{align*}
& \left<\partial_t \left(\mathbf{F}_l-\mathbf{F}\right), \mathbf{F}_l-\mathbf{F} \right>+\underbrace{\int_{\Omega}\left(\bar{\mathbf{v}}_l \cdot \nabla \right)\left(\mathbf{F}_l-\mathbf{F}\right)\colon \left(\mathbf{F}_l-\mathbf{F}\right)}_{=0}+\int_{\Omega}\left(\left(\bar{\mathbf{v}}_l-\bar{\mathbf{v}}\right) \cdot \nabla \right)\mathbf{F}\colon \left(\mathbf{F}_l-\mathbf{F}\right)\\
& -\int_{\Omega}\left(\nabla \bar{\mathbf{v}}_l\right)\left(\mathbf{F}_l-\mathbf{F}\right)\colon \left(\mathbf{F}_l-\mathbf{F}\right)-\int_{\Omega}\nabla\left(\bar{\mathbf{v}}_l-\bar{\mathbf{v}}\right)\mathbf{F}\colon \left(\mathbf{F}_l-\mathbf{F}\right)+\lambda \int_{\Omega}f(\bar{\phi}_l)\nabla \left(\mathbf{F}_l-\mathbf{F}\right)\mathbin{\tensorm} \nabla \left(\mathbf{F}_l-\mathbf{F}\right)\\
& + \lambda \int_{\Omega}\left(f(\bar{\phi}_l)-f(\bar{\phi})\right)\nabla \mathbf{F}\mathbin{\tensorm} \nabla \left(\mathbf{F}_l-\mathbf{F}\right)+ \lambda \int_{\Omega}f'(\bar{\phi}_l) \left(\mathbf{F}_l-\mathbf{F}\right)\otimes \nabla \bar{\phi}_l\mathbin{\tensorm} \nabla \left(\mathbf{F}_l-\mathbf{F}\right)\\
& + \lambda \int_{\Omega}\mathbf{F} \otimes \nabla \left(f(\bar{\phi}_l)-f(\bar{\phi})\right)\mathbin{\tensorm} \nabla \left(\mathbf{F}_l-\mathbf{F}\right)=0,
\end{align*}
where the second term is equal to zero {\rosso due to the boundary condition on $\bar{\mathbf{v}}_l$ and to the fact that $\diver\mathbf{v}_l=0$}. Integrating in time over the interval $(0,t_3)$, using the Cauchy--Schwarz and the Young inequalities, Assumption \textbf{A3}, \eqref{eqn:fspace} and considering that $\mathbf{F}_l(\mathbf{x},0)=\mathbf{F}(\mathbf{x},0)$, we obtain the inequality
\begin{align*}
&\frac{1}{2}\left|\left|\mathbf{F}_l-\mathbf{F}\right|\right|_{L^2(\Omega;\mathbb{R}^{d\times d})}^2+\lambda f_{\text{min}}\int_0^{t_3}\int_{\Omega}\left|\nabla \left(\mathbf{F}_l-\mathbf{F}\right)\right|^2 \\
& \leq \left|\left|\bar{\mathbf{v}}_l-\bar{\mathbf{v}}\right|\right|_{L^{\infty}\left(0,t_3;L^{\infty}(\Omega;\mathbb{R}^{d})\right)}\int_0^{t_3}\left|\left|\nabla \mathbf{F}\right|\right|_{L^2(\Omega;\mathbb{R}^{d\times d\times d})}\left|\left|\mathbf{F}_l-\mathbf{F}\right|\right|_{L^2(\Omega;\mathbb{R}^{d\times d})}\\
&+\left|\left|\nabla \bar{\mathbf{v}}_l\right|\right|_{L^{\infty}\left(0,t_3;L^{\infty}(\Omega;\mathbb{R}^{d\times d})\right)}\int_0^{t_3}\left|\left|\mathbf{F}_l-\mathbf{F}\right|\right|_{L^2(\Omega;\mathbb{R}^{d\times d})}^2\\
&+\left|\left|\nabla \left(\bar{\mathbf{v}}_l-\bar{\mathbf{v}}\right)\right|\right|_{L^{\infty}\left(0,t_3;L^{\infty}(\Omega;\mathbb{R}^{d\times d})\right)}\int_0^{t_3}\left|\left|\mathbf{F}\right|\right|_{L^2(\Omega;\mathbb{R}^{d\times d})}\left|\left|\mathbf{F}_l-\mathbf{F}\right|\right|_{L^2(\Omega;\mathbb{R}^{d\times d})}\\
& +\lambda \left|\left|f(\bar{\phi}_l)-f(\bar{\phi})\right|\right|_{L^{\infty}\left(0,t_3;L^{\infty}(\Omega)\right)}\int_0^{t_3}\left|\left|\nabla \mathbf{F}\right|\right|_{L^2(\Omega;\mathbb{R}^{d\times d\times d})}\left|\left|\nabla \left(\mathbf{F}_l-\mathbf{F}\right)\right|\right|_{L^2(\Omega;\mathbb{R}^{d\times d\times d})}\\
&+\lambda f_1 \left|\left|\nabla \bar{\phi}_l\right|\right|_{L^{\infty}\left(0,t_3;L^{\infty}(\Omega;\mathbb{R}^d)\right)}\int_0^{t_3}\left|\left|\mathbf{F}_l-\mathbf{F}\right|\right|_{L^2(\Omega;\mathbb{R}^{d\times d})}\left|\left|\nabla \left(\mathbf{F}_l-\mathbf{F}\right)\right|\right|_{L^2(\Omega;\mathbb{R}^{d\times d\times d})}\\
& +\lambda \left|\left|\nabla \left(f(\bar{\phi}_l)-f(\bar{\phi})\right)\right|\right|_{L^{\infty}\left(0,t_3;L^{\infty}(\Omega;\mathbb{R}^d)\right)}\int_0^{t_3}\left|\left|\mathbf{F}\right|\right|_{L^2(\Omega;\mathbb{R}^{d\times d})}\left|\left|\nabla \left(\mathbf{F}_l-\mathbf{F}\right)\right|\right|_{L^2(\Omega;\mathbb{R}^{d\times d\times d})} \\
& \leq C\int_0^{t_3}\left|\left|\mathbf{F}_l-\mathbf{F}\right|\right|_{L^2(\Omega;\mathbb{R}^{d\times d})}^2+C||\bar{\mathbf{d}}_l-\bar{\mathbf{d}}||_{L^{\infty}\left(0,t_3\right)}^2\left(\int_0^{t_3}\left|\left|\nabla \mathbf{F}\right|\right|_{L^2(\Omega;\mathbb{R}^{d\times d\times d})}^2+\int_0^{t_3}\left|\left|\mathbf{F}\right|\right|_{L^2(\Omega;\mathbb{R}^{d\times d})}^2\right)\\
& +\frac{\lambda}{2} f_{\text{min}}\int_0^{t_3}\int_{\Omega}\left|\nabla \left(\mathbf{F}_l-\mathbf{F}\right)\right|^2+C||\bar{\mathbf{a}}_l-\bar{\mathbf{a}}||_{L^{\infty}\left(0,t_3\right)}^2\left(\int_0^{t_3}\left|\left|\nabla \mathbf{F}\right|\right|_{L^2(\Omega;\mathbb{R}^{d\times d\times d})}^2+\int_0^{t_3}\left|\left|\mathbf{F}\right|\right|_{L^2(\Omega;\mathbb{R}^{d\times d})}^2\right)\\
& \leq\frac{\lambda}{2} f_{\text{min}}\int_0^{t_3}\int_{\Omega}\left|\nabla \left(\mathbf{F}_l-\mathbf{F}\right)\right|^2+C||\bar{\mathbf{d}}_l-\bar{\mathbf{d}}||_{L^{\infty}\left(0,t_3\right)}^2+C||\bar{\mathbf{a}}_l-\bar{\mathbf{a}}||_{L^{\infty}\left(0,t_3\right)}^2+C\int_0^{t_3}\left|\left|\mathbf{F}_l-\mathbf{F}\right|\right|_{L^2(\Omega;\mathbb{R}^{d\times d})}^2,
\end{align*}
from which, using a Gronwall argument and from the assumed convergence $\bar{\mathbf{d}}_l\to \bar{\mathbf{d}}$ and $\bar{\mathbf{a}}_l\to \bar{\mathbf{a}}$ in $\left[C^0([0,t_3])\right]^{m+1}$, we get \eqref{eqn:fltof}.
\newline
We now show that, if \eqref{eqn:fltof} is valid, then correspondingly
\begin{align}
\label{eqn:abdltoabd}
&\left({d}_i^m\right)_l\to {d}_i^m \quad \text{in} \quad C^0([0,t_3]),\\
& \notag \left({a}_i^m\right)_l\to {a}_i^m \quad \text{in} \quad C^0([0,t_3]),\\
& \notag \left({b}_i^m\right)_l\to {b}_i^m \quad \text{in} \quad C^0([0,t_3]),
\end{align}
for $i=0,\dots, m$, which give the continuity of $\mathcal{L}$ in the proposed topology.
Taking the difference between the equations for $\mathbf{d}_l$ and $\mathbf{d}$ in \eqref{eqn:odegalerkin2} and taking its scalar product with $\mathbf{d}_l-\mathbf{d}$, we get
\begin{align}
\label{eqn:dld}
&\frac{\gamma}{2}\frac{d}{dt}|\mathbf{d}_l-\mathbf{d}|^2+\nu \text{diag}(\boldsymbol{\beta})\left(\mathbf{d}_l-\mathbf{d}\right)\cdot \left(\mathbf{d}_l-\mathbf{d}\right)=\mathbf{D}\left((\mathbf{b}_l-\mathbf{b})\otimes \mathbf{a}\right)\cdot \left(\mathbf{d}_l-\mathbf{d}\right)\\
&\notag +\mathbf{D}\left(\mathbf{b}_l\otimes \left(\mathbf{a}_l-\mathbf{a}\right)\right)\cdot \left(\mathbf{d}_l-\mathbf{d}\right)-\mathbf{E}\left(\mathbf{a}_l,\mathbf{F}_l-\mathbf{F}\right)\mathbf{a}_l\cdot \left(\mathbf{d}_l-\mathbf{d}\right)-\mathbf{E}\left(\mathbf{a}_l,\mathbf{F}\right)\left(\mathbf{a}_l-\mathbf{a}\right)\cdot \left(\mathbf{d}_l-\mathbf{d}\right)\\
&\notag -\mathbf{E}\left(\mathbf{a}_l-\mathbf{a},\mathbf{F}\right)\mathbf{a}\cdot \left(\mathbf{d}_l-\mathbf{d}\right)-\mathbf{p}\left(\mathbf{a}_l,\mathbf{F}_l-\mathbf{F}\right)\cdot \left(\mathbf{d}_l-\mathbf{d}\right)-\mathbf{p}\left(\mathbf{a}_l-\mathbf{a},\mathbf{F}\right)\cdot \left(\mathbf{d}_l-\mathbf{d}\right).
\end{align} 
Taking the difference between the second equations for $\mathbf{a}_l, \mathbf{b}_l$ and $\mathbf{a}, \mathbf{b}$ in \eqref{eqn:odegalerkin2} and taking its scalar product with $\left(\mathbf{a}_l-\mathbf{a}\right)+\frac{d}{dt}\left(\mathbf{a}_l-\mathbf{a}\right)$, we get
{\rosso
\begin{align}
\label{eqn:ala}
&\text{diag}(\boldsymbol{\alpha})\left(\mathbf{a}_l-\mathbf{a}\right)\cdot \left(\mathbf{a}_l-\mathbf{a}\right)+\left(\mathbf{r}(\mathbf{a}_l)-\mathbf{r}(\mathbf{a})\right)\cdot \left(\mathbf{a}_l-\mathbf{a}\right)+\mathbf{q}\left(\mathbf{a}_l-\mathbf{a},\mathbf{F}_l\right)\cdot \left(\mathbf{a}_l-\mathbf{a}\right)+\mathbf{q}\left(\mathbf{a},\mathbf{F}_l-\mathbf{F}\right)\cdot \left(\mathbf{a}_l-\mathbf{a}\right)\\
& \notag +\frac{\delta}{2}\frac{d}{dt}\left|\mathbf{a}_l-\mathbf{a}\right|^2+\text{diag}(\boldsymbol{\alpha})\left(\mathbf{a}_l-\mathbf{a}\right)\cdot \frac{d}{dt}\left(\mathbf{a}_l-\mathbf{a}\right)+\left(\mathbf{r}(\mathbf{a}_l)-\mathbf{r}(\mathbf{a})\right)\cdot \frac{d}{dt}\left(\mathbf{a}_l-\mathbf{a}\right)+\mathbf{q}\left(\mathbf{a}_l-\mathbf{a},\mathbf{F}\right)\cdot \frac{d}{dt}\left(\mathbf{a}_l-\mathbf{a}\right)\\
& \notag +\mathbf{q}\left(\mathbf{a},\mathbf{F}_l-\mathbf{F}\right)\cdot \frac{d}{dt}\left(\mathbf{a}_l-\mathbf{a}\right)+\delta\left|\frac{d}{dt}\left(\mathbf{a}_l-\mathbf{a}\right)\right|^2= \left(\mathbf{b}_l-\mathbf{b}\right)\cdot \left(\mathbf{a}_l-\mathbf{a}\right)+\left(\mathbf{b}_l-\mathbf{b}\right)\cdot \frac{d}{dt}\left(\mathbf{a}_l-\mathbf{a}\right).
\end{align}
}
Taking finally the difference between the third equations for $\mathbf{a}_l, \mathbf{b}_l$ and $\mathbf{a}, \mathbf{b}$ in \eqref{eqn:odegalerkin2} and taking its scalar product with $\left(\mathbf{b}_l-\mathbf{b}\right)$, we get
\begin{align}
\label{eqn:blb}
&\frac{d}{dt}\left(\mathbf{a}_l-\mathbf{a}\right)\cdot \left(\mathbf{b}_l-\mathbf{b}\right)=-\mathbf{B}\left((\mathbf{d}_l-\mathbf{d})\otimes \mathbf{a}_l\right)\cdot \left(\mathbf{b}_l-\mathbf{b}\right)-\mathbf{B}\left(\mathbf{d}\otimes (\mathbf{a}_l-\mathbf{a})\right)\cdot \left(\mathbf{b}_l-\mathbf{b}\right)\\
& \notag -\left(\mathbf{C}(\mathbf{a}_l)-\mathbf{C}(\mathbf{a})\right)\mathbf{b}_l\cdot \left(\mathbf{b}_l-\mathbf{b}\right)-\mathbf{C}(\mathbf{a})\left(\mathbf{b}_l-\mathbf{b}\right)\cdot \left(\mathbf{b}_l-\mathbf{b}\right)-\frac{\theta}{2}\frac{d}{dt}\left|\mathbf{b}_l-\mathbf{b}\right|^2.
\end{align}
We use the identity
\begin{equation}
\label{eqn:ffid1}
\mathbf{F}_l \colon \mathbf{F}_l - \mathbf{F} \colon \mathbf{F}=\left(\mathbf{F}_l-\mathbf{F}\right)\colon \left(\mathbf{F}_l+\mathbf{F}\right) 
\end{equation}
in the expressions for $\mathbf{E}\left(\mathbf{a}_l,\mathbf{F}_l-\mathbf{F}\right)$ and $\mathbf{q}\left(\mathbf{a},\mathbf{F}_l-\mathbf{F}\right)$ from \eqref{eqn:tensors} to get, using \eqref{eqn:fspace}, Assumption \textbf{A3}, the Cauchy--Schwarz and the Young inequalities, that
\[
\left|\mathbf{E}\left(\mathbf{a}_l,\mathbf{F}_l-\mathbf{F}\right)\mathbf{a}_l\cdot \left(\mathbf{d}_l-\mathbf{d}\right)\right|\leq C\left|\left|\mathbf{F}_l-\mathbf{F}\right|\right|_{L^{\infty}(0,t_3;L^2(\Omega;\mathbb{R}^{d\times d}))}^2+C|\mathbf{d}_l-\mathbf{d}|^2,
\]
\[
\mathbf{q}\left(\mathbf{a},\mathbf{F}_l-\mathbf{F}\right)\cdot \left(\left(\mathbf{a}_l-\mathbf{a}\right)+\frac{d}{dt}\left(\mathbf{a}_l-\mathbf{a}\right)\right)
\leq C\left|\left|\mathbf{F}_l-\mathbf{F}\right|\right|_{L^{\infty}(0,t_3;L^2(\Omega;\mathbb{R}^{d\times d}))}^2+C|\mathbf{a}_l-\mathbf{a}|^2+\frac{\delta}{10}\left|\frac{d}{dt}\left(\mathbf{a}_l-\mathbf{a}\right)\right|^2.
\]
We moreover use the identity
\begin{equation}
\label{eqn:ffid2}
\mathbf{F}_l \mathbf{F}_l^T - \mathbf{F} \mathbf{F}^T=\left(\mathbf{F}_l-\mathbf{F}\right)\left(\mathbf{F}_l+\mathbf{F}\right)^T+ \left(\mathbf{F}-\mathbf{F}_l\right)\mathbf{F}^T+\mathbf{F}\left(\mathbf{F}_l-\mathbf{F}\right)^T
\end{equation}
in the expression for $\mathbf{p}\left(\mathbf{a}_l,\mathbf{F}_l-\mathbf{F}\right)$ from \eqref{eqn:tensors} to get analogously
\[
\left|\mathbf{p}\left(\mathbf{a}_l,\mathbf{F}_l-\mathbf{F}\right)\cdot \left(\mathbf{d}_l-\mathbf{d}\right)\right|\leq C\left|\left|\mathbf{F}_l-\mathbf{F}\right|\right|_{L^{\infty}(0,t_3;L^2(\Omega;\mathbb{R}^{d\times d}))}^2+C|\mathbf{d}_l-\mathbf{d}|^2.
\]
Summing \eqref{eqn:dld}, \eqref{eqn:ala} and \eqref{eqn:blb}, using the Assumptions \textbf{A1, A2, A3}, the Cauchy--Schwarz and the Young inequalities and the previous arguments, we obtain
\begin{align*}
&\frac{\gamma}{2}\frac{d}{dt}|\mathbf{d}_l-\mathbf{d}|^2+\nu \text{diag}(\boldsymbol{\beta})\left(\mathbf{d}_l-\mathbf{d}\right)\cdot \left(\mathbf{d}_l-\mathbf{d}\right)+\text{diag}(\boldsymbol{\alpha})\left(\mathbf{a}_l-\mathbf{a}\right)\cdot \left(\mathbf{a}_l-\mathbf{a}\right)+\frac{\delta}{2}\frac{d}{dt}\left|\mathbf{a}_l-\mathbf{a}\right|^2 \\
& +\delta\left|\frac{d}{dt}\left(\mathbf{a}_l-\mathbf{a}\right)\right|^2+\mathbf{C}(\mathbf{a})\left(\mathbf{b}_l-\mathbf{b}\right)\cdot \left(\mathbf{b}_l-\mathbf{b}\right)+\frac{\theta}{2}\frac{d}{dt}\left|\mathbf{b}_l-\mathbf{b}\right|^2\\
& \leq C|\mathbf{d}_l-\mathbf{d}|^2+C|\mathbf{a}_l-\mathbf{a}|^2+C|\mathbf{b}_l-\mathbf{b}|^2+C\left|\left|\mathbf{F}_l-\mathbf{F}\right|\right|_{L^{\infty}(0,t_3;L^2(\Omega;\mathbb{R}^{d\times d}))}^2+5\frac{\delta}{10}\left|\frac{d}{dt}\left(\mathbf{a}_l-\mathbf{a}\right)\right|^2,
\end{align*}
which gives that
\begin{align}
\label{eqn:continuity}
&\frac{d}{dt}\left(\frac{\gamma}{2}|\mathbf{d}_l-\mathbf{d}|^2+\frac{\delta}{2}\left|\mathbf{a}_l-\mathbf{a}\right|^2+\frac{\theta}{2}\left|\mathbf{b}_l-\mathbf{b}\right|^2\right)\\
& \notag \leq C\left(\frac{\gamma}{2}|\mathbf{d}_l-\mathbf{d}|^2+\frac{\delta}{2}\left|\mathbf{a}_l-\mathbf{a}\right|^2+\frac{\theta}{2}\left|\mathbf{b}_l-\mathbf{b}\right|^2\right)+ C\left|\left|\mathbf{F}_l-\mathbf{F}\right|\right|_{L^{\infty}(0,t_3;L^2(\Omega;\mathbb{R}^{d\times d}))}^2.
\end{align}
Integrating \eqref{eqn:continuity} in time {\rosso with respect to $t$, for $t\in [0,t_3]$}, and applying a Gronwall argument, given the convergence result \eqref{eqn:fltof}, we establish the convergence results \eqref{eqn:abdltoabd}, hence the continuity of the map $\mathcal{L}$.
\end{pf}
\begin{rem}
\label{rem:theta}
In the case $\delta=\theta=0$ (which will be considered in the 2-d case), the proof of the continuity of the map $\mathcal{L}$ could be obtained by a simpler argument by testing the difference between the equations for $\mathbf{d}_l$ and $\mathbf{d}$ in \eqref{eqn:odegalerkin2} by $\mathbf{d}_l-\mathbf{d}$, the difference between the second equations for $\mathbf{b}_l$ and $\mathbf{b}$ by $\mathbf{b}_l-\mathbf{b}$, and the difference between the third equations for $\mathbf{a}_l$ and $\mathbf{a}$ by $\mathbf{a}_l-\mathbf{a}$. Since, in the case $\delta>0, \theta=0$, the continuity of the map $\mathcal{L}$ would be true only for sufficiently small values of $\delta$, the term $\theta \partial_t \mu_m$ has been introduced in \eqref{eqn:galerkin} to obtain a more general result.
\end{rem}
Thanks to Lemma \ref{lem:fixedpoint}, we can apply the Schauder fixed point and prove that there exists a fixed point of $\mathcal{L}$, which is a solution of system \eqref{eqn:galerkin} in the time interval $[0,t_3]$.
\subsection{Step $4$: A priori estimates and convergence to a continuous weak solution} 
We now deduce a priori estimates, uniform in the discretization parameter $m$, for the solutions of system \eqref{eqn:galerkin}, which can be rewritten, summing the first, the third and the fourth equations over $i=0, \dots, m$, as
\begin{equation}
\label{eqn:galerkin2}
\begin{cases}
 \displaystyle \gamma\int_{\Omega}\partial_t\mathbf{v}_m\cdot \boldsymbol{\eta}+\nu \int_{\Omega}\nabla \mathbf{v}_m\colon \nabla \boldsymbol{\eta} =\int_{\Omega}\mu_m \nabla \phi_m \cdot \boldsymbol{\eta}-\int_{\Omega}\frac{f'\left(\phi_m\right)}{2}\left(\mathbf{F}_m\colon \mathbf{F}_m-d\right)\nabla \phi_m\cdot \boldsymbol{\eta}\\
 \displaystyle -\int_{\Omega}f\left(\phi_m\right)\mathbf{F}_m\mathbf{F}_m^T\colon \nabla \boldsymbol{\eta},\\
\displaystyle  \left<\partial_t \mathbf{F}_m, \boldsymbol{\Sigma} \right>+\int_{\Omega}\left({\mathbf{v}}_m \cdot \nabla \right)\mathbf{F}_m\colon \boldsymbol{\Sigma}-\int_{\Omega}\left(\nabla {\mathbf{v}}_m\right)\mathbf{F}_m\colon \boldsymbol{\Sigma}+\lambda \int_{\Omega}\nabla \left(f({\phi}_m)\mathbf{F}_m\right)\mathbin{\tensorm} \nabla \boldsymbol{\Sigma}=0,\\
 \displaystyle \int_{\Omega}\partial_t\phi_m \xi+\int_{\Omega}\left(\mathbf{v}_m \cdot \nabla \phi_m\right) \xi+\int_{\Omega}b\left(\phi_m\right)\nabla \mu_m \cdot \nabla \xi +\theta \int_{\Omega}\partial_t\mu_m \xi=0,\\
\displaystyle \int_{\Omega}\mu_m \chi=\int_{\Omega}\nabla \phi_m\cdot \nabla \chi+\int_{\Omega}\psi'\left(\phi_m\right)\chi+\int_{\Omega}\frac{f'\left(\phi_m\right)}{2}\left(\mathbf{F}_m\colon \mathbf{F}_m-d\right)\chi+\delta \int_{\Omega}\partial_t\phi_m \chi,
\end{cases}
\end{equation} 
for a.e. $t \in [0,t_3]$, with $\boldsymbol{\eta} \in \text{span}\left\{\boldsymbol{\eta}_0, \boldsymbol{\eta}_1, \dots, \boldsymbol{\eta}_m\right\}$, $\xi, \chi \in \text{span}\left\{\xi_0, \xi_1, \dots, \xi_m\right\}$ and $\forall \boldsymbol{\Sigma} \in H^1(\Omega;\mathbb{R}^{d\times d})$, with initial conditions defined in \eqref{eqn:galerkin}. We take $\boldsymbol{\eta}=\mathbf{v}_m$, $\Sigma=f(\phi_m)\mathbf{F}_m$, $\xi=\mu_m$, {\rosso $\chi=-\frac{\partial \phi_m}{\partial t}$} in \eqref{eqn:galerkin2}, sum all the equations and integrate in time between $0$ and $t\in [0,t_3]$. Note that, thanks to Assumption \textbf{A3} and \eqref{eqn:fspace}, $f(\phi_m)\mathbf{F}_m$ is a.e. in $H^1(\Omega;\mathbb{R}^{d\times d})$. We get, for any $t\in [0,t_3]$,
\begin{align}
\label{eqn:apriori1}
&\frac{\gamma}{2}||\mathbf{v}_m||_{L^2(\Omega;\mathbb{R}^d)}^2+\nu\int_0^{t}\int_{\Omega}\left|\nabla \mathbf{v}_m\right|^2+\lambda \int_0^{t}\int_{\Omega}\left|\nabla \left(f(\phi_m)\mathbf{F}_m\right)\right|^2+\int_0^{t}\int_{\Omega}b(\phi_m)\left|\nabla \mu_m\right|^2\\
& \notag +\int_{\Omega}\left(\frac{1}{2}f(\phi_m)\left(\mathbf{F}_m\colon \mathbf{F}_m-d\right)+\frac{1}{2}|\nabla \phi_m|^2+\psi(\phi_m)\right)+\frac{\theta}{2}||\mu_m||_{L^2(\Omega)}^2+\delta \int_0^{t}\int_{\Omega}\left|\frac{\partial \phi_m}{\partial t}\right|^2 \\
& \notag = -\int_0^{t}\int_{\Omega}\frac{f'\left(\phi_m\right)}{2}\left(\mathbf{F}_m\colon \mathbf{F}_m-d\right)\nabla \phi_m\cdot \mathbf{v}_m-\int_0^{t}\int_{\Omega}\left({\mathbf{v}}_m \cdot \nabla \right)\mathbf{F}_m\colon f(\phi_m)\mathbf{F}_m\\
& \notag +\frac{\gamma}{2}||\mathbf{v}_m(0)||_{L^2(\Omega;\mathbb{R}^d)}^2+\frac{\theta}{2}||\mu_m(0)||_{L^2(\Omega)}^2\\
& \notag +\int_{\Omega}\left(\frac{1}{2}f(\phi_m(0))\left(\mathbf{F}_m(0)\colon \mathbf{F}_m(0)-d\right)+\frac{1}{2}|\nabla \phi_m(0)|^2+\psi(\phi_m(0))\right).
\end{align}
Similarly to \eqref{eqn:gradgradparts}, we rewrite the first term on the right hand side of \eqref{eqn:apriori1}, integrating by parts, as
\begin{align*}
&-\int_0^{t}\int_{\Omega}\frac{f'\left(\phi_m\right)}{2}\left(\mathbf{F}_m\colon \mathbf{F}_m-d\right)\nabla \phi_m\cdot \mathbf{v}_m=-\frac{1}{2}\underbrace{\int_0^{t}\int_{\Omega}\nabla \left(f(\phi_m)\left(\mathbf{F}_m\colon \mathbf{F}_m-d\right)\right)\cdot \mathbf{v}_m}_{=0}\\
& +\int_0^{t}\int_{\Omega}\left({\mathbf{v}}_m \cdot \nabla \right)\mathbf{F}_m\colon f(\phi_m)\mathbf{F}_m,
\end{align*}
which is valid, for any $t\in [0,t_3]$, thanks to \eqref{eqn:fspace} and the regularity in space and time of $\phi_m$ and $\mathbf{v}_m$. 
Hence, \eqref{eqn:apriori1} becomes
\begin{align}
\label{eqn:apriori}
& \notag \sup_{t\in [0,t_3]}\left(\frac{\gamma}{2}||\mathbf{v}_m||_{L^2(\Omega;\mathbb{R}^d)}^2+\int_{\Omega}\left(\frac{1}{2}f(\phi_m)\left(\mathbf{F}_m\colon \mathbf{F}_m-d\right)+\frac{1}{2}|\nabla \phi_m|^2+\psi(\phi_m)\right)+\frac{\theta}{2}||\mu_m||_{L^2(\Omega)}^2\right)\\
& \notag +\nu\int_0^{t_3}\int_{\Omega}\left|\nabla \mathbf{v}_m\right|^2+\lambda \int_0^{t_3}\int_{\Omega}\left|\nabla \left(f(\phi_m)\mathbf{F}_m\right)\right|^2+\int_0^{t_3}\int_{\Omega}b(\phi_m)\left|\nabla \mu_m\right|^2+\delta \int_0^{t_3}\int_{\Omega}\left|\frac{\partial \phi_m}{\partial t}\right|^2 \\
& \leq C(\mathbf{v}_0,\mathbf{F}_0,\phi_0,\mu_0),
\end{align}
where we used Assumptions \textbf{A2}, \textbf{A3} and {\rosso \textbf{A4 Bis}}. The constant in the right hand side of \eqref{eqn:apriori} depends only on the initial data and not on the discretization parameter $m$ and on the regularization parameters $\gamma, \theta, \delta$. Thanks to the a priori estimate \eqref{eqn:apriori}, we may extend {\rosso by continuity} the local solution of system \eqref{eqn:galerkin} to the interval $[0,T]$.
{\rosso
Next, taking $\xi=1$ (which is a multiple of $\xi_0$) in \eqref{eqn:galerkin2}$_3$, we obtain, using the facts that $\diver\mathbf{v}_m=0$, $\mathbf{v}_m=\mathbf{0}$ on $\partial \Omega$ and integrating by parts in the second term, that
\begin{equation}
\label{eqn:massphimdisc}
\left(\partial_t\phi_m,1\right)=-\theta\left(\partial_t\mu_m,1\right).
\end{equation}
Integrating in time both sides of \eqref{eqn:massphimdisc} and using Assumption \textbf{A4 Bis} and \eqref{eqn:apriori}, we get that
\begin{equation}
\label{eqn:massphimdisc2}
|\left(\phi_m(t),1\right)|= |\left(\phi_m(0),1\right)+\theta\left(\mu_m(0),1\right)-\theta\left(\mu_m(t),1\right)|\leq C(\mathbf{v}_0,\mathbf{F}_0,\phi_0,\mu_0),
\end{equation}
for all $t\in[0,T]$. Hence, from \eqref{eqn:apriori}, \eqref{eqn:massphimdisc2} and the Poincar\'e--Wirtinger inequality we deduce that
\begin{equation}
\label{eqn:phimh1}
||\phi_m||_{L^{\infty}(0,T;H^1(\Omega))}\leq C(\mathbf{v}_0,\mathbf{F}_0,\phi_0,\mu_0).
\end{equation}
Now, taking $\chi= 1$ in \eqref{eqn:galerkin2}$_4$, we obtain, using Assumptions \textbf{A2, A3}, \eqref{eqn:phimh1} and the uniform boundedness of $\mathbf{F}_m$ in $L^{\infty}(0,T;L^2(\Omega,\mathbb{R}^{d\times d}))$ (cf. \eqref{eqn:apriori}), that
\begin{align}
\label{eqn:mumdiscrete}
&\left|(\mu_m,1)\right|=\left|\int_{\Omega}\psi'\left(\phi_m\right)+\int_{\Omega}\frac{f'\left(\phi_m\right)}{2}\left(\mathbf{F}_m\colon \mathbf{F}_m-d\right)+\delta \int_{\Omega}\partial_t \phi_m\right|\\
& \notag \leq C\left(||\phi_m||_{H^1(\Omega)}^p+1\right)+C||\mathbf{F}_m||_{L^2(\Omega,\mathbb{R}^{d\times d})}^2+\delta||\partial_t \phi_m||_{L^2(\Omega)}+C,
\end{align}
where $p$ is as in Assumption \textbf{A2}. Hence, $|(\mu_m,1)|$ is bounded in $L^2(0,T)$ and consequently, by the Poincar\'e--Wirtinger inequality and \eqref{eqn:apriori},
\begin{equation}
\label{eqn:mumh1}
||\mu_m||_{L^{2}(0,T;H^1(\Omega))}\leq C(\mathbf{v}_0,\mathbf{F}_0,\phi_0,\mu_0).
\end{equation}}%

Before studying the convergence properties of the discrete solution of system \eqref{eqn:galerkin2} to a continuous solution, identifying its limit as the discretization parameter $m$ tends to infinity, we identify preliminarly its limit at the discrete level as the regularization parameters $\gamma, \theta \to 0$. Indeed, as described at the beginning of Section \ref{sec:2}, we want to identify the limit problem of \eqref{eqn:galerkin2} in the case $\gamma=\theta=0$, {\rosso with $m\in \mathbb{N}$ and $\delta>0$ being fixed in this limit procedure. We thus multiply the equations in \eqref{eqn:galerkin2} by a function $\omega \in C_0^{\infty}([0,T])$} and integrate in time over the interval $[0,T]$, letting $\gamma,\theta \to 0$. From the a priori estimate \eqref{eqn:apriori}, from \eqref{eqn:mumdiscrete} and \eqref{eqn:fnorm}, which are uniform in the parameters $\gamma,\theta$, we have that, at least for a subsequence which we still label by the index $m$ without reporting the dependence on $\gamma,\theta$,
\begin{align}
\label{eqn:convreg}
&\sqrt{\gamma}\mathbf{v}_m\overset{\ast}{\rightharpoonup} \hat{\mathbf{v}}_m, \quad \text{in} \quad L^{\infty}\left(0,T;\text{span}\left\{\boldsymbol{\eta}_0, \boldsymbol{\eta}_1, \dots, \boldsymbol{\eta}_m\right\}\right),\\
& \notag \mathbf{v}_m\rightharpoonup \tilde{\mathbf{v}}_m, \quad \text{in} \quad L^{2}\left(0,T;\text{span}\left\{\boldsymbol{\eta}_0, \boldsymbol{\eta}_1, \dots, \boldsymbol{\eta}_m\right\}\right),\\
& \notag \sqrt{\theta}\mu_m\overset{\ast}{\rightharpoonup} \hat{\mu}_m, \quad \text{in} \quad L^{\infty}\left(0,T;\text{span}\left\{\xi_0, \xi_1, \dots, \xi_m\right\}\right),
\\
& \notag \mu_m\rightharpoonup \tilde{\mu}_m, \quad \text{in} \quad L^{2}\left(0,T;\text{span}\left\{\xi_0, \xi_1, \dots, \xi_m\right\}\right),
\\
& \notag
\phi_m\overset{\ast}{\rightharpoonup} \tilde{\phi}_m, \quad \text{in} \quad L^{\infty}\left(0,T;\text{span}\left\{\xi_0, \xi_1, \dots, \xi_m\right\}\right),
\\
& \notag
\frac{\partial \phi_m}{\partial t}\rightharpoonup \frac{\partial \tilde{\phi}_m}{\partial t}, \quad \text{in} \quad L^{2}\left(0,T;\text{span}\left\{\xi_0, \xi_1, \dots, \xi_m\right\}\right),
\\
& \notag
\mathbf{F}_m \overset{\ast}{\rightharpoonup} \tilde{\mathbf{F}}_m \quad \text{in} \quad L^{\infty}\left(0,T;L^2(\Omega;\mathbb{R}^{d\times d})\right) \cap L^{2}\left(0,T;H^1(\Omega;\mathbb{R}^{d\times d})\right),
\\
& \notag
\frac{\partial \mathbf{F}_m}{\partial t} {\rightharpoonup} \frac{\partial \tilde{\mathbf{F}}_m}{\partial t} \quad \text{in} \quad L^{2}\left(0,T;\left(H^1(\Omega;\mathbb{R}^{d\times d})\right)'\right).
\end{align}
Hence, 
\begin{equation}
\label{eqn:limgamma}
\sqrt{\gamma}\int_0^T\omega(t)\int_{\Omega}\sqrt{\gamma}\partial_t\mathbf{v}_m\cdot \boldsymbol{\eta}=-\sqrt{\gamma}\int_0^T\frac{d}{dt}\omega(t)\int_{\Omega}\sqrt{\gamma}\mathbf{v}_m\cdot \boldsymbol{\eta}\to 0,
\end{equation}
\begin{equation}
\label{eqn:limtheta}
\sqrt{\theta}\int_0^T\omega(t)\int_{\Omega}\sqrt{\theta}\partial_t\mu_m\xi=-\sqrt{\theta}\int_0^T\frac{d}{dt}\omega(t)\int_{\Omega}\sqrt{\theta}\mu_m\xi \to 0,
\end{equation}
as $\gamma,\theta \to 0$. {\rosso Moreover, note that, as a consequence of the second and fourth convergence properties in \eqref{eqn:convreg}, it turns out that $\hat{\mathbf{v}}_m\equiv \mathbf{0}$ and $\hat{\mu}_m\equiv 0$}. It's easy to identify the limit of the remaining terms in \eqref{eqn:galerkin2}, as $\gamma,\theta \to 0$, thanks to \eqref{eqn:convreg} and exploiting the regularity in space of $\phi_m, \mu_m, \mathbf{v}_m$. We finally obtain the following system, where we do not report the tilde notation for the limit point for ease of notation, 
\begin{equation}
\label{eqn:galerkin3}
\begin{cases}
 \displaystyle \nu \int_{\Omega}\nabla \mathbf{v}_m\colon \nabla \boldsymbol{\eta} =\int_{\Omega}\mu_m \nabla \phi_m \cdot \boldsymbol{\eta}-\int_{\Omega}\frac{f'\left(\phi_m\right)}{2}\left(\mathbf{F}_m\colon \mathbf{F}_m-2\right)\nabla \phi_m\cdot \boldsymbol{\eta}-\int_{\Omega}f\left(\phi_m\right)\mathbf{F}_m\mathbf{F}_m^T\colon \nabla \boldsymbol{\eta},\\
\displaystyle  \left<\partial_t \mathbf{F}_m, \boldsymbol{\Sigma} \right>+\int_{\Omega}\left({\mathbf{v}}_m \cdot \nabla \right)\mathbf{F}_m\colon \boldsymbol{\Sigma}-\int_{\Omega}\left(\nabla {\mathbf{v}}_m\right)\mathbf{F}_m\colon \boldsymbol{\Sigma}+\lambda \int_{\Omega}\nabla \left(f({\phi}_m)\mathbf{F}_m\right)\mathbin{\tensorm} \nabla \boldsymbol{\Sigma}=0,\\
 \displaystyle \int_{\Omega}\partial_t\phi_m \xi+\int_{\Omega}\left(\mathbf{v}_m \cdot \nabla \phi_m\right) \xi+\int_{\Omega}b\left(\phi_m\right)\nabla \mu_m \cdot \nabla \xi=0,\\
\displaystyle \int_{\Omega}\mu_m \chi=\int_{\Omega}\nabla \phi_m\cdot \nabla \chi+\int_{\Omega}\psi'\left(\phi_m\right)\chi+\int_{\Omega}\frac{f'\left(\phi_m\right)}{2}\left(\mathbf{F}_m\colon \mathbf{F}_m-2\right)\chi+\delta \int_{\Omega}\partial_t\phi_m \chi,
\end{cases}
\end{equation} 
for a.e. $t \in [0,T]$, with $\boldsymbol{\eta} \in \text{span}\left\{\boldsymbol{\eta}_0, \boldsymbol{\eta}_1, \dots, \boldsymbol{\eta}_m\right\}$, $\xi, \chi \in \text{span}\left\{\xi_0, \xi_1, \dots, \xi_m\right\}$ and $\forall \, \boldsymbol{\Sigma} \in H^1(\Omega;\mathbb{R}^{2\times 2})$.

We now study the convergence properties of the discrete solution of system \eqref{eqn:galerkin3}, identifying its limit as the discretization parameter $m$ tends to infinity, {\rosso starting from \eqref{eqn:apriori}, which is valid also for the limit case $\gamma,\theta \to 0$. Note that in the latter case the constant in the right hand side of \eqref{eqn:apriori} depends only on the initial data $\mathbf{F}_0$ and $\phi_0$}. Since this study will be different in the cases of $2$ or $3$ space dimensions, we separate the two cases in two different paragraphs.
\subsection{Limit problem in the $2$-d case}
We derive the necessary regularity results to study the limit point of \eqref{eqn:galerkin3}, as $m\to \infty$ and $\delta \to 0$. From the Poincar\'e inequality and from \eqref{eqn:apriori} (with $\gamma=0$), we have that
{\rosso
\begin{equation}
\label{eqn:vm}
\mathbf{v}_m\quad \text{is uniformly bounded in} \; L^2(0,T;H_{0,\diver}^1(\Omega;\mathbb{R}^2))\hookrightarrow L^2(0,T;L_{\diver}^p(\Omega;\mathbb{R}^2)),
\end{equation}
with $p\geq 1$. Moreover, from \eqref{eqn:apriori} and Assumption \textbf{A3}, we have that
\begin{equation}
\label{eqn:fFm}
f(\phi_m)\mathbf{F}_m \quad \text{is uniformly bounded in} \; L^{\infty}(0,T;L^2(\Omega;\mathbb{R}^{2\times 2}))\cap L^{2}(0,T;H^1(\Omega;\mathbb{R}^{2\times 2})),
\end{equation}
which implies
\begin{equation}
\label{eqn:Fm}
||\mathbf{F}_m||_{L^{\infty}(0,T;L^2(\Omega;\mathbb{R}^{2\times 2}))}\leq C.
\end{equation}
We recall \eqref{eqn:phimh1} and take} $\chi=-\Delta \phi_m$ in \eqref{eqn:galerkin3}$_4$. Integrating in time over the interval $[0,T]$ and using the Cauchy--Schwarz and Young inequalities and Assumptions \textbf{A2} and \textbf{A3}, we get
\begin{align}
\label{eqn:deltaphi1}
& \int_0^T||\Delta \phi_m||^2+\int_0^T\psi_+''(\phi_m)|\nabla \phi_m|^2+\frac{\delta}{2}||\nabla \phi_m(T)||^2\leq \int_0^T||\nabla \mu_m||_{L^2(\Omega;\mathbb{R}^2)}||\nabla \phi_m||_{L^2(\Omega;\mathbb{R}^2)}\\
& \notag -\int_0^T\psi_-''(\phi_m)|\nabla \phi_m|^2+\int_0^T\int_{\Omega}\frac{f'\left(\phi_m\right)}{2}\left(\mathbf{F}_m\colon \mathbf{F}_m-2\right)\Delta \phi_m+\frac{\delta}{2}||\nabla \phi_m(0)||^2\\
& \notag  \leq \frac{1}{2}\int_0^T||\nabla \mu_m||_{L^2(\Omega;\mathbb{R}^2)}^2+\frac{1}{2}\int_0^T||\nabla \phi_m||_{L^2(\Omega;\mathbb{R}^2)}^2+C\int_0^T\int_{\Omega}\left(1+|\phi_m|^q\right)|\nabla \phi_m|^2\\
& \notag +\frac{f_1}{2}\int_0^T||\Delta \phi_m||\, ||\mathbf{F}_m||_{L^4(\Omega;\mathbb{R}^{2\times 2})}^2+\frac{1}{6}\int_0^T||\Delta \phi_m||^2+C.
\end{align} 
We use \eqref{eqn:agmon2d}, elliptic regularity and the Young inequality, observing that $\frac{2}{\theta q}>1$ when $q<\frac{2}{\theta}$, which complies with Assumption \textbf{A2} since $\theta$ can be taken arbitrarily small, to write
\begin{align*}
&\int_0^T\int_{\Omega}|\phi_m|^q|\nabla \phi_m|^2\leq \int_0^T||\phi_m||_{L^{\infty}(\Omega)}^q||\nabla \phi_m||_{L^2(\Omega;\mathbb{R}^2)}^2 \\
& \leq C\int_0^T||\phi_m||_{H^1(\Omega;\mathbb{R}^2)}^{2+q(1-\theta)}\left(||\phi_m||^{q\theta}+||\Delta \phi_m||^{q\theta}\right)\\
& \leq C||\phi_m||_{L^{\infty}(0,T;H^1(\Omega;\mathbb{R}^2))}^{2+q}+||\phi_m||_{L^{\infty}(0,T;H^1(\Omega;\mathbb{R}^2))}^{\frac{2(2+q(1-\theta))}{2-q\theta}}+\frac{1}{6}\int_0^T||\Delta \phi_m||^2.
\end{align*}
Choosing $h=2$ in \eqref{eqn:gninterpolation}, then
{\rosso
\begin{align}
\label{eqn:Fm4}
& f(\phi_m)\mathbf{F}_m \quad \text{is uniformly bounded in }\\
& \notag L^{\infty}(0,T;L^2(\Omega;\mathbb{R}^{2\times 2}))\cap L^{2}(0,T;H^1(\Omega;\mathbb{R}^{2\times 2}))\hookrightarrow L^{4}(0,T;L^4(\Omega;\mathbb{R}^{2\times 2})),\\
& \notag \text{which implies} \quad ||\mathbf{F}_m||_{L^{4}(0,T;L^4(\Omega;\mathbb{R}^{2\times 2}))}\leq C. 
\end{align}
}%
Using the preceding results in \eqref{eqn:deltaphi1}, using moreover the Young inequality, \eqref{eqn:apriori} and Assumption \textbf{A1}, we obtain that
\begin{equation*}
\int_0^T||\Delta \phi_m||^2\leq C+C||\mathbf{F}_m||_{L^{4}(0,T;L^4(\Omega;\mathbb{R}^{2\times 2}))}^4+\frac{1}{2}\int_0^T||\Delta \phi_m||^2,
\end{equation*}
which implies that
\begin{equation}
\label{eqn:deltaphi}
\int_0^T||\Delta \phi_m||^2\leq C.
\end{equation}
Hence, {\rosso from \eqref{eqn:phimh1} and \eqref{eqn:deltaphi}} we obtain that
{\rosso
\begin{equation}
\label{eqn:phim}
\phi_m\quad \text{is uniformly bounded in } \, L^{\infty}(0,T;H^1(\Omega))\cap L^{2}(0,T;H^2(\Omega)).
\end{equation}
}%
Note that, thanks again to \eqref{eqn:gninterpolation} with $h=2$,
{\rosso
\begin{equation}
\label{eqn:nablaphim4}
||\nabla \phi_m||_{L^{4}(0,T;L^4(\Omega;\mathbb{R}^2))}\leq C.
\end{equation}
}
In the $2$-d case, thanks to \eqref{eqn:Fm4} and \eqref{eqn:nablaphim4}, $\sqrt{f(\phi_m)}\mathbf{F}_m\in L^{2}(0,T;H^1(\Omega;\mathbb{R}^{2\times 2}))$ implies that 
{\rosso
\begin{equation}
\label{eqn:gradFm}
||\mathbf{F}_m||_{L^{2}(0,T;H^1(\Omega;\mathbb{R}^{2\times 2}))}\leq C.
\end{equation}
}
Indeed, writing
\begin{align*}
&\int_0^T\int_{\Omega}\nabla \left(f({\phi}_m)\mathbf{F}_m\right)\mathbin{\tensorm} \nabla \left(f({\phi}_m)\mathbf{F}_m\right)=\int_0^T\int_{\Omega}f({\phi}_m)^2\left|\nabla \mathbf{F}_m\right|^2+\int_0^T\int_{\Omega}f'({\phi}_m)^2\left|\mathbf{F}_m\otimes \nabla \phi_m\right|^2\\
&+2\int_0^T\int_{\Omega}f({\phi}_m)f'({\phi}_m)\nabla \mathbf{F}_m\mathbin{\tensorm} \left(\mathbf{F}_m\otimes \nabla \phi_m\right),
\end{align*}
and, due to \eqref{eqn:apriori}, we have that
\begin{align*}
&\int_0^T\int_{\Omega}f({\phi}_m)^2\left|\nabla \mathbf{F}_m\right|^2+\int_0^T\int_{\Omega}f'({\phi}_m)^2\left|\mathbf{F}_m\otimes \nabla \phi_m\right|^2\\
& +2\int_0^T\int_{\Omega}f({\phi}_m)f'({\phi}_m)\nabla \mathbf{F}_m\mathbin{\tensorm} \left(\mathbf{F}_m\otimes \nabla \phi_m\right)\leq C,
\end{align*}
from which we infer, using Assumption \textbf{A3}, \eqref{eqn:Fm4} and \eqref{eqn:nablaphim4}, the Cauchy-Schwarz and Young inequalities,
\begin{align}
\label{eqn:gradFmcalccalc}
&f_{\text{min}}\int_0^T\int_{\Omega}\left|\nabla \mathbf{F}_m\right|^2\leq \int_0^T\int_{\Omega}f({\phi}_m)^2\left|\nabla \mathbf{F}_m\right|^2 \\
& \notag \leq C+2f_1\int_0^T\left|\left|\nabla \mathbf{F}_m\right|\right|_{L^2(\Omega;\mathbb{R}^{2\times 2\times 2})}\left|\left|\mathbf{F}_m\otimes \nabla \phi_m\right|\right|_{L^2(\Omega;\mathbb{R}^{2\times 2\times 2})}\\
& \notag \leq C+\frac{f_{\text{min}}}{2}\int_0^T\int_{\Omega}\left|\nabla \mathbf{F}_m\right|^2+C\int_0^T\left|\left|\mathbf{F}_m\right|\right|_{L^4(\Omega;\mathbb{R}^{2\times 2})}^2\left|\left|\nabla \phi_m\right|\right|_{L^4(\Omega;\mathbb{R}^{2})}^2\\
& \notag \leq C+\frac{f_{\text{min}}}{2}\int_0^T\int_{\Omega}\left|\nabla \mathbf{F}_m\right|^2+C\left|\left|\mathbf{F}_m\right|\right|_{L^{4}\left(0,T;L^4(\Omega;\mathbb{R}^{2\times 2})\right)}^2\left|\left|\nabla \phi_m\right|\right|_{L^4\left(0,T;L^4(\Omega;\mathbb{R}^{2})\right)}^2,
\end{align}
which implies that
\begin{equation}
\label{eqn:gradFmcalc}
\frac{f_{\text{min}}}{2}\int_0^T\int_{\Omega}\left|\nabla \mathbf{F}_m\right|^2\leq C.
\end{equation}

We now deduce a priori estimates, uniform in $m$ and $\delta$, for the time derivatives of $\mathbf{F}_m$ and $\phi_m$ in \eqref{eqn:galerkin3}, which are needed to obtain compactness results for the sequence of discrete solutions of the system. We start from the second equation of \eqref{eqn:galerkin3}. We rewrite the second term on the left hand side, integrating by parts and considering the properties of the variable $\mathbf{v}_m$, as
\begin{equation}
\label{eqn:fmparts}
\int_{\Omega}\left({\mathbf{v}}_m \cdot \nabla \right)\mathbf{F}_m\colon \boldsymbol{\Sigma}=\underbrace{\int_{\partial \Omega}\left(\mathbf{F}_m\colon \boldsymbol{\Sigma}\right)\mathbf{v}_m\cdot \mathbf{n}}_{=0}- \underbrace{\int_{\Omega}\diver \mathbf{v}_m\mathbf{F}_m\colon \boldsymbol{\Sigma}}_{=0}-\int_{\Omega}\mathbf{F}_m\otimes \mathbf{v}_m\mathbin{\tensorm} \nabla \boldsymbol{\Sigma},
\end{equation}
which is valid, for a fixed $m$, thanks to \eqref{eqn:fspace} and the regularity of $\mathbf{v}_m$. 
We multiply \eqref{eqn:fmparts} by a time function $\zeta \in L^{2+h}(0,T)$, with $h>0$, and integrate in time over the interval $(0,T)$, to infer that
\begin{align}
\label{eqn:mfdualcalc1}
&\int_0^T\int_{\Omega}\left({\mathbf{v}}_m \cdot \nabla \right)\mathbf{F}_m\colon \boldsymbol{\Sigma}\zeta=-\int_0^T\int_{\Omega}\mathbf{F}_m\otimes \mathbf{v}_m\mathbin{\tensorm} \nabla \boldsymbol{\Sigma}\zeta \\
& \notag \leq \int_0^T ||\mathbf{F}_m||_{L^{2+h}(\Omega;\mathbb{R}^{2\times 2})}||\mathbf{v}_m||_{L^{\frac{2(2+h)}{h}}(\Omega;\mathbb{R}^{2})}||\nabla \boldsymbol{\Sigma}||_{L^{2}(\Omega;\mathbb{R}^{2\times 2 \times 2})}|\zeta| \\
& \notag \leq C||\mathbf{F}_m||_{L^{\frac{2(2+h)}{h}}\left(0,T;L^{2+h}(\Omega;\mathbb{R}^{2\times 2})\right)}||\mathbf{v}_m||_{L^2\left(0,T;L^{\frac{2(2+h)}{h}}(\Omega;\mathbb{R}^{2})\right)}||\nabla \boldsymbol{\Sigma}||_{L^{2}(\Omega;\mathbb{R}^{2\times 2 \times 2})}||\zeta||_{L^{2+h}(0,T)}.
\end{align}
Then, thanks to \eqref{eqn:gninterpolation} and Sobolev embedding, we have that
{\rosso
\begin{equation}
\label{eqn:fmdual1}
\left({\mathbf{v}}_m \cdot \nabla \right)\mathbf{F}_m\quad \text{is uniformly bounded in}\, L^{2-\frac{h}{1+h}}\left(0,T;\left(H^1(\Omega;\mathbb{R}^{2\times 2})\right)'\right),
\end{equation}}%
for arbitrary $h>0$.
Using \eqref{eqn:fmdual1}, multiplying \eqref{eqn:galerkin3}$_2$ by a time function $\zeta \in L^{2+h}(0,T)$, with $h>0$, and integrating in time over the interval $(0,T)$, we get
\begin{align}
\label{eqn:mfdualcalc2}
&\int_0^T\left<\partial_t \mathbf{F}_m, \boldsymbol{\Sigma} \right>\zeta \leq \int_0^T ||\left({\mathbf{v}}_m \cdot \nabla \right)\mathbf{F}_m||_{\left(H^1(\Omega;\mathbb{R}^{2\times 2})\right)'}||\boldsymbol{\Sigma}||_{H^{1}(\Omega;\mathbb{R}^{2\times 2})}|\zeta|\\
& \notag +\int_0^T ||\nabla \mathbf{v}_m||_{L^{2}(\Omega;\mathbb{R}^{2})}||\mathbf{F}_m||_{L^{2+h}(\Omega;\mathbb{R}^{2\times 2})}||\boldsymbol{\Sigma}||_{L^{\frac{2(2+h)}{h}}(\Omega;\mathbb{R}^{2\times 2})}|\zeta|\\
& \notag +\lambda\int_0^T||\nabla \left(f({\phi}_m)\mathbf{F}_m\right)||_{L^{2}(\Omega;\mathbb{R}^{2\times 2\times 2})}||\nabla \boldsymbol{\Sigma}||_{L^{2}(\Omega;\mathbb{R}^{2\times 2\times 2})}|\zeta|\leq C || \boldsymbol{\Sigma}||_{H^{1}(\Omega;\mathbb{R}^{2\times 2})}||\zeta||_{L^{2+h}(0,T)}\\
& \notag +C||\nabla \mathbf{v}_m||_{L^2\left(0,T;L^{2}(\Omega;\mathbb{R}^{2\times 2})\right)}||\mathbf{F}_m||_{L^{\frac{2(2+h)}{h}}\left(L^{2+h}(\Omega;\mathbb{R}^{2\times 2})\right)}||\boldsymbol{\Sigma}||_{H^{1}(\Omega;\mathbb{R}^{2\times 2})}||\zeta||_{L^{2+h}(0,T)}.
\end{align}
Hence, we have that
{\rosso
\begin{equation}
\label{eqn:fmdual2}
\partial_t \mathbf{F}_m \quad \text{is uniformly bounded in}\; L^{2-\frac{h}{1+h}}\left(0,T;\left(H^1(\Omega;\mathbb{R}^{2\times 2})\right)'\right).
\end{equation}}%
We now consider the third equation of \eqref{eqn:galerkin3}. For what concerns the advective term, thanks to \eqref{eqn:phim} and \eqref{eqn:vm}, we have the estimate
\begin{equation}
\label{eqn:phimdualcalc1}
\int_0^T\left(\int_{\Omega}\left|\mathbf{v}_m \cdot \nabla \phi_m\right|^{2-h}\right)^{\frac{2}{2-h}}\leq \int_0^T||\mathbf{v}_m||_{L^{\frac{2(2-h)}{h}}(\Omega;\mathbb{R}^2)}^2||\nabla \phi_m||_{L^2(\Omega;\mathbb{R}^2)}^2\leq C,
\end{equation}
hence
\begin{equation}
\label{eqn:phimdual1}
||\mathbf{v}_m \cdot \nabla \phi_m||_{L^{2}\left(0,T;L^{2-h}(\Omega)\right)}\leq C,
\end{equation}
for $h\in (0,1]$.
Multiplying \eqref{eqn:galerkin3}$_3$ by a time function $\zeta \in L^{2}(0,T)$, integrating in time over the interval $(0,T)$ and choosing $\xi=P_m^L(\pi)$, with generic $\pi \in H^1(\Omega)$, using \eqref{eqn:phimdual1}, \eqref{eqn:apriori} and Assumption \textbf{A1}, we get
\begin{align*}
& \int_0^T\int_{\Omega}\partial_t\phi_m \pi \zeta = \int_0^T\int_{\Omega}\partial_t\phi_m \left(\pi+ (P_m^L(\pi)-\pi)\right) \zeta \\
&\leq \int_0^T||\mathbf{v}_m \cdot \nabla \phi_m||_{L^{2-h}(\Omega)}||P_m^L(\pi)||_{L^{2+\frac{h}{1-h}}(\Omega)}|\zeta|\\
& +\int_0^T||b\left(\phi_m\right)\nabla \mu_m||_{L^{2}(\Omega;\mathbb{R}^2)} ||\nabla P_m^L(\pi)||_{L^{2}(\Omega;\mathbb{R}^2)}|\zeta|\leq C||\pi||_{H^1(\Omega)}||\zeta||_{L^2(0,T)}.
\end{align*}
Hence, we have that
\begin{equation}
\label{eqn:phimdual2}
||\partial_t \phi_m||_{L^{2}\left(0,T;\left(H^1(\Omega)\right)'\right)}\leq C.
\end{equation}
Collecting the estimates \eqref{eqn:vm}, \eqref{eqn:fFm}, \eqref{eqn:Fm}, \eqref{eqn:gradFm}, \eqref{eqn:phim}, \eqref{eqn:mumh1},  \eqref{eqn:fmdual2} and \eqref{eqn:phimdual2}, which are uniform in $m$ and $\delta$, from the Banach--Alaoglu and the Aubin--Lions lemma, we finally obtain the convergence properties, up to subsequences of the solutions, which we still label by the index $m$,
\begin{align}
\label{eqn:conv2d1} & \mathbf{v}_m \rightharpoonup \mathbf{v} \quad \text{in} \quad L^{2}\left(0,T;H_{0,\diver}^1\left(\Omega;\mathbb{R}^{2}\right)\right),\\
\label{eqn:conv2d2} & f(\phi_m)\mathbf{F}_m \overset{\ast}{\rightharpoonup} \boldsymbol{\beta} \quad \text{in} \quad L^{\infty}(0,T;L^2(\Omega;\mathbb{R}^{2\times 2}))\cap L^{2}(0,T;H^1(\Omega;\mathbb{R}^{2\times 2})),\\
\label{eqn:conv2d3} & \mathbf{F}_m \overset{\ast}{\rightharpoonup} \mathbf{F} \quad \text{in} \quad L^{\infty}(0,T;L^2(\Omega;\mathbb{R}^{2\times 2}))\cap L^{2}(0,T;H^1(\Omega;\mathbb{R}^{2\times 2})),\\
\label{eqn:conv2d4} & \partial_t \mathbf{F}_m \rightharpoonup \partial_t \mathbf{F} \quad \text{in} \quad L^{2-h}\left(0,T;\left(H^1(\Omega;\mathbb{R}^{2\times 2})\right)'\right), \quad {\rosso h\in (0,1)},\\
\label{eqn:conv2d5} & \mathbf{F}_m \to \mathbf{F} \quad \text{in} \quad L^{2}(0,T;L^p(\Omega;\mathbb{R}^{2\times 2})), \quad p\geq 1, \quad \text{and} \quad \text{a.e. in} \; \; \Omega_T,\\
\label{eqn:conv2d6} & \phi_m \overset{\ast}{\rightharpoonup} \phi \quad \text{in} \quad L^{\infty}(0,T;H^1(\Omega))\cap L^{2}(0,T;H^2(\Omega)),\\
\label{eqn:conv2d7} & \partial_t \phi_m \rightharpoonup \partial_t \phi \quad \text{in} \quad L^{2}\left(0,T;\left(H^1(\Omega)\right)'\right),\\
\label{eqn:conv2d8} & \phi_m \to \phi \quad \text{in} \quad C^{0}(0,T;L^p(\Omega)) \cap L^{2}(0,T;W^{1,p}(\Omega)), \quad p\geq 1, \quad \text{and} \quad \text{a.e. in} \; \; \Omega_T,\\
\label{eqn:conv2d9} & \mu_m \rightharpoonup \mu \quad \text{in} \quad L^{2}\left(0,T;H^1\left(\Omega\right)\right),
\end{align}
as $m\to \infty$, $\delta \to 0$, for some limit function $\boldsymbol{\beta}$ in $L^{\infty}(0,T;L^2(\Omega;\mathbb{R}^{2\times 2}))\cap L^{2}(0,T;H^1(\Omega;\mathbb{R}^{2\times 2}))$, which may be identified as $\boldsymbol{\beta}=f(\phi)\mathbf{F}$ thanks to \eqref{eqn:conv2d5} and \eqref{eqn:conv2d8} (see the later calculations). Moreover, from \eqref{eqn:apriori}, we have the convergence result
\begin{equation}
\label{eqn:conv2d10}\sqrt{\delta} \partial_t \phi_m \rightharpoonup \alpha \quad \text{in} \quad L^{2}(\Omega_T),
\end{equation}
as $m\to \infty$, $\delta \to 0$, for some function $\alpha \in L^{2}(\Omega_T)$, {\rosso which should be zero as a consequence of \eqref{eqn:conv2d7}}.

With the convergence results \eqref{eqn:conv2d1}--\eqref{eqn:conv2d9}, we can pass to the limit in the system \eqref{eqn:galerkin3} as $m\to \infty$ and $\delta\to 0$. Let's take $\boldsymbol{\eta}=\boldsymbol{\eta}_m=P_m^{S}(\mathbf{u})$, for arbitrary $\mathbf{u}\in H_{0,\diver}^1(\Omega;\mathbb{R}^2)$, $\xi=\xi_m=P_m^L(q)$, for arbitrary $q\in H^1(\Omega)$, $\chi=\chi_m=P_m^L(r)$, for arbitrary $r\in H^1(\Omega)$, multiply the equations by a function $\omega \in C_0^{\infty}([0,T])$ and integrate over the time interval $[0,T]$. This gives
\begin{equation}
\label{eqn:limit2d1}
\begin{cases}
\displaystyle \nu \int_0^T\omega \int_{\Omega}\nabla \mathbf{v}_m\colon \nabla \boldsymbol{\eta}_m= \int_0^T\omega \int_{\Omega}\mu_m \nabla \phi_m \cdot \boldsymbol{\eta}_m-\int_0^T\omega\int_{\Omega}\frac{f'\left(\phi_m\right)}{2}\left(\mathbf{F}_m\colon \mathbf{F}_m-2\right)\nabla \phi_m\cdot \boldsymbol{\eta}_m\\ 
\displaystyle -\int_0^T\omega\int_{\Omega}f\left(\phi_m\right)\mathbf{F}_m\mathbf{F}_m^T\colon \nabla \boldsymbol{\eta}_m,\\ \\
\displaystyle  \int_0^T\omega\left<\partial_t \mathbf{F}_m, \boldsymbol{\Sigma} \right>+\int_0^T\omega\int_{\Omega}\left({\mathbf{v}}_m \cdot \nabla \right)\mathbf{F}_m\colon \boldsymbol{\Sigma}-\int_0^T\omega\int_{\Omega}\left(\nabla {\mathbf{v}}_m\right)\mathbf{F}_m\colon \boldsymbol{\Sigma}\\ 
\displaystyle +\lambda \int_0^T\omega\int_{\Omega}\nabla \left(f({\phi}_m)\mathbf{F}_m\right)\mathbin{\tensorm} \nabla \boldsymbol{\Sigma}=0,\\ \\
 \displaystyle \int_0^T\omega\int_{\Omega}\partial_t\phi_m \xi_m+\int_0^T\omega\int_{\Omega}\left(\mathbf{v}_m \cdot \nabla \phi_m\right) \xi_m+\int_0^T\omega\int_{\Omega}b\left(\phi_m\right)\nabla \mu_m \cdot \nabla \xi_m=0,\\ \\
\displaystyle \int_0^T\omega\int_{\Omega}\mu_m \chi_m=\int_0^T\omega\int_{\Omega}\nabla \phi_m\cdot \nabla \chi_m+\int_0^T\omega\int_{\Omega}\psi'\left(\phi_m\right)\chi_m\\
\displaystyle +\int_0^T\omega\int_{\Omega}\frac{f'\left(\phi_m\right)}{2}\left(\mathbf{F}_m\colon \mathbf{F}_m-2\right)\chi_m+\delta \int_0^T\omega\int_{\Omega}\partial_t\phi_m \chi_m.
\end{cases}
\end{equation} 
Thanks to \eqref{eqn:conv2d1}, we have
\[
\nu \int_0^T\omega \int_{\Omega}\nabla \mathbf{v}_m\colon \nabla \boldsymbol{\eta}_m\to \nu \int_0^T\omega \int_{\Omega}\nabla \mathbf{v}\colon \nabla \mathbf{u},
\]
as $m\to \infty$, $\delta \to 0$. Thanks to \eqref{eqn:conv2d8} and the fact that $P_m^{S}(\mathbf{u})$ converges strongly to $\mathbf{u}$ in $H^1(\Omega;\mathbb{R}^2)$, we have that $\omega \nabla \phi_m \cdot \boldsymbol{\eta}_m\to \omega \nabla \phi \cdot \mathbf{u}$ strongly in $L^2(0,T;L^p(\Omega))$, with arbitrary $p\geq 1$. Indeed,
\begin{align}
\label{eqn:limitcalc1}
& \int_0^T\left(\int_{\Omega}|\omega \nabla \phi_m \cdot \boldsymbol{\eta}_m - \omega \nabla \phi \cdot \mathbf{u}|^p\right)^{\frac{2}{p}}=\int_0^T\left(\int_{\Omega}|\omega (\nabla \phi_m -\nabla \phi) \cdot \boldsymbol{\eta}_m + \omega \nabla \phi \cdot (\boldsymbol{\eta}_m - \mathbf{u})|^p\right)^{\frac{2}{p}}\\
& \notag \leq C \int_0^T|\omega|^2\,||\nabla \phi_m -\nabla \phi||_{L^{2p}(\Omega;\mathbb{R}^2)}^2||\boldsymbol{\eta}_m||_{L^{2p}(\Omega;\mathbb{R}^2)}^2\\
& \notag +C \int_0^T|\omega|^2\,||\nabla \phi||_{L^{2p}(\Omega;\mathbb{R}^2)}^2||\boldsymbol{\eta}_m-\mathbf{u}||_{L^{2p}(\Omega;\mathbb{R}^2)}^2 \rightarrow 0,
\end{align}
as $m\to \infty$, $\delta \to 0$.
Hence, using \eqref{eqn:conv2d9}, which implies that $\mu_m \rightharpoonup \mu \quad \text{in} \quad L^{2}\left(0,T;L^p\left(\Omega\right)\right)$, with arbitrary $p\geq 1$, by the product of weak--strong convergence we have that 
\[
\int_0^T\omega\int_{\Omega}\mu_m \nabla \phi_m \cdot \boldsymbol{\eta}_m\to \int_0^T\omega\int_{\Omega}\mu \nabla \phi \cdot \mathbf{u},
\] 
as $m\to \infty$, $\delta \to 0$. For what concerns the second term on the right hand side of the first equation of \eqref{eqn:limit2d1}, we observe that $f'\left(\phi_m\right)\mathbf{F}_m\colon \mathbf{F}_m\to f'\left(\phi\right)\mathbf{F}\colon \mathbf{F}$ strongly in $L^2(0,T;L^{2-h}(\Omega))$, for arbitrary $h \in (0,1]$. Indeed, using an inequality similar to \eqref{eqn:phimdualcalc1}, using moreover Assumption \textbf{A3} and \eqref{eqn:ffid1}, we get
\begin{align*}
& \int_0^T\left(\int_{\Omega}\left|f'\left(\phi_m\right)\mathbf{F}_m\colon \mathbf{F}_m - f'\left(\phi\right)\mathbf{F}\colon \mathbf{F}\right|^{2-h}\right)^{\frac{2}{2-h}}\\
& =\int_0^T\left(\int_{\Omega}\left|f'\left(\phi_m\right)\left((\mathbf{F}_m-\mathbf{F})\colon (\mathbf{F}_m+\mathbf{F})\right)+\left(f'\left(\phi_m\right)-f'\left(\phi \right)\right)\mathbf{F}\colon \mathbf{F}\right|^{2-h}\right)^{\frac{2}{2-h}} \\
&\leq f_1\int_0^T||\mathbf{F}_m-\mathbf{F}||_{L^{\frac{2(2-h)}{h}}(\Omega;\mathbb{R}^{2\times 2})}^2||\mathbf{F}_m+\mathbf{F}||_{L^2(\Omega;\mathbb{R}^{2\times 2})}^2\\
&+c_f\int_0^T||\phi_m-\phi||_{L^{\frac{2(2-h)}{h}}(\Omega)}^2||\mathbf{F}\colon \mathbf{F}||_{L^2(\Omega;\mathbb{R}^{2\times 2})}^2\to 0,
\end{align*}
as $m\to \infty$, $\delta \to 0$, where $c_f$ is the Lipschitz constant of $f'$ and where we used \eqref{eqn:conv2d3}, \eqref{eqn:conv2d5}, \eqref{eqn:conv2d8} and \eqref{eqn:Fm4}, which implies that $\mathbf{F}\colon \mathbf{F}\in L^2\left(0,t;L^2(\Omega;\mathbb{R}^{2\times 2})\right)$. Then, considering the strong convergence of $\boldsymbol{\eta}_m=P_m^{S}(\mathbf{u})$ to $\mathbf{u}$ in $H^1(\Omega)$ and making similar calculations to \eqref{eqn:limitcalc1}, we have that $f'\left(\phi_m\right)\left(\mathbf{F}_m\colon \mathbf{F}_m-2\right)\boldsymbol{\eta}_m\to f'\left(\phi\right)\left(\mathbf{F}\colon \mathbf{F}-2\right)\mathbf{u}$ strongly in $L^2(0,T;L^{2-h_1}(\Omega;\mathbb{R}^2))$, with $0<h<h_1\leq 1$, where we also used Assumption \textbf{A3}, \eqref{eqn:conv2d8} and the Lebesgue convergence theorem to obtain $f'\left(\phi_m\right)\to f'\left(\phi\right)$ strongly in $L^{2}(0,T;L^2(\Omega))$. Hence, using \eqref{eqn:conv2d6} and the product of weak--strong convergence we have that 
\[
\int_0^T\omega\int_{\Omega}\frac{f'\left(\phi_m\right)}{2}\left(\mathbf{F}_m\colon \mathbf{F}_m-2\right)\nabla \phi_m\cdot \boldsymbol{\eta}_m\to \int_0^T\omega\int_{\Omega}\frac{f'\left(\phi\right)}{2}\left(\mathbf{F}\colon \mathbf{F}-2\right)\nabla \phi\cdot \mathbf{u},
\] 
as $m\to \infty$, $\delta \to 0$. For what concerns the third term on the right hand side of the first equation of \eqref{eqn:limit2d1}, we observe that, thanks to Assumption \textbf{A3} and to \eqref{eqn:conv2d8}, $ f(\phi_m)\to f(\phi) $ a.e. in $\Omega_T$ and is uniformly bounded. Hence, due to \eqref{eqn:conv2d5}, a generalized version of the Lebesgue convergence theorem {\rosso (see e.g. \cite[Theorem 3.25, p. 60]{alt})} gives that $f(\phi_m)\mathbf{F}_m\to f(\phi)\mathbf{F}$ strongly in $L^2(0,T;L^{p}(\Omega;\mathbb{R}^{2\times 2}))$, for arbitrary $p\geq 1$, and a.e. in $\Omega_T$. In particular, thanks to the uniqueness of the weak limit and the a.e. convergence,
\begin{equation}
\label{eqn:beta}
\boldsymbol{\beta}=f(\phi)\mathbf{F}.
\end{equation}
Then, considering again the strong convergence of $\boldsymbol{\eta}_m=P_m^{S}(\mathbf{u})$ to $\mathbf{u}$ in $H^1(\Omega)$, we have that \\
$f(\phi_m)\nabla \boldsymbol{\eta}_m\mathbf{F}_m\to f(\phi)\nabla \mathbf{u}\mathbf{F}$ strongly in $L^2(0,T;L^{2-h}(\Omega;\mathbb{R}^{2\times 2}))$, for arbitrary $h\in (0,1]$. Hence, using \eqref{eqn:conv2d3} and the product of weak--strong convergence we have that 
\[
\int_0^T\omega\int_{\Omega}f\left(\phi_m\right)\mathbf{F}_m\mathbf{F}_m^T\colon \nabla \boldsymbol{\eta}_m \to \int_0^T\omega\int_{\Omega}f\left(\phi\right)\mathbf{F}\mathbf{F}^T\colon \nabla \mathbf{u},
\]
as $m\to \infty$, $\delta \to 0$.

Now, if we consider the second equation in \eqref{eqn:limit2d1}, since
\[\omega \boldsymbol{\Sigma}\in C^0(0,T;H^1(\Omega;\mathbb{R}^{2\times 2}))\hookrightarrow L^{2+h}\left(0,T;H^1(\Omega;\mathbb{R}^{2\times 2})\right),\] 
we get from \eqref{eqn:conv2d4} that
\[
\int_0^T\omega\left<\partial_t \mathbf{F}_m, \boldsymbol{\Sigma} \right>\to \int_0^T\omega\left<\partial_t \mathbf{F}, \boldsymbol{\Sigma} \right>,
\]
as $m\to \infty$, $\delta \to 0$. Due to the fact that we only have weak convergence for both $\mathbf{v}_m$ and $\nabla \mathbf{F}_m$ in \eqref{eqn:conv2d1}--\eqref{eqn:conv2d5}, we rewrite the second term in the second equation of \eqref{eqn:limit2d1} employing the formula \eqref{eqn:fmparts}. Thanks to \eqref{eqn:conv2d5} and the fact that $\nabla \boldsymbol{\Sigma}\in L^2(\Omega;\mathbb{R}^{2\times 2\times 2})$, we have that $\omega \left(\nabla \boldsymbol{\Sigma}\right)^T\mathbf{F}_m \to \omega \left(\nabla \boldsymbol{\Sigma}\right)^T\mathbf{F}$ strongly in $L^2(0,T;L^{2-h}(\Omega;\mathbb{R}^{2}))$, for arbitrary $h\in (0,1]$. Hence, using \eqref{eqn:conv2d1} and the product of weak--strong convergence we have that 
\[
\int_0^T\omega\int_{\Omega}\left({\mathbf{v}}_m \cdot \nabla \right)\mathbf{F}_m\colon \boldsymbol{\Sigma}=-\int_0^T\omega \int_{\Omega}\mathbf{F}_m\otimes \mathbf{v}_m\mathbin{\tensorm} \nabla \boldsymbol{\Sigma}\to -\int_0^T\omega \int_{\Omega}\mathbf{F}\otimes \mathbf{v}\mathbin{\tensorm} \nabla \boldsymbol{\Sigma},
\]
as $m \to \infty$, $\delta \to 0$.
Moreover, thanks to \eqref{eqn:conv2d5} and the fact that $\boldsymbol{\Sigma}\in H^1(\Omega;\mathbb{R}^{2\times 2})\hookrightarrow L^p(\Omega;\mathbb{R}^{2\times 2})$, for arbitrary $p\geq 1$, we have that $\omega\boldsymbol{\Sigma}\mathbf{F}_m^T \to \omega\boldsymbol{\Sigma}\mathbf{F}^T$ strongly in $L^2(0,T;L^{p}(\Omega;\mathbb{R}^{2}))$, for arbitrary $p\geq 1$. Hence, using \eqref{eqn:conv2d1} and the product of weak--strong convergence we have that 
\[
\int_0^T\omega\int_{\Omega}\left(\nabla {\mathbf{v}}_m\right)\mathbf{F}_m\colon \boldsymbol{\Sigma}\to \int_0^T\omega\int_{\Omega}\left(\nabla {\mathbf{v}}\right)\mathbf{F}\colon \boldsymbol{\Sigma},
\]
as $m \to \infty$, $\delta \to 0$. Finally, given \eqref{eqn:conv2d2} and \eqref{eqn:beta},
\[
\lambda \int_0^T\omega\int_{\Omega}\nabla \left(f({\phi}_m)\mathbf{F}_m\right)\mathbin{\tensorm} \nabla \boldsymbol{\Sigma}\to \lambda \int_0^T\omega\int_{\Omega}\nabla \left(f({\phi})\mathbf{F}\right)\mathbin{\tensorm} \nabla \boldsymbol{\Sigma},
\]
as $m \to \infty$, $\delta \to 0$.

We now consider the third equation of \eqref{eqn:limit2d1}. Since
\[\omega \xi_m =\omega P_m^L(q)\in C^0(0,T;H^1(\Omega))\hookrightarrow L^{2}\left(0,T;H^1(\Omega)\right),\] 
we get from \eqref{eqn:conv2d7} that
\[
\int_0^T\omega\int_{\Omega}\partial_t\phi_m \xi_m \to \int_0^T\omega<\partial_t\phi, q>,
\]
as $m\to \infty$, $\delta \to 0$. 
Thanks to \eqref{eqn:conv2d8} and the fact that {\rosso $\xi_m=P_m^{L}(q)$} converges strongly to $q$ in $H^1(\Omega)$, with similar considerations as in \eqref{eqn:limitcalc1} we have that $\omega \nabla \phi_m \xi_m \to \omega \nabla \phi q$ strongly in $L^2(0,T;L^p(\Omega))$, with arbitrary $p\geq 1$.
Hence, using \eqref{eqn:conv2d1}, by the product of weak--strong convergence we have that 
\[
\int_0^T\omega\int_{\Omega}\left(\mathbf{v}_m \cdot \nabla \phi_m\right) \xi_m\to \int_0^T\omega\int_{\Omega}\left(\mathbf{v} \cdot \nabla \phi \right) q,
\] 
as $m\to \infty$, $\delta \to 0$. Finally, considering \eqref{eqn:conv2d8} and Assumption \textbf{A1}, $b(\phi_m)\to b(\phi)$ a.e. in $\Omega_T$ and is uniformly bounded, hence by applying a generalized form of the Lebesgue convergence theorem, considering also that $P_m^{L}(q)$ converges strongly to $q$ in $H^1(\Omega)$, we may obtain that $b(\phi_m)\nabla \xi_m\to b(\phi)\nabla q$ strongly in $L^{2}\left(0,T;L^2(\Omega;\mathbb{R}^2)\right)$. Then, by \eqref{eqn:conv2d9} and the product of weak--strong convergence we have that 
\[
\int_0^T\omega\int_{\Omega}b\left(\phi_m\right)\nabla \mu_m \cdot \nabla \xi_m\to \int_0^T\omega\int_{\Omega}b\left(\phi\right)\nabla \mu \cdot \nabla q,
\]
as $m\to \infty$, $\delta \to 0$.

{\rosso As for \eqref{eqn:limit2d1}$_4$,} from the convergence properties \eqref{eqn:conv2d6} and \eqref{eqn:conv2d9}, we obtain that
\[
\int_0^T\omega\int_{\Omega}\mu_m \chi_m\to \int_0^T\omega\int_{\Omega}\mu r,
\]
\[
\int_0^T\omega\int_{\Omega}\nabla \phi_m\cdot \nabla \chi_m\to \int_0^T\omega\int_{\Omega}\nabla \phi\cdot \nabla r,
\]
as $m\to \infty$, $\delta \to 0$. Moreover, using \eqref{eqn:conv2d8} and the continuity of $\psi'$, we have that $\psi'(\phi_m)\to \psi'(\phi)$ a.e. in $\Omega_T$. Then, thanks to Assumption \textbf{A2}, in particular to the growth law $|\psi'(\phi_m)|\leq c_1(|\phi_m|^p+1)$, for arbitrary $p\in \mathbb{R}^+$, and again to \eqref{eqn:conv2d8}, applying a generalized form of the Lebesgue convergence theorem we have that
\[
\int_0^T\omega\int_{\Omega}\psi'\left(\phi_m\right)\chi_m\to \int_0^T\omega\int_{\Omega}\psi'\left(\phi\right)r,
\]
as $m\to \infty$, $\delta \to 0$. Using the properties that
$f'\left(\phi_m\right)\mathbf{F}_m\colon \mathbf{F}_m\to f'\left(\phi\right)\mathbf{F}\colon \mathbf{F}$ strongly in $L^2(0,T;L^{2-h}(\Omega))$ and
$f'\left(\phi_m\right)\to f'\left(\phi\right)$ strongly in $L^{2}(0,T;L^2(\Omega))$, which we previously derived, we conclude that 
\[
\int_0^T\omega\int_{\Omega}\frac{f'\left(\phi_m\right)}{2}\left(\mathbf{F}_m\colon \mathbf{F}_m-2\right)\chi_m\to \int_0^T\omega\int_{\Omega}\frac{f'\left(\phi \right)}{2}\left(\mathbf{F}\colon \mathbf{F}-2\right)r,
\]
as $m\to \infty$, $\delta \to 0$. Lastly, from \eqref{eqn:conv2d10}, we get that
\[
\sqrt{\delta}\int_0^T\omega\int_{\Omega}\sqrt{\delta}\partial_t\phi_m \chi_m \to 0,
\]
as $m\to \infty$, $\delta \to 0$.

We need finally to prove that the initial conditions hold. Due to \eqref{eqn:conv2d8}, in particular to the fact that $\phi_m\to \phi$ strongly in $C^0([0,T];L^2(\Omega))$,  and the fact that $\phi_m(0)=P_m^L(\phi_0)\to \phi_0$ strongly in $L^2(\Omega)$, we have that $\phi(0)=\phi_0$ a.e. in $\Omega$. Moreover, due to \eqref{eqn:conv2d3} and \eqref{eqn:conv2d4}, 
\[
\mathbf{F}\in C^0\left([0,T];\left(H^1(\Omega;\mathbb{R}^{2\times 2})\right)'\right) \cap L^{\infty}(0,T;L^2(\Omega;\mathbb{R}^{2\times 2}))\hookrightarrow C_w^0\left([0,T];L^2(\Omega;\mathbb{R}^{2\times 2})\right).
\]
Then, using the same reasoning as in \cite{benesova}, we obtain from the convexity of $w(\phi, \cdot)$ and from the energy estimate \eqref{eqn:apriori} that $\mathbf{F}(\cdot,t)\to \mathbf{F}_0$ strongly in $L^2(\Omega;\mathbb{R}^{2\times 2})$ as $t\to 0^+$.

All the aforementioned results prove Theorem \ref{thm:2d}.

\subsection{Limit problem in the $3$-d case}
As we will see in the following calculations, in the $3$-d case the viscous regularization term $\delta \frac{\partial \phi_m}{\partial t}$ will be needed to enhance the regularity of $\partial_t \phi_m$ and to recover a compactenss result for the sequence $\mathbf{F}_m$, without disposing of a direct control on $\nabla \mathbf{F}_m$ in some B{\"o}chner space. Hence, the limit problem of \eqref{eqn:galerkin} in $3$ space dimension will be identified in the case $\gamma=\theta=0$. We thus consider system \eqref{eqn:galerkin3} with $\delta>0$.

We now derive the necessary regularity results to study the limit point of \eqref{eqn:galerkin3}, as $m\to \infty$, starting from \eqref{eqn:apriori}. From the Poincar\'e inequality and from \eqref{eqn:apriori} (with $\gamma=0$), we have that
{\rosso
\begin{equation}
\label{eqn:vm3d}
\mathbf{v}_m\quad \text{is uniformly bounded in} \; L^2(0,T;H_{0,\diver}^1(\Omega;\mathbb{R}^3))\hookrightarrow L^2(0,T;L_{\diver}^p(\Omega;\mathbb{R}^3)),
\end{equation}
with $p\in [1,6]$. Moreover, from \eqref{eqn:apriori} and Assumption \textbf{A3}, we have that
\begin{equation}
\label{eqn:fFm3d}
f(\phi_m)\mathbf{F}_m \quad \text{is uniformly bounded in} \; L^{\infty}(0,T;L^2(\Omega;\mathbb{R}^{3\times 3}))\cap L^{2}(0,T;H^1(\Omega;\mathbb{R}^{3\times 3})),
\end{equation}
which implies
\begin{equation}
\label{eqn:Fm3d}
||\mathbf{F}_m||_{L^{\infty}(0,T;L^2(\Omega;\mathbb{R}^{3\times 3}))}\leq C.
\end{equation}
Taking $h=2$ in \eqref{eqn:gninterpolation3d},
\begin{align}
\label{eqn:Fm43d}
& f(\phi_m)\mathbf{F}_m \quad \text{is uniformly bounded in} \\
& \notag L^{\infty}(0,T;L^2(\Omega;\mathbb{R}^{3\times 3}))\cap L^{2}(0,T;H^1(\Omega;\mathbb{R}^{3\times 3}))\hookrightarrow L^{\frac{8}{3}}(0,T;L^4(\Omega;\mathbb{R}^{3\times 3})),\\
& \notag \text{which implies} \quad ||\mathbf{F}_m||_{L^{\frac{8}{3}}(0,T;L^4(\Omega;\mathbb{R}^{3\times 3}))}\leq C. 
\end{align}
}%
Taking $\chi=-\Delta \phi_m$ in the fourth equation of \eqref{eqn:galerkin3}, using the Cauchy--Schwarz and Young inequalities and Assumptions \textbf{A2} and \textbf{A3}, we get
\begin{align}
\label{eqn:deltaphi13d}
& ||\Delta \phi_m||^2+\int_{\Omega}\psi_+''(\phi_m)|\nabla \phi_m|^2\leq ||\nabla \mu_m||_{L^2(\Omega;\mathbb{R}^3)}||\nabla \phi_m||_{L^2(\Omega;\mathbb{R}^3)}-\int_{\Omega}\psi_-''(\phi_m)|\nabla \phi_m|^2\\
& \notag +\int_{\Omega}\frac{f'\left(\phi_m\right)}{2}\left(\mathbf{F}_m\colon \mathbf{F}_m-2\right)\Delta \phi_m+\delta ||\partial_t \phi_m||\,||\Delta \phi_m||\leq \frac{1}{2}||\nabla \mu_m||_{L^2(\Omega;\mathbb{R}^3)}^2+\frac{1}{2}||\nabla \phi_m||_{L^2(\Omega;\mathbb{R}^3)}^2\\
& \notag +C\int_{\Omega}\left(1+|\phi_m|^q\right)|\nabla \phi_m|^2+\frac{f_1}{2}||\Delta \phi_m||\, ||\mathbf{F}_m||_{L^4(\Omega;\mathbb{R}^{3\times 3})}^2+\frac{1}{6}||\Delta \phi_m||^2+C||\partial_t \phi_m||^2+C.
\end{align} 
We use \eqref{eqn:agmon3d}, elliptic regularity and the Young inequality, observing that $\frac{4}{q}>1$ when $q<4$, to write
\begin{align*}
&\int_{\Omega}|\phi_m|^q|\nabla \phi_m|^2\leq ||\phi_m||_{L^{\infty}(\Omega)}^q||\nabla \phi_m||_{L^2(\Omega;\mathbb{R}^3)}^2\leq C||\phi_m||_{H^1(\Omega;\mathbb{R}^3)}^{\frac{4+q}{2}}\left(||\phi_m||^{\frac{q}{2}}+||\Delta \phi_m||^{\frac{q}{2}}\right) \\
&\leq C||\phi_m||_{H^1(\Omega;\mathbb{R}^3)}^{2+q}+||\phi_m||_{H^1(\Omega;\mathbb{R}^3)}^{\frac{2(q+4)}{4-q}}+\frac{1}{6}||\Delta \phi_m||^2.
\end{align*}
Using this inequality in \eqref{eqn:deltaphi13d}, together with the Young inequality, we obtain that
\begin{align}
\label{eqn:deltaphi23d}
&||\Delta \phi_m||^2\leq \frac{1}{2}||\nabla \mu_m||_{L^2(\Omega;\mathbb{R}^3)}^2+C||\phi_m||_{H^1(\Omega;\mathbb{R}^3)}^2+C||\phi_m||_{H^1(\Omega;\mathbb{R}^3)}^{2+q}+||\phi_m||_{H^1(\Omega;\mathbb{R}^3)}^{\frac{2(q+4)}{4-q}}\\
& \notag +C||\mathbf{F}_m||_{L^4(\Omega;\mathbb{R}^{3\times 3})}^4+\frac{1}{2}||\Delta \phi_m||^2+C||\partial_t \phi_m||^2+C.
\end{align}
Taking the power $\frac{2}{3}$ of both sides of \eqref{eqn:deltaphi23d}, integrating in time over the interval $(0,T)$, using \eqref{eqn:apriori}, Assumption \textbf{A1} and \eqref{eqn:Fm43d}, we infer that
\begin{equation}
\label{eqn:deltaphi3d}
\int_0^T||\Delta \phi_m||^\frac{4}{3}\leq C+C||\mathbf{F}_m||_{L^{\frac{8}{3}}(0,T;L^4(\Omega;\mathbb{R}^{3\times 3}))}^{\frac{8}{3}}\leq C,
\end{equation}
from which we have, {\rosso considering also \eqref{eqn:phimh1} and \eqref{eqn:deltaphi3d}, that
\begin{equation}
\label{eqn:phim3d}
\phi_m\quad \text{is uniformly bounded in} \; L^{\infty}(0,T;H^1(\Omega))\cap L^{\frac{4}{3}}(0,T;H^2(\Omega)).
\end{equation}
}%
We observe that, taking the Gagliardo--Nirenberg inequality \eqref{eqn:gagliardoniremberg}, specified for $j=0$, $m=1$, $p=2+h$, $r=2$, $q=2$, to the power of $\frac{2}{3}$ and integrating over $(0,T)$, thus obtaining \eqref{eqn:gninterpolation3d} with the integrands to the power of  $\frac{2}{3}$, choosing $h=2$, we get
{\rosso
\begin{equation}
\label{eqn:nablaphim43d}
||\nabla \phi_m||_{L^{\frac{16}{9}}(0,T;L^4(\Omega;\mathbb{R}^3))}\leq C,
\end{equation}
and choosing $h=1$,
\begin{equation}
\label{eqn:nablaphim33d}
||\nabla \phi_m||_{L^{\frac{8}{3}}(0,T;L^3(\Omega;\mathbb{R}^3))}\leq C.
\end{equation}
}%
Note that, unlike the $2$-d case, in the $3$-d case, given \eqref{eqn:Fm43d} and \eqref{eqn:nablaphim43d}, we cannot use the arguments in \eqref{eqn:gradFmcalc} to deduce \eqref{eqn:gradFm}, which is not valid in $3$ space dimentions.

We now deduce a priori estimates, uniform in $m$, for the time derivative of $\mathbf{F}_m$ in \eqref{eqn:galerkin3}. We note that, from the energy estimate \eqref{eqn:apriori}, with $\delta>0$, we already have
{\rosso
\begin{equation}
\label{eqn:phimdual3d}
||\partial_t \phi_m||_{L^{2}\left(0,T;L^2(\Omega)\right)}\leq C.
\end{equation}
}%
Multiplying \eqref{eqn:fmparts} by a time function $\zeta \in L^{4}(0,T)$ and integrating in time over the interval $(0,T)$, considering \eqref{eqn:gninterpolation3d} with $h=1$ and \eqref{eqn:vm3d}, we get
\begin{align}
\label{eqn:mfdualcalc13d}
&\int_0^T\int_{\Omega}\left({\mathbf{v}}_m \cdot \nabla \right)\mathbf{F}_m\colon \boldsymbol{\Sigma}\,\zeta=-\int_0^T\int_{\Omega}\mathbf{F}_m\otimes \mathbf{v}_m\mathbin{\tensorm} \nabla \boldsymbol{\Sigma}\,\zeta \\
& \notag \leq  \int_0^T ||\mathbf{F}_m||_{L^{3}(\Omega;\mathbb{R}^{3\times 3})}||\mathbf{v}_m||_{L^{6}(\Omega;\mathbb{R}^{3})}||\nabla \boldsymbol{\Sigma}||_{L^{2}(\Omega;\mathbb{R}^{3\times 3 \times 3})}|\zeta| \\
& \notag \leq C||\mathbf{F}_m||_{L^{4}\left(0,T;L^{3}(\Omega;\mathbb{R}^{3\times 3})\right)}||\mathbf{v}_m||_{L^2\left(0,T;L^{6}(\Omega;\mathbb{R}^{3})\right)}||\nabla \boldsymbol{\Sigma}||_{L^{2}(\Omega;\mathbb{R}^{3\times 3 \times 3})}||\zeta||_{L^{4}(0,T)}.
\end{align}
Then, we have that
\begin{equation}
\label{eqn:fmdual13d}
\left|\left|\left({\mathbf{v}}_m \cdot \nabla \right)\mathbf{F}_m\right|\right|_{L^{\frac{4}{3}}\left(0,T;\left(H^1(\Omega;\mathbb{R}^{3\times 3})\right)'\right)}\leq C.
\end{equation}
Using \eqref{eqn:fmdual13d}, multiplying \eqref{eqn:galerkin3}$_2$ by a time function $\zeta \in L^{4}(0,T)$ and integrating in time over the interval $(0,T)$, using moreover \eqref{eqn:apriori}, we get
\begin{align}
\label{eqn:mfdualcalc23d}
&\int_0^T\left<\partial_t \mathbf{F}_m, \boldsymbol{\Sigma} \right>\zeta \leq \int_0^T ||\left({\mathbf{v}}_m \cdot \nabla \right)\mathbf{F}_m||_{\left(H^1(\Omega;\mathbb{R}^{3\times 3})\right)'}||\boldsymbol{\Sigma}||_{H^{1}(\Omega;\mathbb{R}^{3\times 3})}|\zeta|\\
& \notag +\int_0^T ||\nabla \mathbf{v}_m||_{L^{2}(\Omega;\mathbb{R}^{3})}||\mathbf{F}_m||_{L^{3}(\Omega;\mathbb{R}^{3\times 3})}||\boldsymbol{\Sigma}||_{L^{6}(\Omega;\mathbb{R}^{3\times 3})}|\zeta|\\
& \notag +\lambda\int_0^T||\nabla \left(f({\phi}_m)\mathbf{F}_m\right)||_{L^{2}(\Omega;\mathbb{R}^{3\times 3})}||\nabla \boldsymbol{\Sigma}||_{L^{2}(\Omega;\mathbb{R}^{3\times 3\times 3})}|\zeta|\leq C || \boldsymbol{\Sigma}||_{H^{1}(\Omega;\mathbb{R}^{3\times 3})}||\zeta||_{L^{4}(0,T)}\\
& \notag +C||\nabla \mathbf{v}_m||_{L^2\left(0,T;L^{2}(\Omega;\mathbb{R}^{3\times 3})\right)}||\mathbf{F}_m||_{L^{4}\left(0,T;L^{3}(\Omega;\mathbb{R}^{3\times 3})\right)}||\boldsymbol{\Sigma}||_{H^{1}(\Omega;\mathbb{R}^{3\times 3})}||\zeta||_{L^{4}(0,T)}.
\end{align}
Hence, we deduce that
{\rosso
\begin{equation}
\label{eqn:fmdual23d}
||\partial_t \mathbf{F}_m||_{L^{\frac{4}{3}}\left(0,T;\left(H^1(\Omega;\mathbb{R}^{3\times 3})\right)'\right)}\leq C.
\end{equation}
Given \eqref{eqn:fmdual23d} and \eqref{eqn:phimdual3d}, we can infer that
\begin{equation}
\label{eqn:phimfmdual}
||\partial_t \left(f(\phi_m) \mathbf{F}_m\right)||_{L^{\frac{4}{3}}\left(0,T;\left(H_N^2(\Omega;\mathbb{R}^{3\times 3})\right)'\right)}\leq C.
\end{equation}
}%
The latter result is fundamental to study the limit problem of \eqref{eqn:galerkin3}, in particular to deduce compactness results for the sequence $\mathbf{F}_m$ from a compactness result for the sequence $f(\phi_m) \mathbf{F}_m$, which can be obtained from \eqref{eqn:fFm3d} and \eqref{eqn:phimfmdual}. In order to prove \eqref{eqn:phimfmdual}, let us take generic functions $\boldsymbol{\Pi} \in H_N^2(\Omega;\mathbb{R}^{3\times 3})$ and $\zeta \in L^4(0,T)$, and consider the integral in time over $(0,T)$ of the dual product $<\partial_t \left(\phi_m \mathbf{F}_m\right),\boldsymbol{\Pi}>_2$ multiplied by $\zeta$. We get, using Assumption \textbf{A3},
\begin{align*}
& \int_0^T\zeta<\partial_t \left(f(\phi_m) \mathbf{F}_m\right), \boldsymbol{\Pi}>_2=\int_0^T\int_{\Omega}f'(\phi_m)(\partial_t \phi_m)\mathbf{F}_m\colon \boldsymbol{\Pi}\, \zeta+\int_0^T\zeta<\partial_t \mathbf{F}_m, f(\phi_m)\boldsymbol{\Pi}> \\
& \leq C||\partial_t \phi_m||_{L^2(\Omega_T)}||\mathbf{F}_m||_{L^{4}\left(0,T;L^{3}(\Omega;\mathbb{R}^{3\times 3})\right)}||\zeta||_{L^{4}(0,T)}||\boldsymbol{\Pi}||_{H^{1}(\Omega;\mathbb{R}^{3\times 3})}\\
& +C||\partial_t \mathbf{F}_m||_{L^{\frac{4}{3}}\left(0,T;\left(H^1(\Omega;\mathbb{R}^{3\times 3})\right)'\right)}||f(\phi_m)\boldsymbol{\Pi}||_{L^{\infty}(0,T;H^1(\Omega;\mathbb{R}^{3\times 3}))}||\zeta||_{L^{4}(0,T)},
\end{align*}
and from the bound
\begin{align*}
&||f(\phi_m)\boldsymbol{\Pi}||_{L^{\infty}(0,T;H^1(\Omega;\mathbb{R}^{3\times 3}))}\\
&=\underset{(0,T)}{\text{ess.sup}}\left(\int_{\Omega}\left(|f(\phi_m)\boldsymbol{\Pi}|^2+|\boldsymbol{\Pi}\otimes \nabla f(\phi_m)|^2+|f(\phi_m)\nabla \boldsymbol{\Pi}|^2\right)+2f(\phi_m)\nabla \boldsymbol{\Pi}\mathbin{\tensorm} \boldsymbol{\Pi}\otimes \nabla f(\phi_m)\right)^{\frac{1}{2}} \\
& \leq C\left(1+||\phi_m||_{L^{\infty}(0,T;H^1(\Omega))}\right)||\boldsymbol{\Pi}||_{H_N^2(\Omega;\mathbb{R}^{3\times 3})},
\end{align*}
where we used the embedding $H_N^2(\Omega;\mathbb{R}^{3\times 3})\hookrightarrow L^{\infty}(\Omega;\mathbb{R}^{3\times 3})$, we obtain \eqref{eqn:phimfmdual}.

Collecting the estimates \eqref{eqn:vm3d}, \eqref{eqn:fFm3d}, \eqref{eqn:Fm3d}, \eqref{eqn:phim3d}, \eqref{eqn:mumh1},  \eqref{eqn:phimfmdual} and \eqref{eqn:phimdual3d}, which are uniform in $m$, from the Banach--Alaoglu and the Aubin--Lions lemma, we finally obtain the convergence properties, up to subsequences of the solutions, which we still label by the index $m$,
\begin{align}
\label{eqn:conv3d1} & \mathbf{v}_m \rightharpoonup \mathbf{v} \quad \text{in} \quad L^{2}\left(0,T;H_{0,\diver}^1\left(\Omega;\mathbb{R}^{3}\right)\right),\\
\label{eqn:conv3d2} & f(\phi_m)\mathbf{F}_m \overset{\ast}{\rightharpoonup} \boldsymbol{\beta} \quad \text{in} \quad L^{\infty}(0,T;L^2(\Omega;\mathbb{R}^{3\times 3}))\cap L^{2}(0,T;H^1(\Omega;\mathbb{R}^{3\times 3})),\\
\label{eqn:conv3d3} & \mathbf{F}_m \overset{\ast}{\rightharpoonup} \mathbf{F} \quad \text{in} \quad L^{\infty}(0,T;L^2(\Omega;\mathbb{R}^{3\times 3})),\\
\label{eqn:conv3d4} & \partial_t \mathbf{F}_m \rightharpoonup \partial_t \mathbf{F} \quad \text{in} \quad L^{\frac{4}{3}}\left(0,T;\left(H^1(\Omega;\mathbb{R}^{3\times 3})\right)'\right),\\
\label{eqn:conv3d4bis} & \partial_t \left(f(\phi_m)\mathbf{F}_m\right) \rightharpoonup \partial_t \boldsymbol{\beta} \quad \text{in} \quad L^{\frac{4}{3}}\left(0,T;\left(H_N^2(\Omega;\mathbb{R}^{3\times 3})\right)'\right),\\
\label{eqn:conv3d5} & f(\phi_m)\mathbf{F}_m \to \boldsymbol{\beta} \quad \text{in} \quad L^{2}(0,T;L^p(\Omega;\mathbb{R}^{3\times 3})), \quad p\in [1,6), \quad \text{and} \quad \text{a.e. in} \; \; \Omega_T,\\
\label{eqn:conv3d6} & \phi_m \overset{\ast}{\rightharpoonup} \phi \quad \text{in} \quad L^{\infty}(0,T;H^1(\Omega))\cap L^{\frac{4}{3}}(0,T;H^2(\Omega)),\\
\label{eqn:conv3d7} & \partial_t \phi_m \rightharpoonup \partial_t \phi \quad \text{in} \quad L^{2}(\Omega_T),\\
\label{eqn:conv3d8} & \phi_m \to \phi \quad \text{in} \quad C^{0}(0,T;L^p(\Omega)) \cap L^{\frac{4}{3}}(0,T;W^{1,p}(\Omega)), \quad p\in [1,6), \quad \text{and} \quad \text{a.e. in} \; \; \Omega_T,\\
\label{eqn:conv3d9} & \mu_m \rightharpoonup \mu \quad \text{in} \quad L^{2}\left(0,T;H^1\left(\Omega\right)\right),
\end{align}
for some limit function $\boldsymbol{\beta}$ in $L^{\infty}(0,T;L^2(\Omega;\mathbb{R}^{3\times 3}))\cap L^{2}(0,T;H^1(\Omega;\mathbb{R}^{3\times 3}))$, which will be identified as $\boldsymbol{\beta}=f(\phi)\mathbf{F}$ in the following calculations. Moreover, thanks to \eqref{eqn:conv3d8}, the Assumption \textbf{A3} and {\rosso the Lebesgue dominated convergence theorem},
\begin{equation}
\label{eqn:conv3d10}
f(\phi_m)\to f(\phi), \; \frac{1}{f(\phi_m)}\to \frac{1}{f(\phi)}, \; f'(\phi_m)\to f'(\phi) \quad \text{in} \quad L^{p}\left(\Omega_T\right), \, \forall p\geq 1.
\end{equation}

We now derive a useful interpolation result for $\phi_m$ from \eqref{eqn:conv3d6} and \eqref{eqn:conv3d8}, which will be used to identify the limit of system \eqref{eqn:galerkin3} as $m\to \infty$. 
From \eqref{eqn:lpinterpolation} with $s=2, q=6$ we have that
\[
L^{\infty}(0,T;L^2(\Omega,\mathbb{R}^3))\cap L^{\frac{4}{3}}(0,T;L^6(\Omega,\mathbb{R}^3))\hookrightarrow L^{\frac{26}{9}}(0,T;L^{\frac{26}{9}}(\Omega,\mathbb{R}^3)),
\]
hence, from \eqref{eqn:conv3d6} and \eqref{eqn:conv3d8} we have that
\begin{equation}
\label{eqn:phiminterp1}
\nabla \phi_m \to \nabla \phi \quad \text{in} \quad  L^{p}(0,T;L^{p}(\Omega,\mathbb{R}^3)), \; p\in [1,26/9).
\end{equation}
We also observe that, from \eqref{eqn:conv3d2} and \eqref{eqn:lpinterpolation} with $q=6$ we have that
\begin{equation}
\label{eqn:fFminterp1}
f(\phi_m)\mathbf{F}_m \to \boldsymbol{\beta} \quad \text{in} \quad  L^{p}(0,T;L^{p}(\Omega,\mathbb{R}^{3\times 3})), \; p\in [1,10/3),
\end{equation}
and from \eqref{eqn:gninterpolation3d} with $h\leftarrow 4-h$,
\begin{equation}
\label{eqn:fFminterp2}
\boldsymbol{\beta} \in L^{2+\frac{2}{3}\frac{h}{4-h}}(0,T;L^{6-h}(\Omega,\mathbb{R}^{3\times 3})).
\end{equation}
We prove now a compactness result for the sequence $\mathbf{F}_m$. Given \eqref{eqn:conv3d5}, \eqref{eqn:fFminterp2}, \eqref{eqn:conv3d10} and writing $\mathbf{F}_m$ as
\[
\mathbf{F}_m=\frac{1}{f(\phi_m)}f(\phi_m)\mathbf{F}_m,
\]
we have, for arbitrary small $h\in (0,5/2]$, using also Assumption \textbf{A3},
\begin{align}
\label{eqn:fmconv1}
& \int_0^T\left(\int_{\Omega}\left|\mathbf{F}_m-\frac{\boldsymbol{\beta}}{f(\phi)}\right|^{6-2h}\right)^{\frac{2}{6-2h}}\leq \int_0^T\left(\int_{\Omega}\left|\frac{1}{f(\phi_m)}\left(f(\phi_m)\mathbf{F}_m-\boldsymbol{\beta}\right)\right|^{6-2h}\right)^{\frac{2}{6-2h}}\\
& \notag +\int_0^T\left(\int_{\Omega}\left|\left(\frac{1}{f(\phi_m)}-\frac{1}{f(\phi)}\right)\boldsymbol{\beta}\right|^{6-2h}\right)^{\frac{2}{6-2h}}\leq C||f(\phi_m)\mathbf{F}_m-\boldsymbol{\beta}||_{L^{2}(0,T;L^{6-2h}(\Omega;\mathbb{R}^{3\times 3}))}^{2}\\
& \notag +C\left|\left|\frac{1}{f(\phi_m)}-\frac{1}{f(\phi)}\right|\right|_{L^{\frac{24-4h}{h}}\left(0,T;L^{\frac{(6-2h)(6-h)}{h}}(\Omega)\right)}^{2}||\boldsymbol{\beta}||_{L^{2+\frac{2}{3}\frac{h}{4-h}}(0,T;L^{6-h}(\Omega;\mathbb{R}^{3\times 3}))}^{2}\to 0,
\end{align}
as $m\to \infty$. Hence,
\begin{equation}
\label{eqn:fmstrong1}
\mathbf{F}_m\to \frac{\boldsymbol{\beta}}{f(\phi)} \quad \text{strongly in} \quad L^{2}(0,T;L^{6-h}(\Omega;\mathbb{R}^{3\times 3})) \; \text{and a.e. in} \; \Omega_T,
\end{equation}
with $h\in (0,5]$
From \eqref{eqn:fmstrong1}, \eqref{eqn:conv3d3} and and the uniqueness of the weak limit and the a.e. convergence, $\mathbf{F}=\frac{\boldsymbol{\beta}}{f(\phi)}$, whence
\begin{equation}
\label{eqn:beta3d}
\boldsymbol{\beta}=f(\phi)\mathbf{F}.
\end{equation}
Analogously to \eqref{eqn:fFminterp1}, we thus have that
\begin{equation}
\label{eqn:fminterp1}
\mathbf{F}_m \to \mathbf{F} \quad \text{in} \quad  L^{p}(0,T;L^{p}(\Omega,\mathbb{R}^{3\times 3})), \; p\in [1,10/3).
\end{equation}
With the convergence results \eqref{eqn:conv3d1}--\eqref{eqn:conv3d10}, \eqref{eqn:phiminterp1}, \eqref{eqn:beta3d} and \eqref{eqn:fmstrong1} we can pass to the limit in the system \eqref{eqn:galerkin3} as $m\to \infty$. We proceed in a different way with respect to the $2-d$ case.
We fix $0\leq i \leq m$, and take $\boldsymbol{\eta}=\boldsymbol{\eta}_i$, $\xi=\xi_i$, $\chi=\xi_i$, multiply the equations by a function $\omega \in C_0^{\infty}([0,T])$ and integrate over the time interval $[0,T]$. This gives
\begin{equation}
\label{eqn:limit3d1}
\begin{cases}
\displaystyle \nu \int_0^T\omega \int_{\Omega}\nabla \mathbf{v}_m\colon \nabla \boldsymbol{\eta}_i= \int_0^T\omega\int_{\Omega}\mu_m \nabla \phi_m \cdot \boldsymbol{\eta}_i-\int_0^T\omega\int_{\Omega}\frac{f'\left(\phi_m\right)}{2}\left(\mathbf{F}_m\colon \mathbf{F}_m-3\right)\nabla \phi_m\cdot \boldsymbol{\eta}_i\\ 
\displaystyle -\int_0^T\omega\int_{\Omega}f\left(\phi_m\right)\mathbf{F}_m\mathbf{F}_m^T\colon \nabla \boldsymbol{\eta}_i,\\ \\
\displaystyle  \int_0^T\omega\left<\partial_t \mathbf{F}_m, \boldsymbol{\Sigma} \right>+\int_0^T\omega\int_{\Omega}\left({\mathbf{v}}_m \cdot \nabla \right)\mathbf{F}_m\colon \boldsymbol{\Sigma}-\int_0^T\omega\int_{\Omega}\left(\nabla {\mathbf{v}}_m\right)\mathbf{F}_m\colon \boldsymbol{\Sigma}\\
\displaystyle +\lambda \int_0^T\omega\int_{\Omega}\nabla \left(f({\phi}_m)\mathbf{F}_m\right)\mathbin{\tensorm} \nabla \boldsymbol{\Sigma}=0,\\ \\
 \displaystyle \int_0^T\omega\int_{\Omega}\partial_t\phi_m \xi_i+\int_0^T\omega\int_{\Omega}\left(\mathbf{v}_m \cdot \nabla \phi_m\right) \xi_i+\int_0^T\omega\int_{\Omega}b\left(\phi_m\right)\nabla \mu_m \cdot \nabla \xi_i=0,\\ \\
\displaystyle \int_0^T\omega\int_{\Omega}\mu_m \xi_i=\int_0^T\omega\int_{\Omega}\nabla \phi_m\cdot \nabla \xi_i+\int_0^T\omega\int_{\Omega}\psi'\left(\phi_m\right)\xi_i\\
\displaystyle +\int_0^T\omega\int_{\Omega}\frac{f'\left(\phi_m\right)}{2}\left(\mathbf{F}_m\colon \mathbf{F}_m-3\right)\xi_i+\delta \int_0^T\omega\int_{\Omega}\partial_t\phi_m \xi_i.
\end{cases}
\end{equation} 
Thanks to \eqref{eqn:conv3d1}, we have
\[
\nu \int_0^T\omega \int_{\Omega}\nabla \mathbf{v}_m\colon \nabla \boldsymbol{\eta}_i\to \nu \int_0^T\omega \int_{\Omega}\nabla \mathbf{v}\colon \nabla \boldsymbol{\eta}_i,
\]
as $m\to \infty$. Thanks to \eqref{eqn:phiminterp1} and the fact that $\boldsymbol{\eta}_i\in H^2(\Omega;\mathbb{R}^3)\hookrightarrow L^{\infty}(\Omega;\mathbb{R}^3)$, we have that $\omega \nabla \phi_m \cdot \boldsymbol{\eta}_i\to \omega \nabla \phi \cdot \boldsymbol{\eta}_i$ strongly in $L^2(0,T;L^2(\Omega))$.
Hence, using \eqref{eqn:conv3d9}, by the product of weak--strong convergence we have that 
\[
\int_0^T\omega\int_{\Omega}\mu_m \nabla \phi_m \cdot \boldsymbol{\eta}_i\to \int_0^T\omega\int_{\Omega}\mu \nabla \phi \cdot \boldsymbol{\eta}_i,
\] 
as $m\to \infty$. For what concerns the second term on the right hand side of the first equation of \eqref{eqn:limit3d1}, we observe that $f'\left(\phi_m\right)\mathbf{F}_m\colon \mathbf{F}_m\to f'\left(\phi\right)\mathbf{F}\colon \mathbf{F}$ strongly in $L^{\frac{4}{3}-h}(0,T;L^{2}(\Omega))$, with $h\in (0,1/3]$. In order to prove this convergence result, we recall the following regularity result,
\begin{equation}
\label{eqn:interp45}
\mathbf{F}_m, \mathbf{F} \in L^{\frac{8}{3}-\frac{4}{3}\frac{h}{2+h}}(0,T;L^{4+h}(\Omega,\mathbb{R}^{3\times 3}))\cap L^{\frac{20}{9}-\frac{2}{3}h}\left(0,T;L^{5+\frac{9h}{4-3h}}(\Omega,\mathbb{R}^{3\times 3})\right),
\end{equation}
with $h\in [0,1/3]$, which are obtained from \eqref{eqn:fFm3d} and \eqref{eqn:gninterpolation3d} with $h\leftarrow 2+h$ and $h\leftarrow 3+\frac{9h}{4-3h}$. Introducing the parameter $k:=\frac{4}{3}\frac{h}{2+h}$, with $h\in (0,2/3]\to k\in (0,1/3]$, using Assumption \textbf{A3}, \eqref{eqn:ffid1}, \eqref{eqn:conv3d10} and \eqref{eqn:interp45}, we have that
\begin{align*}
& \int_0^T\left(\int_{\Omega}\left|f'\left(\phi_m\right)\mathbf{F}_m\colon \mathbf{F}_m - f'\left(\phi\right)\mathbf{F}\colon \mathbf{F}\right|^{2}\right)^{\frac{4-3k}{6}}\\
& =\int_0^T\left(\int_{\Omega}\left|f'\left(\phi_m\right)\left((\mathbf{F}_m-\mathbf{F})\colon (\mathbf{F}_m+\mathbf{F})\right)+\left(f'\left(\phi_m\right)-f'\left(\phi \right)\right)\mathbf{F}\colon \mathbf{F}\right|^{2}\right)^{\frac{4-3k}{6}} \\
& \leq f_1\int_0^T||\mathbf{F}_m-\mathbf{F}||_{L^{\frac{10}{3}-k}(\Omega;\mathbb{R}^{3\times 3})}^{\frac{4}{3}-k}||\mathbf{F}_m+\mathbf{F}||_{L^{5+\frac{9k}{4-3k}}(\Omega;\mathbb{R}^{3\times 3})}^{\frac{4}{3}-k}\\
& +\int_0^T||f'(\phi_m)-f'(\phi)||_{L^{\frac{2(4+h)}{h}}(\Omega)}^{\frac{4}{3}-k}||\mathbf{F}||_{L^{4+h}(\Omega;\mathbb{R}^{3\times 3})}^{\frac{8}{3}-2k} \\
&  \leq f_1||\mathbf{F}_m-\mathbf{F}||_{L^{\frac{10}{3}-k}\left(0,T;L^{\frac{10}{3}-k}(\Omega;\mathbb{R}^{3\times 3})\right)}^{\frac{4}{3}-k}||\mathbf{F}_m+\mathbf{F}||_{L^{\frac{20}{9}-\frac{k(14-3k)}{6}}\left(0,T;L^{5+\frac{9k}{4-3k}}(\Omega;\mathbb{R}^{3\times 3})\right)}^{\frac{4}{3}-k}\\
& +||f'(\phi_m)-f'(\phi)||_{L^{\frac{(4-3k)(8-3k)}{9k}}\left(0,T;L^{\frac{2(4+h)}{h}}(\Omega)\right)}^{\frac{4}{3}-k}||\mathbf{F}||_{L^{\frac{8}{3}-k}\left(0,T;L^{4+h}(\Omega;\mathbb{R}^{3\times 3})\right)}^{\frac{8}{3}-2k}\to 0,
\end{align*}
as $m\to \infty$, where we used the fact that $\frac{2}{3}k<\frac{k(14-3k)}{6}$ for $k\in (0,1/3]$. 
Then, considering again that $\boldsymbol{\eta}_i\in L^{\infty}(\Omega;\mathbb{R}^3)$ we have that $f'\left(\phi_m\right)\left(\mathbf{F}_m\colon \mathbf{F}_m-3\right)\boldsymbol{\eta}_i\to f'\left(\phi\right)\left(\mathbf{F}\colon \mathbf{F}-3\right)\boldsymbol{\eta}_i$ strongly in $L^{\frac{4}{3}-h}(0,T;L^{2}(\Omega;\mathbb{R}^3))$, where we also used Assumption \textbf{A3}, \eqref{eqn:conv3d8} and the Lebesgue convergence theorem to obtain $f'\left(\phi_m\right)\to f'\left(\phi\right)$ strongly in $L^{\frac{4}{3}-h}(0,T;L^2(\Omega))$. Hence, using \eqref{eqn:conv3d6} and the product of weak--strong convergence we have that 
\[
\int_0^T\omega\int_{\Omega}\frac{f'\left(\phi_m\right)}{2}\left(\mathbf{F}_m\colon \mathbf{F}_m-3\right)\nabla \phi_m\cdot \boldsymbol{\eta}_i\to \int_0^T\omega\int_{\Omega}\frac{f'\left(\phi\right)}{2}\left(\mathbf{F}\colon \mathbf{F}-3\right)\nabla \phi\cdot \boldsymbol{\eta}_i,
\] 
as $m\to \infty$. For what concerns the third term on the right hand side of the first equation of \eqref{eqn:limit3d1}, since, from \eqref{eqn:conv3d5}, $f(\phi_m)\mathbf{F}_m\to f(\phi)\mathbf{F}$ strongly in $L^2(0,T;L^{6-h}(\Omega;\mathbb{R}^{3\times 3}))$, $h\in (0,5]$, and considering that $\boldsymbol{\eta}_i\in H^2(\Omega;\mathbb{R}^3)\hookrightarrow W^{1,6}(\Omega;\mathbb{R}^3)$, we have that
$f(\phi_m)\nabla \boldsymbol{\eta}_i\mathbf{F}_m\to f(\phi)\nabla \boldsymbol{\eta}_i\mathbf{F}$ strongly in $L^2(0,T;L^{3-h}(\Omega;\mathbb{R}^{3\times 3}))$. Hence, using \eqref{eqn:conv3d3} and the product of weak--strong convergence we have that 
\[
\int_0^T\omega\int_{\Omega}f\left(\phi_m\right)\mathbf{F}_m\mathbf{F}_m^T\colon \nabla \boldsymbol{\eta}_i \to \int_0^T\omega\int_{\Omega}f\left(\phi\right)\mathbf{F}\mathbf{F}^T\colon \nabla \boldsymbol{\eta}_i,
\]
as $m\to \infty$.

Now considering the second equation in \eqref{eqn:limit3d1}, since
\[\omega \boldsymbol{\Sigma}\in C^0(0,T;H^1(\Omega;\mathbb{R}^{3\times 3}))\hookrightarrow L^{4}\left(0,T;H^1(\Omega;\mathbb{R}^{3\times 3})\right),\] 
we get from \eqref{eqn:conv3d4} that
\[
\int_0^T\omega\left<\partial_t \mathbf{F}_m, \boldsymbol{\Sigma} \right>\to \int_0^T\omega\left<\partial_t \mathbf{F}, \boldsymbol{\Sigma} \right>,
\]
as $m\to \infty$. Thanks to \eqref{eqn:fmstrong1} and the fact that $\nabla \boldsymbol{\Sigma}\in L^2(\Omega;\mathbb{R}^{3\times 2\times 3})$, we have that $\omega \left(\nabla \boldsymbol{\Sigma}\right)^T\mathbf{F}_m \to \omega \left(\nabla \boldsymbol{\Sigma}\right)^T\mathbf{F}$ strongly in $L^2(0,T;L^{\frac{3}{2}-h}(\Omega;\mathbb{R}^{3}))$. Hence, using \eqref{eqn:conv3d1} and the product of weak--strong convergence we have that 
\[
\int_0^T\omega\int_{\Omega}\left({\mathbf{v}}_m \cdot \nabla \right)\mathbf{F}_m\colon \boldsymbol{\Sigma}=-\int_0^T\omega \int_{\Omega}\mathbf{F}_m\otimes \mathbf{v}_m\mathbin{\tensorm} \nabla \boldsymbol{\Sigma}\to -\int_0^T\omega \int_{\Omega}\mathbf{F}\otimes \mathbf{v}\mathbin{\tensorm} \nabla \boldsymbol{\Sigma},
\]
as $m \to \infty$.

Moreover, thanks to \eqref{eqn:fmstrong1}, we have that $\omega\boldsymbol{\Sigma}\mathbf{F}_m^T \to \omega\boldsymbol{\Sigma}\mathbf{F}^T$ strongly in $L^2(0,T;L^{3-h}(\Omega;\mathbb{R}^{3\times 3}))$. Hence, using \eqref{eqn:conv3d1} and the product of weak--strong convergence we have that 
\[
\int_0^T\omega\int_{\Omega}\left(\nabla {\mathbf{v}}_m\right)\mathbf{F}_m\colon \boldsymbol{\Sigma}\to \int_0^T\omega\int_{\Omega}\left(\nabla {\mathbf{v}}\right)\mathbf{F}\colon \boldsymbol{\Sigma},
\]
as $m \to \infty$. Finally, given \eqref{eqn:conv3d2} and \eqref{eqn:beta3d},
\[
\lambda \int_0^T\omega\int_{\Omega}\nabla \left(f({\phi}_m)\mathbf{F}_m\right)\mathbin{\tensorm} \nabla \boldsymbol{\Sigma}\to \lambda \int_0^T\omega\int_{\Omega}\nabla \left(f({\phi})\mathbf{F}\right)\mathbin{\tensorm} \nabla \boldsymbol{\Sigma},
\]
as $m \to \infty$.

The limit forms of the terms in the third and fourth equations of \eqref{eqn:limit2d1} can be identified with similar calculations as in the $2$-d case, given the convergence properties obtained in the $3$-d case. Hence, we omit the details. Also, the same arguments as in the $2-d$ case can be used to verify the attainment of the initial conditions.

We have just shown that the limit point $(\mathbf{v},\mathbf{F},\phi,\mu)$ satisfies the same weak formulation \eqref{eqn:continuous2d}, with the addition of the viscous term in the equation for $\mu$, with $\mathbf{u}=\boldsymbol{\eta}_i$, $q=\xi_i$, $r=\xi_i$, for all $i\geq 0$. Observing that $\{\boldsymbol{\eta}_i\}_{i \in \mathbb{N}}$ is an Hilbert basis of $H^2(\Omega;\mathbb{R}^3)\cap H_{0,\diver}^1(\Omega;\mathbb{R}^3)$, with $H^2(\Omega;\mathbb{R}^3)\cap H_{0,\diver}^1(\Omega;\mathbb{R}^3)$ dense in $H_{0,\diver}^1(\Omega;\mathbb{R}^3)$, and also that $\{\xi_i\}_{i \in \mathbb{N}}$ is an Hilbert basis of {\rosso $H_N^2(\Omega)$, with $H_N^2(\Omega)$} dense in $H^1(\Omega)$, we have proved Theorem \ref{thm:3d}.

\section{Conclusion}
In this paper we have proposed a class of phase field models of Cahn--Hilliard type coupled with finite viscoelasticity, which considers the microscopic coupling between the mixture constituents, their viscoelastic behaviour and the friction between the phases, with a general compressible polyconvex elastic energy which depends on both the elastic deformation and the phase field variable within a nonlinear coupling. The model was formulated in the Eulerian configuration and it was derived by imposing the mass balance for the mixture components and the momentum balance implied by a generalized form of the principle of virtual powers, given constitutive assumptions complying with a mechanical version of the second law of thermodynamics in isothermal situations. We then studied analytically a simplified and regularized incompressible version of the model, where the regularization substantially depends on the coupling between the phase field variable and the elastic deformation gradient in the elastic energy density. Through a Galerkin approximation of the model equations and the definition of a fixed point problem, where the velocity field is approximated in terms of basis functions of the Stokes operator, the phase field variables are approximated in terms of basis functions of the Laplace operator, and the transport equation for the elastic deformation gradient is solved in a continuous space, we proved existence of a global in time weak solution in two and three space dimensions by studying the limit problem of the Galerkin approximation. In the case of three space dimensions, we added a viscous regularization for the partial time derivative of the phase field variable in the chemical potential in order to obtain the necessary compactness results to study the limit problem of the Galerkin approximation.
Future developments of the present work will investigate the inclusion in the elasticity dynamics of plasticity and growth, {\rosso together with the design of stable and convergent numerical approximations and the investigation of numerical simulations of related image-informed models for neurological diseases and tumor growth models. Also, future developments will concern the analysis for the polyconvex case.}

\section{Acknowledgements}
This research was supported by the Italian Ministry of Education,
University and Research (MIUR): Dipartimenti di Eccellenza Program (2018--2022)
-- Dept.~of Mathematics ``F.~Casorati'', University of Pavia.
In addition, AA, PC and ER  gratefully mention some other support
from the MIUR-PRIN Grant 2020F3NCPX ``Mathematics for industry 4.0 (Math4I4)'' and
their affiliation to the GNAMPA (Gruppo Nazionale per l'Analisi Matematica,
la Probabilit\`a e le loro Applicazioni) of INdAM (Isti\-tuto
Nazionale di Alta Matematica). 

%
\bibliographystyle{plain}
\bibliography{biblio_CHVE} 

\begin{thebibliography}{10}

\bibitem{agmon}
S.~Agmon.
\newblock {\em Lectures on elliptic boundary value problems}.
\newblock Providence, RI: AMS Chelsea Publishing, 2010.

\bibitem{allen}
M.~B. Allen, I.~Herrera, and G.~F. Pinder.
\newblock {\em Numerical modeling in science and engineering}.
\newblock John Wiley \& Sons, 1988.

\bibitem{alt}
H.~W. Alt.
\newblock {\em {L}inear {F}unctional {A}nalysis. An application-oriented
  introduction. Translated from the German edition by Robert N\"urnberg.}
\newblock Universitext. Springer-Verlag London, Ltd., London, 2016, 2016.
\newblock ISBN 978-1-4471-7279-6.

\bibitem{balbi}
V.~Balbi, A.~Trotta, M.~Destrade, and A.~N. Annaidh.
\newblock Poynting effect of brain matter in torsion.
\newblock {\em Soft Matter}, 15(25):5147--5153, 2019.
\newblock doi: \url{https://doi.org/10.1039/C9SM00131J}.

\bibitem{ball}
J.~M. Ball.
\newblock Convexity conditions and existence theorems in nonlinear elasticity.
\newblock {\em Arch. Rat. Mech. Anal.}, 63:337--403, 1976.

\bibitem{benesova}
B.~Benesova, J.~Forster, C.~Liu, and A.~Schl\"{o}merkemper.
\newblock Existence of weak solutions to an evolutionary model for
  magnetoelasticity.
\newblock {\em SIAM J. Math. Anal.}, 50(1):1200--1236, 2016.
\newblock doi: \url{https://doi.org/10.1137/17M1111486}.

\bibitem{collidreyer}
E.~Bonetti, P.~Colli, W.~Dreyer, G.~Gilardi, G.~Schimperna, and J.~Sprekels.
\newblock On a model for phase separation in binary alloys driven by mechanical
  effects.
\newblock {\em Phys. D}, 165(1-2):48--65, 2002.
\newblock doi: \url{https://doi.org/10.1016/S0167-2789(02)00373-1}.

\bibitem{bonetti1}
E.~Bonetti, P.~Colli, and M.~Fremond.
\newblock 2{D} motion with large deformations.
\newblock {\em Boll. Unione Mat. Ital.}, 7:19--44, 2014.
\newblock doi: \url{https://doi.org/10.1007/s40574-014-0002-0}.

\bibitem{bonetti2}
E.~Bonetti and M.~Fremond.
\newblock Motion of an incompressible solid with large deformations.
\newblock {\em C. R. Acad. Sci. Paris, Ser. I}, 356:345--350, 2018.
\newblock doi: \url{https://doi.org/10.1016/j.crma.2018.01.016}.

\bibitem{brezis}
H.~Brezis and P.~Mironescu.
\newblock {G}agliardo-{N}irenberg inequalities and non-inequalities: the full
  story.
\newblock {\em Ann. Inst. H. Poincar\'e - Anal. Non Lin\'eaire}, 35:1355--1376,
  2018.
\newblock doi: \url{https://doi.org/10.1016/j.anihpc.2017.11.007},.

\bibitem{lukacova2}
A.~Brunk, Y.~Lu, and M.~Luk{\'a}cov{\'a}-Medvidov{\'a}.
\newblock Existence, regularity and weak-strong uniqueness for the
  three-dimensional peterlin viscoelastic model.
\newblock 2021.
\newblock doi: \url{https://arxiv.org/abs/2102.02422}.

\bibitem{lukacova1}
A.~Brunk and M.~Luk{\'a}cov{\'a}-Medvidov{\'a}.
\newblock Global existence of weak solutions to viscoelastic phase separation:
  Part i regular case.
\newblock 2019.
\newblock doi: \url{https://arxiv.org/abs/1907.03480}.

\bibitem{neuro3}
N.~Bryniarska-Kubiak, A.~Kubiak, M.~Lekka, and A.~Basta-Kaim.
\newblock The emerging role of mechanical and topographical factors in the
  development and treatment of nervous system disorders: dark and light sides
  of the force.
\newblock {\em Pharmacological Reports}, 73:1626--1641, 2021.
\newblock doi: \url{https://doi.org/10.1007/s43440-021-00315-2}.

\bibitem{neuro5}
P.G. Carlier, B.~Marty, O.~Scheidegger, P.~Loureiro De~Sousa, P.Y. Baudin,
  E.~Snezhko, and D.~Vlodavets.
\newblock Skeletal muscle quantitative nuclear magnetic resonance imaging and
  spectroscopy as an outcome measure for clinical trials.
\newblock {\em J. Neuromuscul. Dis.}, 3:1--28, 2016.
\newblock \url{https://doi.org/10.3233/JND-160145},.

\bibitem{carman}
P.C. Carman.
\newblock Fluid flow through granular beds.
\newblock {\em Trans. Inst. Chem. Eng.}, 15:150166, 1937.
\newblock doi: \url{https://doi.org/10.1016/S0263-8762(97)80003-2}.

\bibitem{miranville}
M.~Carrive, A.~Miranville, and A.~Pi\'etrus.
\newblock The {C}ahn--{H}illiard equation for deformable elastic continua.
\newblock {\em Adv. Math. Sci. Appl.}, 10:539--569, 2000.

\bibitem{dileo}
C.~V {Di Leo}, E.~Rejovitzky, and L.~Anand.
\newblock A cahn–hilliard-type phase-field theory for species diffusion
  coupled with large elastic deformations: Application to phase-separating
  li-ion electrode materials.
\newblock {\em J. Mech. Phys. Solids}, 70:1–--29, 2014.
\newblock doi: \url{https://doi.org/10.1016/j.jmps.2014.05.001}.

\bibitem{frieboes}
H.~B. Frieboes, F.~Jin, Y.~L. Chuang, S.~M. Wise, J.~S. Lowengrub, and
  V.~Cristini.
\newblock Three--dimensional multispecies nonlinear tumor growth--{II}: Tumor
  invasion and angiogenesis.
\newblock {\em J Theor Biol.}, 264(4):1254--1278, 2010.
\newblock doi: \url{https://doi.org/10.1016/j.jtbi.2010.02.036}.

\bibitem{gagliardo}
E.~Gagliardo.
\newblock Ulteriori propriet\'a di alcune classi di funzioni in pi\'u
  variabili.
\newblock {\em Ricerche Mat.}, 8:24--51, 1959.

\bibitem{garckeel}
H.~Garcke.
\newblock On mathematical models for phase separation in elastically stressed
  solids.
\newblock {\em Habilitation Thesis, University of Bonn}, 2000.

\bibitem{garckeel2}
H.~Garcke.
\newblock On {C}ahn--{H}illiard systems with elasticity.
\newblock {\em Proc. Roy. Soc.}, 133A:307--331, 2003.

\bibitem{garcke4}
H.~Garcke, P.~Knopf, S.~Mitra, and A.~Schl\"{o}merkemper.
\newblock Strong well-posedness, stability and optimal control theory for a
  mathematical model for magneto-viscoelastic fluids.
\newblock 2021.
\newblock doi: \url{https://arxiv.org/abs/2108.03094}.

\bibitem{garcke2}
H.~Garcke, P.~Knopf, S.~Mitra, and A.~Schl\"{o}merkemper.
\newblock Strong well-posedness, stability and optimal control theory for a
  mathematical model for magneto-viscoelastic fluids.
\newblock 2021.
\newblock \url{https://arxiv.org/abs/2108.03094}.

\bibitem{garcke5}
H.~Garcke, B.~Kov\'{a}cs, and D.~Trautwein.
\newblock {V}iscoelastic {C}ahn-{H}illiard models for tumour growth.
\newblock 2022.
\newblock In preparation.

\bibitem{garckeel3}
H.~Garcke, K.F. Lam, and A.~Signori.
\newblock On a phase field model of {C}ahn--{H}illiard type for tumour growth
  with mechanical effects.
\newblock {\em Nonlinear Anal. Real World Appl.}, 57:103192, 2021.
\newblock doi: \url{https://doi.org/10.1016/j.nonrwa.2020.103192}.

\bibitem{grs}
M.~Grasselli, J.~E. Mu\v{n}oz~Rivera, and M.~Squassina.
\newblock Asymptotic behavior of a thermoviscoelastic plate with memory
  effects.
\newblock {\em Asymptotic Analysis 63}, 63(1), 2008.
\newblock doi: \url{https://doi.org/10.3233/ASY-2008-0928}.

\bibitem{gurtin1}
M.~E. Gurtin.
\newblock Generalized {G}inzburg--{L}andau and {C}ahn--{H}illiard equations
  based on a microforce balance.
\newblock {\em Phys. D}, 92(3--4):178–--192, 1996.
\newblock doi: \url{https://doi.org/10.1016/0167-2789(95)00173-5}.

\bibitem{gurtin2}
M.~E. Gurtin.
\newblock A gradient theory of single-crystal viscoplasticity that accounts for
  geometrically necessary dislocations.
\newblock {\em J. Mech. Phys. Solids}, 50(1):5--32, 2002.
\newblock doi: \url{https://doi.org/10.1016/S0022-5096(01)00104-1}.

\bibitem{heinemann}
C.~Heinemann, C.~Kraus, E.~Rocca, and R.~Rossi.
\newblock A temperature-dependent phase-field model for phase separation and
  damage.
\newblock {\em Arch. Ration. Mech. Anal.}, 225:177--247, 2017.
\newblock doi: \url{https://doi.org/10.1007/s00205-017-1102-7}.

\bibitem{holmes}
M.~H. Holmes and V.~C. Mow.
\newblock The nonlinear characteristics of soft gels and hydrated connective
  tissues in ultrafiltration.
\newblock {\em J Biomech.}, 23(11):1145--1156, 1990.
\newblock doi: \url{https://doi.org/10.1016/0021-9290(90)90007-p}.

\bibitem{neuro4}
R.~K. Jain, J.~D. Martin, and T.~Stylianopoulos.
\newblock The role of mechanical forces in tumor growth and therapy.
\newblock {\em Annual Review of Biomedical Engineering}, 16:321--346, 2014.
\newblock doi: \url{https://doi.org/10.1146/annurev-bioeng-071813-105259}.

\bibitem{kozeny}
J.~Kozeny.
\newblock {\"U}ber kapillare {L}eitung des {W}assers im {B}oden.
\newblock {\em Sitzungsber Akad Wiss Wien}, 136(2a):271306, 1927.

\bibitem{krejci}
P.~Krej\v{c}\'{i}, E.~Rocca, and J.~Sprekels.
\newblock Analysis of a tumor model as a multicomponent deformable porous
  medium.
\newblock 2021.
\newblock To appear on Interfaces and Free Boundaries. doi:
  \url{https://doi.org/10.4171/IFB/472}.

\bibitem{leoni}
G.~Leoni.
\newblock {\em A First Course in {S}obolev {S}paces: {S}econd {E}dition}.
\newblock Graduate Studies in Mathematics. 181. American Mathematical Society,
  2017.
\newblock ISBN 978-1-4704-2921-8.

\bibitem{mielke}
A.~Mielke and T.~Roub\'{\i}\v{c}ek.
\newblock Thermoviscoelasticity in {K}elvin-{V}oigt rheology at large strains.
\newblock {\em Arch. Ration. Mech. Anal.}, 238(1):1--45, 2020.
\newblock doi: \url{https://doi.org/10.1007/s00205-020-01537-z}.

\bibitem{neuro2}
C.~T. Motz, V.~Kabat, T.~Saxena, R.~V. Bellamkonda, and C.~Zhu.
\newblock Neuromechanobiology: An expanding field driven by the force of
  greater focus.
\newblock {\em Advanced Healthcare Materials}, 19(10):e2100102, 2021.
\newblock doi: \url{https://doi.org/10.1002/adhm.202100102}.

\bibitem{neuro1}
M.~C. Murphy, J.~3rd Huston, and R.~L. Ehman.
\newblock {M}{R} elastography of the brain and its application in neurological
  diseases.
\newblock {\em Neuroimage}, 187:176--183, 2019.
\newblock doi: \url{https://doi.org/10.1016/j.neuroimage.2017.10.008}.

\bibitem{nirenberg}
L.~Nirenberg.
\newblock On elliptic partial differential equations.
\newblock {\em Ann. Scuola Norm. Sup. Pisa}, 3(13):115--162, 1959.

\bibitem{ogden}
R.~W. Ogden.
\newblock Large deformation isotropic elasticity -- on the correlation of
  theory and experiment for incompressible rubberlike solids.
\newblock {\em Proc. R. Soc. Lond. A}, 326(1567):565--584, 1972.
\newblock doi: \url{https://doi.org/10.1098/rspa.1972.0026}.

\bibitem{rubicek}
T.~Roub{\'i}\v{c}ek.
\newblock Thermodynamics of viscoelastic solids, its eulerian formulation, and
  existence of weak solutions.
\newblock 2022.
\newblock doi: \url{https://arxiv.org/abs/2203.06080}.

\bibitem{temam}
R.~Temam.
\newblock {\em {N}avier-{S}tokes {E}quations: Theory and Numerical Analysis}.
\newblock {A}{M}{S} {C}helsea {P}ublishing. {V}olume 343., 1984.
\newblock ISBNs 978-0-8218-2737-6 (print); 978-1-4704-2994-2 (online).

\end{thebibliography}
\end{document}